\documentclass[letterpaper]{article}
\usepackage{authblk}
\usepackage{url}
\usepackage[margin=1in,footskip=0.25in]{geometry}
\usepackage{amssymb}
\usepackage[labelfont=bf]{caption}
\usepackage{graphicx}
\usepackage{amsmath}
\usepackage{algorithmic} 
\usepackage{algorithm}
\usepackage{float}
\usepackage{lipsum}
\usepackage{subcaption}
\usepackage{wrapfig}
\usepackage{cite}
\usepackage{color}
\usepackage{easylist}
\usepackage{lineno}
\usepackage{stmaryrd}
\usepackage{booktabs}
\usepackage{hyperref}
\usepackage{cleveref}
\usepackage{multirow}
\usepackage{hhline}
\usepackage{soul}
\usepackage{blindtext}
\usepackage{threeparttable}
\usepackage[bbgreekl]{mathbbol}

\usepackage{calc}
\usepackage[export]{adjustbox}

\usepackage{makecell}

\usepackage{tikz}
\newcommand{\Cross}{\mathbin{\tikz [x=1.4ex,y=1.4ex,line width=.2ex] \draw (0,0) -- (1,1) (0,1) -- (1,0);}}%

\usepackage{calc}
\usepackage[export]{adjustbox}

\def \grad{\nabla}

\def \p   {\partial}

\newcommand{\D}[2]{\frac{\p #1}{\p #2}}

\renewcommand{\vec}[1]{\bm{\mathrm{#1}}}

\def \CC{\mathbb{C}}

\def \FF{\mathbb{F}}

\def \PP{\mathbb{P}}

\def \F{\vec{F}}

\def \N{\vec{N}}

\def \U{\vec{U}}
\def \V{\vec{V}}
\def \X{\vec{X}}

\def \e{\vec{e}}
\def \f{\vec{f}}

\def \n{\vec{n}}

\def \f{\vec{f}}
\def \u{\vec{u}}

\def \x{\vec{x}}
\def \grad{\nabla}

\def \p{\partial}

\def \vchi{\vec{\chi}}

\def \u{\vec{u}}
\def \x{\vec{x}}

\def \Dx{{\mathrm d} \x}

\def \DX{{\mathrm d} \X}

\def \Dx{{\mathrm d} \x}

\newcommand{\euleriandx}{{\ensuremath{\Delta x}}}
\newcommand{\lagrangiandx}{{\ensuremath{\Delta X}}}
\def \mfac{M_\text{FAC}}
\def \efac{E_\text{FAC}}
\def \CL{C_\text{L}}
\def \CD{C_\text{D}}

\def \L2{L^2}

\def \grad{\nabla}

\def \p   {\partial}
\def \Mfac {M_{\text{fac}}}

\newcommand{\Ptwo}{\mathcal{P}^2}

\renewcommand{\vec}[1]{\ensuremath\boldsymbol{#1}}

\title{On the Lagrangian-Eulerian Coupling in the Immersed Finite Element/Difference Method}

\author[1,2,3]{Jae H.~Lee\footnote{Present address: Center for Drug Evaluation and Research, U.S. Food and Drug Administration, Silver Spring, MD, USA}} %\corref{cor1}}
\author[4,5,6,7]{Boyce E.~Griffith} %\corref{cor2}}
\affil[1]{Department of Mathematics, University of North Carolina, Chapel Hill, NC, USA}
\affil[2]{Department of Mechanical Engineering, Johns Hopkins University, Baltimore, MD, USA}
\affil[3]{Institute for Computational Medicine, Johns Hopkins University, Baltimore, MD, USA}
\affil[4]{Departments of Mathematics, Applied Physical Sciences, and Biomedical Engineering, University of North Carolina, Chapel Hill, NC, USA}
\affil[5]{Carolina Center for Interdisciplinary Applied Mathematics, University of North Carolina, Chapel Hill, NC, USA}
\affil[6]{Computational Medicine Program, University of North Carolina School of Medicine, Chapel Hill, NC, USA}
\affil[7]{McAllister Heart Institute, University of North Carolina School of Medicine, Chapel Hill, NC, USA}
\affil[ ]{\texttt{jaeholee@jhu.edu}}

\date{}

\begin{document}

\maketitle

\begin{abstract}
The immersed boundary (IB) method is a non-body conforming approach to fluid-structure interaction (FSI) that uses an Eulerian description of the momentum, viscosity, and incompressibility of a coupled fluid-structure system and a Lagrangian description of the deformations, stresses, and resultant forces of the immersed structure. 
Integral transforms with Dirac delta function kernels couple the Eulerian and Lagrangian variables, and in practice, discretizations of these integral transforms use regularized delta function kernels.
Many different kernel functions have been proposed, but prior numerical work investigating the impact of the choice of kernel function on the accuracy of the methodology has often been limited to simplified test cases or Stokes flow conditions that may not reflect the method's performance in applications, particularly at intermediate-to-high Reynolds numbers, or under different loading conditions.
This work systematically studies the effect of the choice of regularized delta function in several fluid-structure interaction benchmark tests using the immersed finite element/difference (IFED) method, which is an extension of the IB method that uses a finite element structural discretizations combined with a Cartesian grid finite difference method for the incompressible Navier-Stokes equations.
Whereas the conventional IB method spreads forces from the nodes of the structural mesh and interpolates velocities to those nodes, the IFED formulation evaluates the regularized delta function on a collection of interaction points that can be chosen to be denser than the nodes of the Lagrangian mesh.
This opens the possibility of using structural discretizations with wide node spacings that would produce gaps in the Eulerian force in nodally coupled schemes (e.g., if the node spacing is comparable to or broader than the support of the regularized delta function).
Earlier work with this methodology suggested that such coarse structural meshes can yield improved accuracy for shear-dominated cases and, further, found that accuracy improves when the structural mesh spacing is \textit{increased}.
However, these results were limited to simple test cases that did not include substantial pressure loading on the structure.
This study investigates the effect of varying the relative mesh widths of the Lagrangian and Eulerian discretizations in a broader range of tests.
Our results indicate that kernels satisfying a commonly imposed even--odd condition require higher resolution to achieve similar accuracy as kernels that do not satisfy this condition. 
We also find that narrower kernels are more robust, in the sense that they yield results that are less sensitive to relative changes in the Eulerian and Lagrangian mesh spacings, and that structural meshes that are substantially coarser than the Cartesian grid can yield high accuracy for shear-dominated cases but not for cases with large normal forces.
We verify our results in a large-scale FSI model of a bovine pericardial bioprosthetic heart valve in a pulse duplicator.
%Although this study is performed within the context of the IFED method, we argue that these results underscore the need to evaluate the impact of the choice of kernel function for other IB-type methods that use regularized delta functions to mediate fluid-structure interaction.
\end{abstract}


\newpage

%%%%%%%%%%%%%%%%%%%%%%%%%%%%%%%%%%%%%%%%%%%%%%%%%%%%%%%%%%%%%%%%%%%%%%%%%%%%%%%%%%%%%%%%%%%%%%%%%%%%%%%%%%%%%%%%%%%%%%%%%%%
\section{Introduction}
\label{sec:introduction}
The immersed boundary (IB) method~\cite{Peskin2002} is a non-body conforming approach to fluid-structure interaction (FSI) introduced by Peskin to model heart valves~\cite{Peskin1972, Peskin1977}. 
The IB approach to FSI uses an Eulerian description of the momentum, viscosity, and incompressibility of the coupled fluid-structure system, and it uses a Lagrangian description of the deformations, stresses, and resultant forces of the immersed structure.
In the continuous formulation, integral transforms with Dirac delta function\footnote{In fact, the singular Dirac delta function is not a function that is defined pointwise but instead is a generalized function or distribution. It is commonly referred to as the delta function within the IB literature, however, and we retain that usage herein.} kernels couple Eulerian and Lagrangian variables. 
When these equations are discretized, it is common to replace the singular delta function by a regularized delta function~\cite{Peskin2002}.
This coupling strategy eliminates the need for body-conforming discretizations and thereby facilitates models with very large structural deformations~\cite{Griffith2017, Griffith2020}.  
The IB method and its extensions have enabled simulation studies in a broad range of applications, including cardiac dynamics~\cite{Griffith2009, Griffith2012, Luo2012, Gao2014, Flamini2016, Chen2016, Gao2017, Hasan2017, Feng2018, Lee2020, LeeJTCVS}, platelet adhesion~\cite{Skorczewski2014}, esophageal transport~\cite{Kou2015, Kou2017, Kou2018}, heart development~\cite{Battista2018}, insect flight~\cite{Jones2015, Santhanakrishnan2018}, and undulatory swimming~\cite{Alben2013, Bhalla2014, Tytell2014, Bale2015, Hoover2017, Nangia2017}.

Despite the popularity of the IB method, most prior studies to examine the impact of the form of the regularized delta function on the accuracy of the method~\cite{Stockie1997, Roma1999, Peskin2002, Mori2008, Yang2009, Liu2012, Hosseini2016, Bao2016, Griffith2017, Bao2017, Saito2018, Heltai2020} have been limited to simplified test cases (e.g., two-dimensional Stokes problems) that may not reflect the method's performance in applications, particularly at intermediate-to-high Reynolds numbers, or under various loading conditions.
Peskin~\cite{Peskin2002} constructed a four-point regularized delta function that appears to be among the kernels most commonly used with the IB method. 
This function satisfies a certain set of properties, including an even--odd condition that is designed to avoid the well-known ``checkerboard'' instability that occurs with collocated discretizations of the incompressible Navier-Stokes equations. 
Roma et al.~\cite{Roma1999} introduced a three-point kernel function that satisfies the same properties as Peskin's four-point function except for the even--odd condition, which is not clearly needed for the staggered-grid fluid solver employed in that work. 
Stockie~\cite{Stockie1997} introduced a six-point IB kernel that yields higher-order accuracy than the three- and four-point IB kernels for problems with smooth solutions, albeit at expense of additional computational cost.
Yang et al.~\cite{Yang2009} developed smoothed $\mathcal{C}^2$ IB kernels that can suppress non-physical high-frequency force oscillations that can occur with the standard IB kernels.
Bao et al.~\cite{Bao2016, Bao2017} developed a new $\mathcal{C}^3$ six-point kernel that improves grid translational invariance and regularity compared to the standard three- and four-point kernels and the smoothed kernels of Yang et al.
Griffith and Luo~\cite{Griffith2017} used the benchmark problem of viscous flow past a cylinder to compare the standard three- and four-point kernels as well as the new six-point kernel by Bao et al.~\cite{Bao2016, Bao2017} and demonstrated that the choice of kernel function impacts the accuracy of the methodology. 
Mori~\cite{Mori2008} analyzed the convergence for the Stokes problem and showed that satisfying the even--odd condition improves the convergence properties of the method by eliminating high-frequency errors in the far field. 
Liu and Mori~\cite{Liu2012} extended the work of Mori to analyze convergence for elliptic problems and showed that the \textit{smoothing order}, which generalizes the even--odd condition, of a given delta function influences the convergence for the Stokes problem.
Hosseini et al.~\cite{Hosseini2016} analyzed the convergence of regularization for various PDEs with a singular source and demonstrated the substantial impact of regularization of the source term on the solutions to these problems.
Saito and Sugitani~\cite{Saito2018} studied the convergence of regularization error for a model Stokes problem in the context of finite element method.
Heltai and Lei~\cite{Heltai2020} provided a priori error estimates of regularization for elliptic problems compared to the non-regularized counterpart in the context of finite element formulations.
However, with the exception of the work by Griffith and Luo~\cite{Griffith2017}, none of these focus on tests in the intermediate-to-high Reynolds number regimes in which the IB method is commonly used in practice.
Here, we consider both the widely used IB kernels as well as B-spline kernels, which also are widely used delta function kernels~\cite{Hieber2008, He2021} but which, to our knowledge, have not been systematically compared against kernels that follow the construction approach of Peskin~\cite{Peskin2002, Roma1999, Bao2016, Bao2017} in the context of the IB method.

Herein we examine the impact of different choices of kernels on the dynamics using the immersed finite element/difference (IFED) method~\cite{Griffith2017, Griffith2020}, which is an extension of the IB method that uses finite element structural discretizations combined with a Cartesian grid finite difference method for the incompressible Navier-Stokes equations.
An important difference between the IFED method and conventional IB methods is that discrete IFED coupling operators use \textit{interaction points} that can be chosen to be distinct from the \textit{control points} that determine the configuration of the structure (e.g., the nodes of the Lagrangian mesh).
In this study, we follow the approach of Griffith and Luo~\cite{Griffith2017} and construct the interaction points via adaptively chosen Gaussian quadrature rules that distribute the interaction points in the interiors of the structural elements.
In contrast, the conventional IB method spreads forces from the nodes of the structural mesh and interpolates velocities to those nodes~\cite{Peskin2002}.
In nodally coupled IB methods, catastrophic leakage flows \textit{through the structure} can occur if the node spacing is comparable to or larger than the support of the regularized delta function because in such cases, there will be gaps in the Eulerian structural force density. 
(This issue is distinct from the question of the fundamental volume conservation of the IB method, which has been the subject of numerous studies, including the work of Peskin and Prinz~\cite{Peskin1993} along with more recent work for both immersed boundary and immersed finite element-type methods~\cite{Newren2007, Stockie2009, Wang2010, GriffithPPM2012, Bao2017, Vadala-Roth2020}.)
The IFED approach of using distinct collections of control and interaction points opens the possibility of using structural discretizations with wide node spacings while maintaining a contiguous Eulerian structural force density.
However, prior studies on the impact of the relative node spacing on the impact of the IFED method have been limited and, in particular, considered only cases with negligible normal forces along the fluid-structure interface~\cite{Griffith2017}.
At least in those tests, however, it was found that the accuracy of the method actually \textit{increases} with increasing Lagrangian mesh spacing.
In this study, we systematically investigate the impact of the relative spacings of the Lagrangian and Eulerian discretizations for a broader range of test problems in the intermediate-to-high Reynolds number regimes ranging from 70 to 15000.
The results in this study are concordant with earlier work for shear-dominated cases in that narrower kernels are more robust and that a broad range of relatively coarse structural meshes can be used, but here we also identify that the structural mesh spacing must be comparable to or finer than the background Cartesian grid for cases involving large pressure loads.
Our results also indicate that kernels satisfying a commonly imposed even--odd condition require higher resolution to achieve similar accuracy as kernels that do not satisfy this condition.
We then apply and verify our key findings in a large-scale FSI model of bovine pericaridal bioprosthetic heart valve (BHV) in a pulse duplicator~\cite{Lee2020, LeeJTCVS}.
Although these investigations are all done within the context of the IFED version of the IB method, the large effect of the choice of kernel function on the results suggests the need for similar studies for other IB-type methods that use regularized delta functions to mediate fluid-structure interaction.

\section{Methods}
\label{sec:methods}

This section describes the continuous formulation of the IFED method and the numerical discretization and implementation of the method.
We also define the key factors that impact the interaction between the Lagrangian mesh and the Eulerian grid such as different types of regularized delta functions, as well as the Lagrangian mesh spacing.

\subsection{Immersed finite element/difference method}
\label{subsec:IBM}

The continuous IFED formulation considers fluid-structure system occupying a fixed three-dimensional Eulerian computational domain $\Omega$ that is partitioned into time-dependent fluid ($\Omega^\text{f}_t$) and solid ($\Omega^\text{s}_t$) subdomains, so that $\Omega = \Omega^\text{f}_t \cup \Omega^\text{s}_t$.
Here, $\x = (x_1, x_2, x_3) \in \Omega$ are physical coordinates, $\X = (X_1, X_2, X_3) \in \Omega^\text{s}_0$ are reference coordinates attached to the structure, $\N(\X)$ is the outward unit normal to $\partial\Omega^\text{s}_0$ at material position $\X$, and $\vchi(\X,t) \in \Omega^\text{s}_t$ is the physical position of material point $\X$ at time $t$.
The dynamics of the coupled system are described by
 \begin{align}
    \label{eq:momentum}  \rho\frac{{\mathrm D}\u}{{\mathrm D}t}(\x,t) &= - \grad p(\x,t) + \mu \grad^2 \u(\x,t) + \f(\x,t), \\
    \label{eq:continuity}     \grad \cdot \u(\x,t) &= 0, \\
    \label{eq:fsiconstraint} \f(\x,t) &= \int_{\Omega^s_0}\F(\X,t)\,\delta(\x - \vchi(\X,t)) \, \DX,  \\
    \label{eq:noslip}           \D{\vchi}{t}(\X,t) &= \U(\X,t) = \int_\Omega \u(\x,t) \, \delta(\x - \vchi(\X,t)) \, \Dx = \u(\vchi(\X,t),t),
  \end{align}
 in which $\frac{{\mathrm D}}{{\mathrm D}t} = \frac{\p}{\p t} + \u \cdot \nabla$ is the material derivative, $\u(\x,t)$ and $p(\x,t)$ are the Eulerian velocity and pressure fields, $\f(\x,t)$ is the Eulerian structural force density, $\F(\X,t)$ is the Lagrangian force density, $\U(\X,t)$ is the Lagrangian velocity of the immersed structure, and $\delta(\x) = \prod_{i=1}^{3}\delta(x_{i})$ is the three-dimensional Dirac delta function. 
For simplicity, we assume a uniform mass density $\rho$ and viscosity $\mu$.
Eq.~\eqref{eq:fsiconstraint} implies that the Eulerian and Lagrangian force densities are equivalent as densities, and Eq.~\eqref{eq:noslip} implies that the no-slip condition is satisfied along the fluid-structure interface.
Note that because $\D{\vchi}{t}(\X,t) = \u(\vchi(\X,t),t)$ and $\grad \cdot \u(\x,t) = 0$, the immersed structure is incompressible~\cite{Vadala-Roth2020}.

In our numerical tests, we consider both rigid and elastic immersed structures.
For stationary structures considered in our examples, $\F(\X,t)$ in Eq.~\eqref{eq:fsiconstraint} is a Lagrange multiplier for the constraint $\U(\X,t)\equiv \vec{0}$.
We use a penalty formulation~\cite{Goldstein1993} that yields an approximate Lagrange multiplier force,
\begin{equation}
\F(\X,t) = \kappa(\X - \vchi(\X,t)) - \eta\U(\X,t),
\label{eq:tether_force}
\end{equation}
in which $\kappa$ is a stiffness penalty parameter and $\eta$ is a body damping penalty parameter.
Note that as $\kappa \to \infty$, $\vchi(\X,t) \to \X$ and $\D{\vchi}{t}(\X,t) \to {\bf 0}$.
We include a damping term in the penalty force to reduce spurious oscillations that can occur in practice for finite $\kappa$.

We also consider immersed elastic structures in Sections~\ref{subsec:turek-hron},~\ref{subsec:pressurized_band}, and~\ref{subsec:BHV}. 
In the simplest version of this methodology, the immersed structure is modeled as a viscoelastic solid, in which the viscous stresses in the solid are typically small compared to elastic stresses~\cite{Zhang2004, Zhang2007, Wang2012, Griffith2017}.
In our IFED formulation, the elastic response is that of a hyperelastic material, for which the first Piola--Kirchhoff stress $\PP$ of the immersed structure is related to a strain-energy functional $\Psi(\FF)$ via $\PP = \D{\Psi}{\FF},$ in which $\FF = \p\vchi/\p\X$ is the deformation gradient tensor.
The resultant structural force $\F(\X,t)$ generated by deformations of the elastic structure is determined in a weak sense by satisfying
\begin{align}
\int_{\Omega^\text{s}_0}\F(\X,t)\cdot\V(\X) \, \DX 
		&= -\int_{\Omega^\text{s}_0}  \PP(\X,t) : \grad_{\X}\V(\X) \, \DX
\end{align}
for all smooth $\V(\X)$~\cite{Boffi2008,Griffith2017}.
This is the so-called \textit{unified} weak formulation~\cite{Griffith2017}, which incorporates both internal and transmission forces~\cite{Boffi2008}.
Consequently, $\F$ can be a generalized function or distribution, with force concentrations along the fluid-solid interface that are singular like a one-dimensional delta function.
As in earlier work using finite element-based structural discretizations with the IB framework, these singularities are effectively regularized by projecting them onto the finite element shape functions~\cite{Boffi2008, Zhang2004, Devendran2012, Griffith2017}.
(Griffith and Luo~\cite{Griffith2017} also considered a \textit{partitioned} formulation that separately approximated the (regular) interior force density and the (singular) transmission force density.
In practice, we have not found cases in which that approach yields substantially improved accuracy, but we do find that it yields poorer stability in many cases. 
Consequently, we focus on the unified formulation in this work.)
We also use this approach as a penalty formulation to model rigid structures by treating the structure as an elastic material with a large stiffness parameter.

\subsection{Eulerian and Lagrangian discretizations}
\label{subsec:eul_discretization}
\begin{figure}[t!!]
 \centering
 \includegraphics[scale = 0.25]{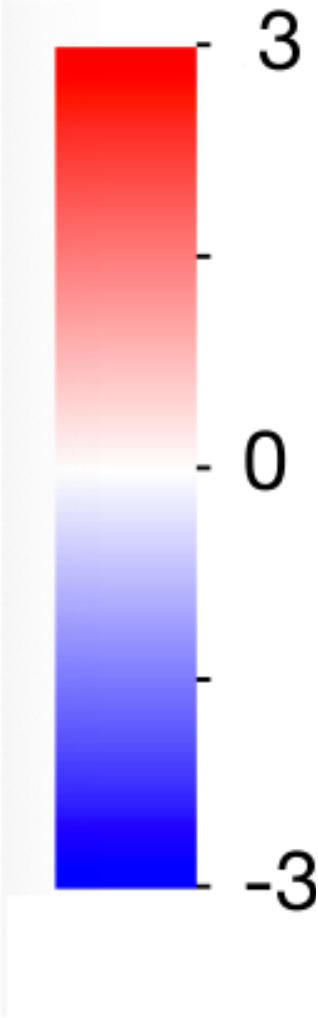}
 \includegraphics[scale = 0.25]{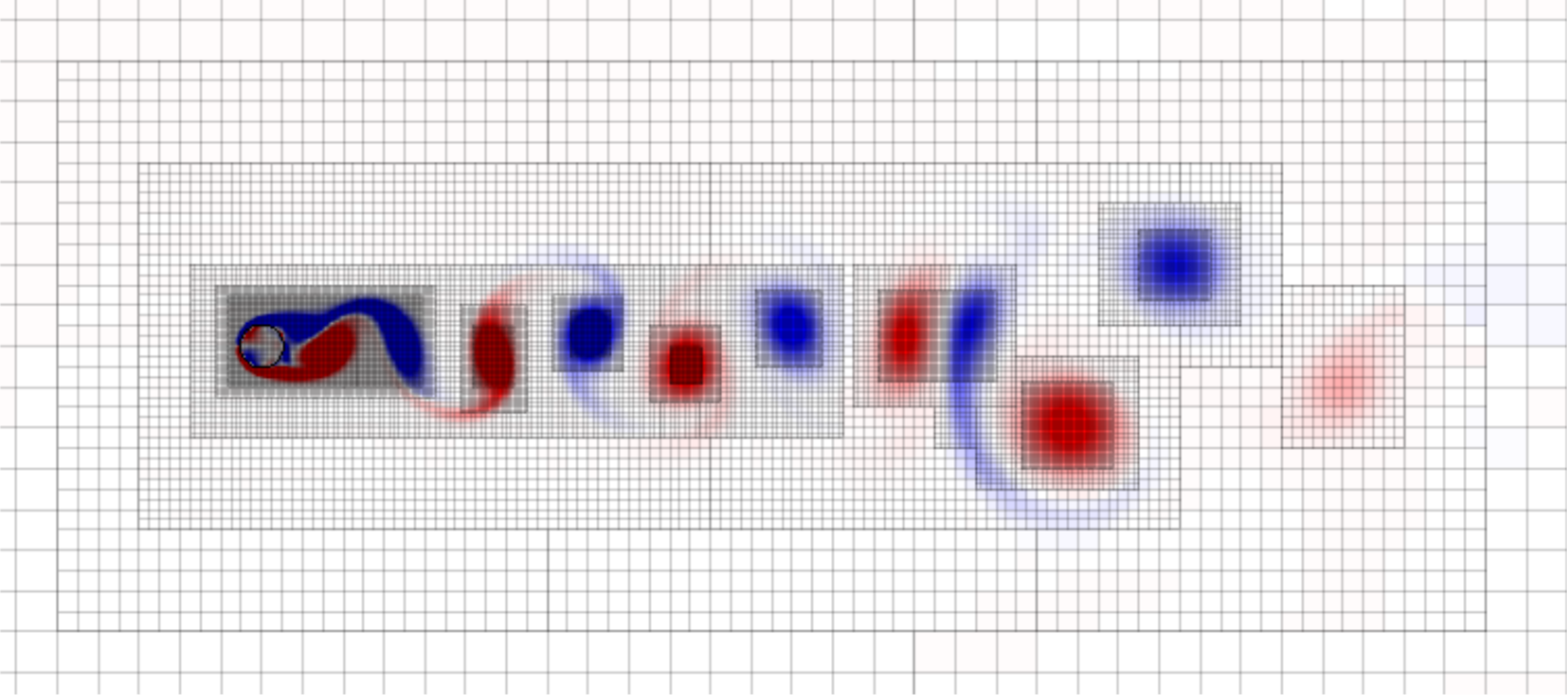}
 \caption{Vortices shed from a stationary circular cylinder at $Re = 200$.
 The computational domain is described by block-structured adaptively refined Cartesian grid that dynamically tracks vortices shed from the immersed structure.}
 \label{fig:AMR}
\end{figure}

The Eulerian variables are solved on the computational domain $\Omega$, which includes both the solid and fluid subregions, and this domain is described using a block-structured locally refined Cartesian grids consisting of nested levels of Cartesian grid patches~\cite{Griffith2012}.
This allows high spatial resolution to be deployed dynamically near fluid-structure interfaces and near flow features that are identified by feature detection criteria (e.g., local magnitude of the vorticity) for enhanced spatial resolution. 
Figure~\ref{fig:AMR} provides an example of the adaptive mesh refinement in the test case of flow past a cylinder.
We use a second-order accurate staggered-grid discretization~\cite{GriffithPPM2009, GriffithPPM2012} of the incompressible Navier-Stokes that includes a version of the piecewise parabolic method (PPM)~\cite{Colella1984} to approximate the convective term.

The Lagrangian variables are solved on the immersed structure, which is discretized with $\mathcal{C}^0$ finite elements as described in Griffith and Luo~\cite{Griffith2017}. 
Briefly, we construct a triangulation, $\mathcal{T}^h$, with $m$ nodes, in which we define the $3m$-dimensional vector-valued approximation space as $X^h\subset H^1(\mathcal{T}^h)^3$.
We then define $\{\phi_\ell\}$ to be the standard nodally interpolating finite element basis of $X^h$.
We track deformation, velocity, and force at the nodes and use the same shape functions for each component, which can be written as
\begin{align}
	\vchi_h(\X,t) &=\sum_{\ell=1}^{m}\vchi_\ell(t)\,\phi_\ell(\X),\\
	\U_h(\X,t) &= \sum_{\ell=1}^{m} \U_\ell(t)\,\phi_\ell(\X), \text{ and}\\
	\F_h(\X,t) &= \sum_{\ell=1}^{m} \F_\ell(t)\,\phi_\ell(\X).
\end{align}
For the rest of this discussion, we drop the subscript ``$h$'' from the numerical approximations to the Lagrangian variables to simplify notation.

\subsection{Lagrangian-Eulerian coupling}
\label{subsec:lag_eul_coupling}
As briefly described in Section~\ref{sec:introduction}, the coupling between Eulerian and Lagrangian variables is mediated by integral transforms with delta function kernels as shown in Eq.~\eqref{eq:fsiconstraint} and~\eqref{eq:noslip}.
To approximate $\f=(f_1,f_2,f_3)$ in Eq.~\eqref{eq:fsiconstraint} on the Cartesian grid, we construct a Gaussian quadrature rule with $N^e$ quadrature (or interaction) points $\X_Q^e \in K^e$ and weights $w_Q^e$, $Q = 1, \ldots, N^e$ for each element $K^e \in \mathcal{T}^h$. 
Then $f_1$, $f_2$, and $f_3$ on the faces of the Cartesian grid cells are computed as~\cite{Griffith2017}
\begin{align}
(f_1)_{i-\frac{1}{2},j,k} &= \sum_{K^e\in\mathcal{T}^h}\sum_{Q=1}^{N^e}F_1(\X_Q^e, t)\,\delta_h(\x_{i-\frac{1}{2},j,k}-\vchi(\X_Q^e,t))w_Q^e,\\
(f_2)_{i,j-\frac{1}{2},k} &= \sum_{K^e\in\mathcal{T}^h}\sum_{Q=1}^{N^e}F_2(\X_Q^e, t)\,\delta_h(\x_{i,j-\frac{1}{2},k}-\vchi(\X_Q^e,t))w_Q^e,\\
(f_3)_{i,j,k-\frac{1}{2}} &= \sum_{K^e\in\mathcal{T}^h}\sum_{Q=1}^{N^e}F_3(\X_Q^e, t)\,\delta_h(\x_{i,j,k-\frac{1}{2}}-\vchi(\X_Q^e,t))w_Q^e,
\end{align}
in which $\F(\X,t)=(F_1(\X,t),F_2(\X,t),F_3(\X,t))$ are the Lagrangian force densities and $\delta_h(\x)$ is a regularized delta function. 
We use the compact notation
\begin{equation}
	\f(\x,t) = \boldsymbol{\mathcal{S}}[\vchi(\cdot,t)]\,\F(\X,t),
\end{equation}
in which $\boldsymbol{\mathcal{S}}[\vchi(\cdot,t)]$ is the force-prolongation operator. 
Similarly, the velocity of the structure, $\D{\vchi}{t}(\X,t)$ in Eq.~\eqref{eq:noslip}, can be approximated by using the Cartesian grid velocity $\u(\x,t)$,
\begin{equation}
	 \D{\vchi}{t}(\X,t) = \boldsymbol{\mathcal{J}}[\vchi(\cdot,t)]\,\u(\x,t),
\end{equation}
in which $\boldsymbol{\mathcal{J}}[\vchi(\cdot,t)]$ is the velocity-restriction operator that is constructed to satisfy the adjoint condition, $\boldsymbol{\mathcal{J}} = \boldsymbol{\mathcal{S}}^*$~\cite{Griffith2017}.
It is clear that the coupling operators $\boldsymbol{\mathcal{S}}$ and $\boldsymbol{\mathcal{J}}$  depend on the spatial discretization and the choice of regularized delta function kernel, $\delta_h$.

\subsection{Regularized delta functions}
\label{subsec:numerical}

In our computations, we use a regularized delta function $\delta_h(\x)$ in our discrete approximations to the integral transforms in Eq.~\eqref{eq:fsiconstraint} and~\eqref{eq:noslip}. 
\begin{figure}[t!!]
	\begin{center}
		\begin{subfigure}[t]{0.19\textwidth}
			\centering
			\includegraphics[scale = 0.1575]{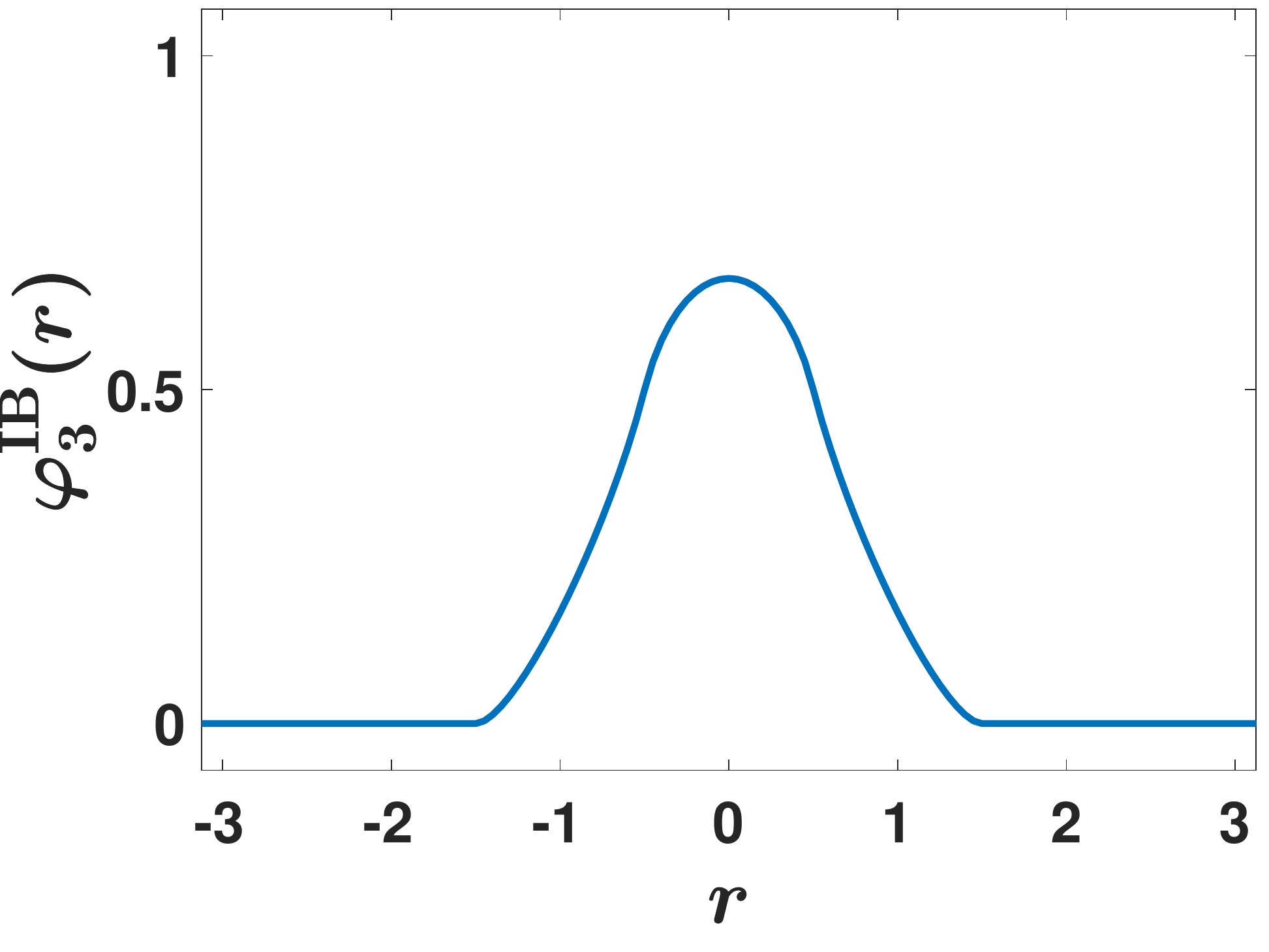}
			\caption{\parbox[t]{1.6cm}{\centering three-point\\IB}} \label{fig:ib_3}
		\end{subfigure}
		\begin{subfigure}[t]{0.19\textwidth}
			\centering
			\includegraphics[scale = 0.1575]{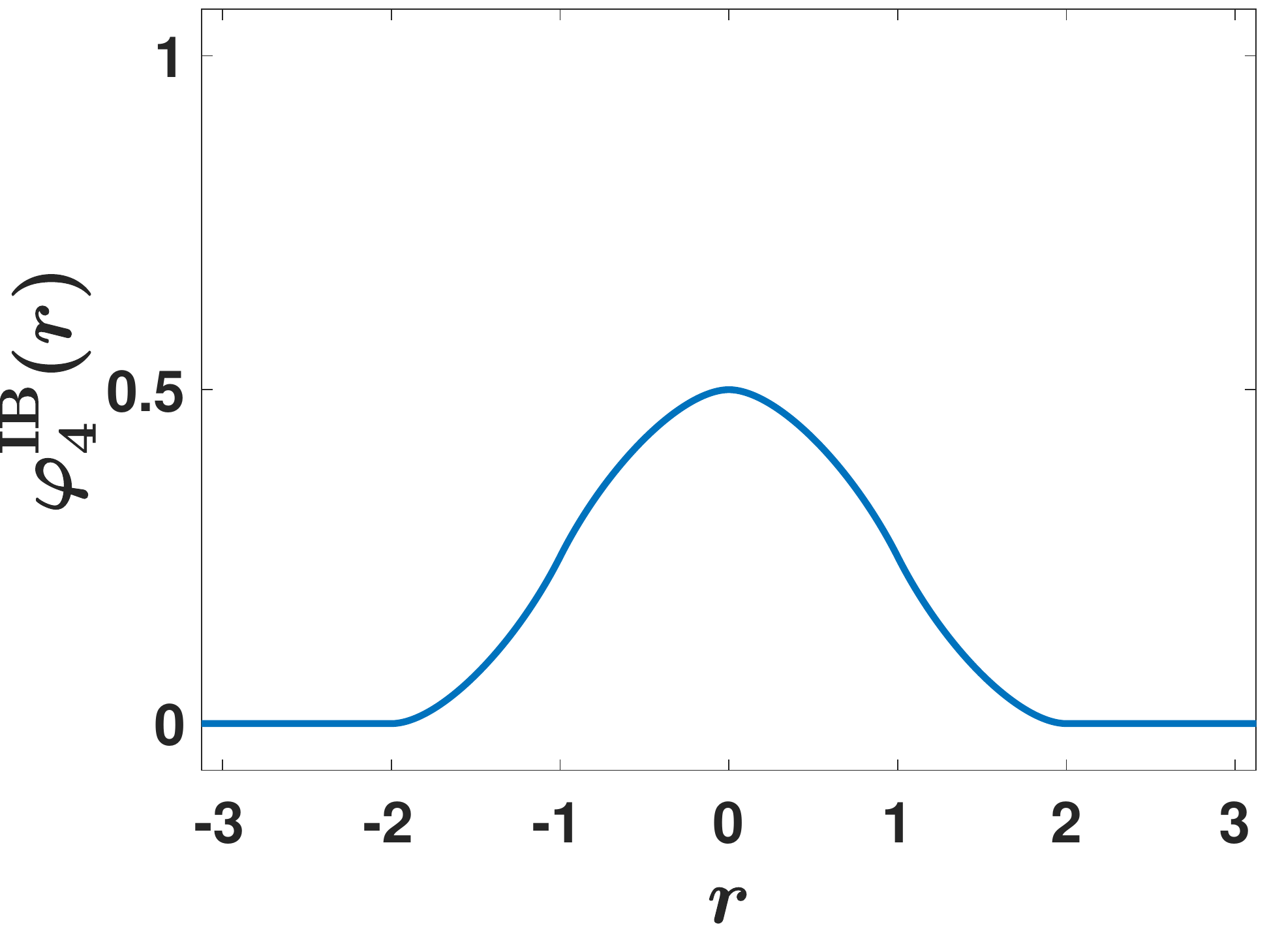}
			\caption{\parbox[t]{1.6cm}{\centering four-point\\IB}} \label{fig:ib_4}
		\end{subfigure} 
		\begin{subfigure}[t]{0.19\textwidth}
			\centering
			\includegraphics[scale = 0.1575]{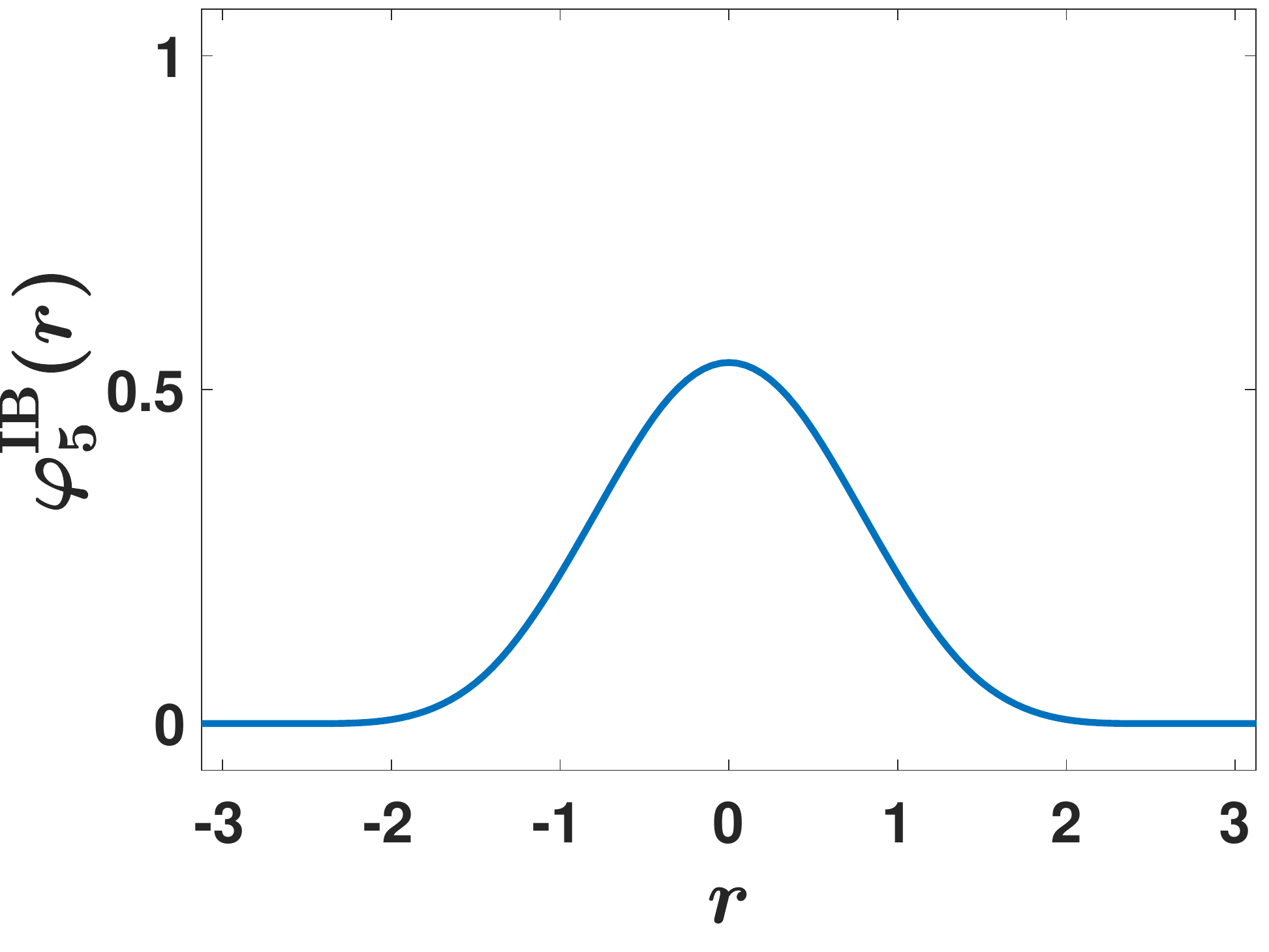}
			\caption{\parbox[t]{1.6cm}{\centering five-point\\IB}} \label{fig:ib_5}
		\end{subfigure}
		\begin{subfigure}[t]{0.19\textwidth}
			\centering
			\includegraphics[scale = 0.1575]{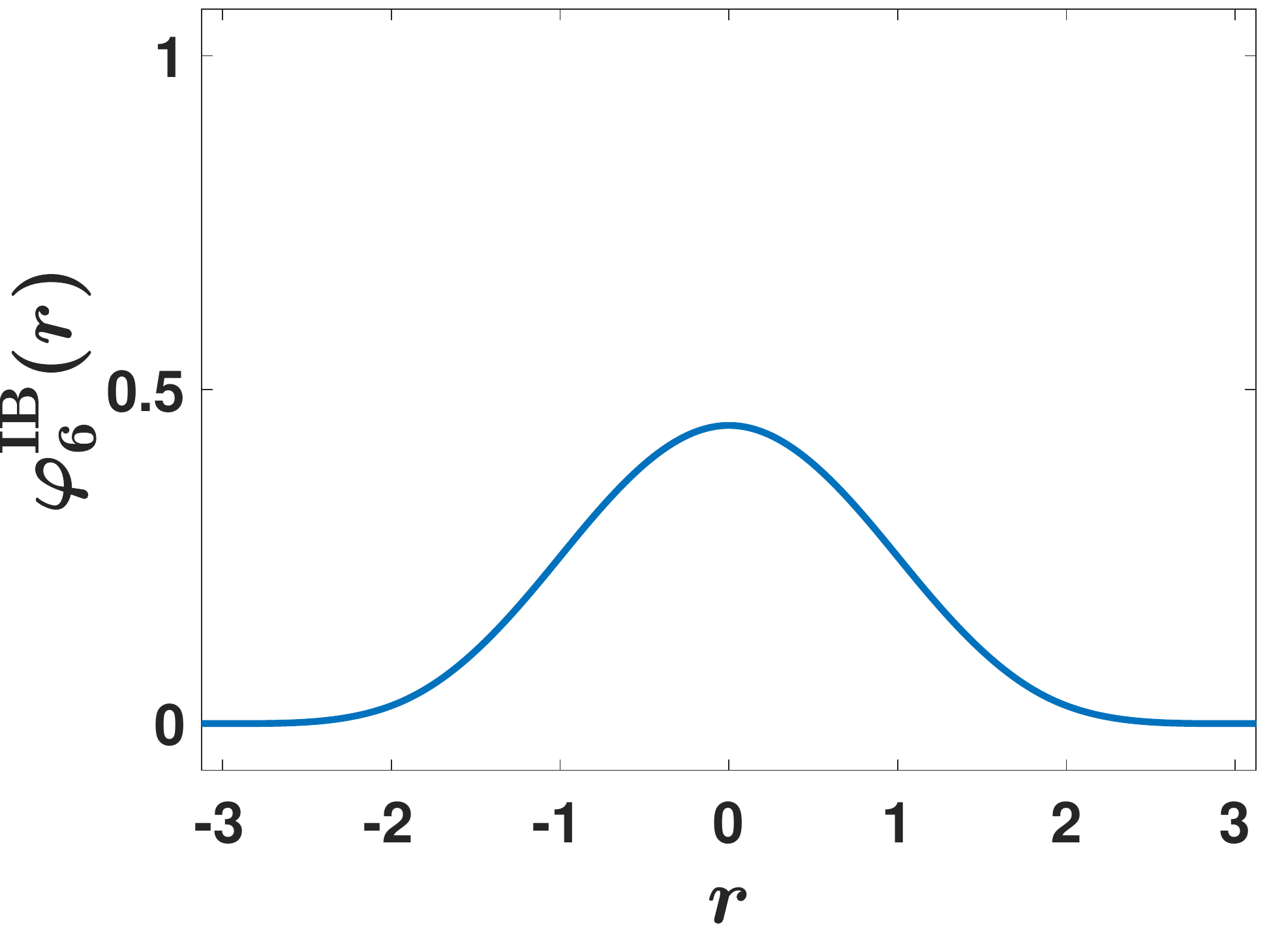}
			\caption{\parbox[t]{1.6cm}{\centering six-point\\IB}} \label{fig:ib_6}
		\end{subfigure}\\
		\vspace{0.1in}
		\begin{subfigure}[t]{0.19\textwidth}
			\centering
			\includegraphics[scale = 0.1575]{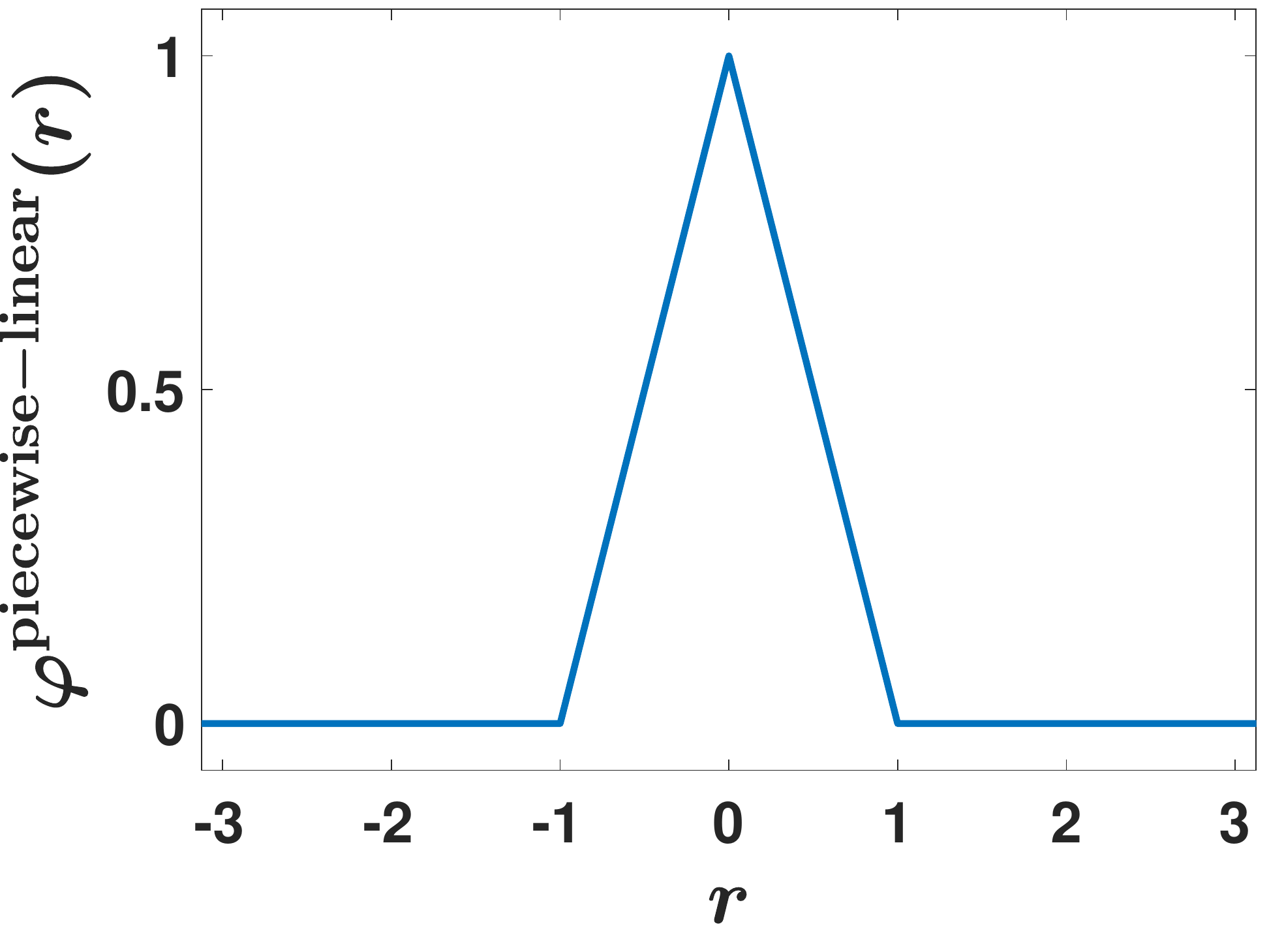}
			\caption{\parbox[t]{1.5cm}{\centering two-point\\B-spline}} \label{fig:piecewise_linear}
		\end{subfigure}
		\begin{subfigure}[t]{0.19\textwidth}
			\centering
			\includegraphics[scale = 0.1575]{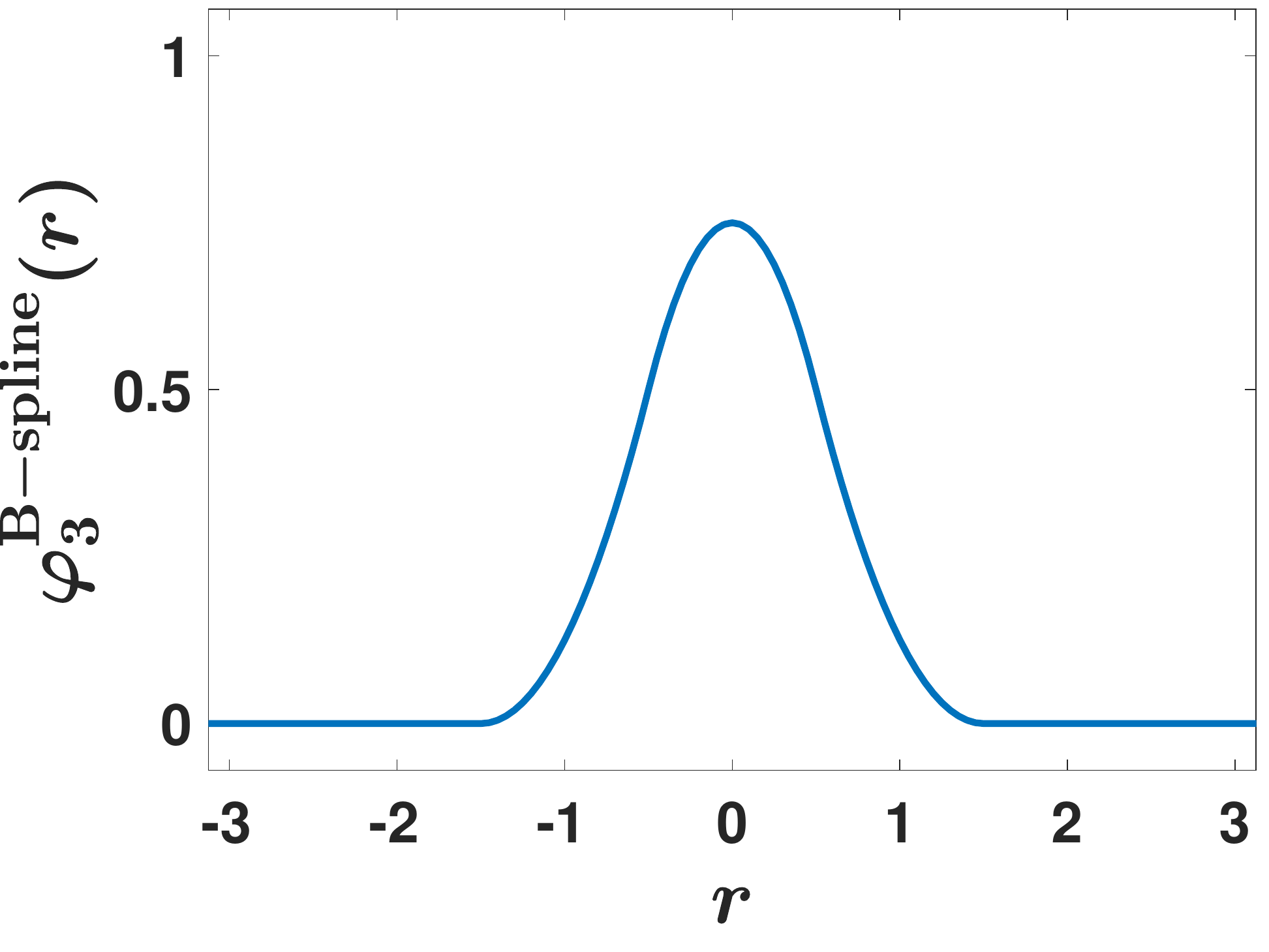}
			\caption{\parbox[t]{1.6cm}{\centering three-point\\B-spline}} \label{fig:bspline_3}
		\end{subfigure} 
		\begin{subfigure}[t]{0.19\textwidth}
			\centering
			\includegraphics[scale = 0.1575]{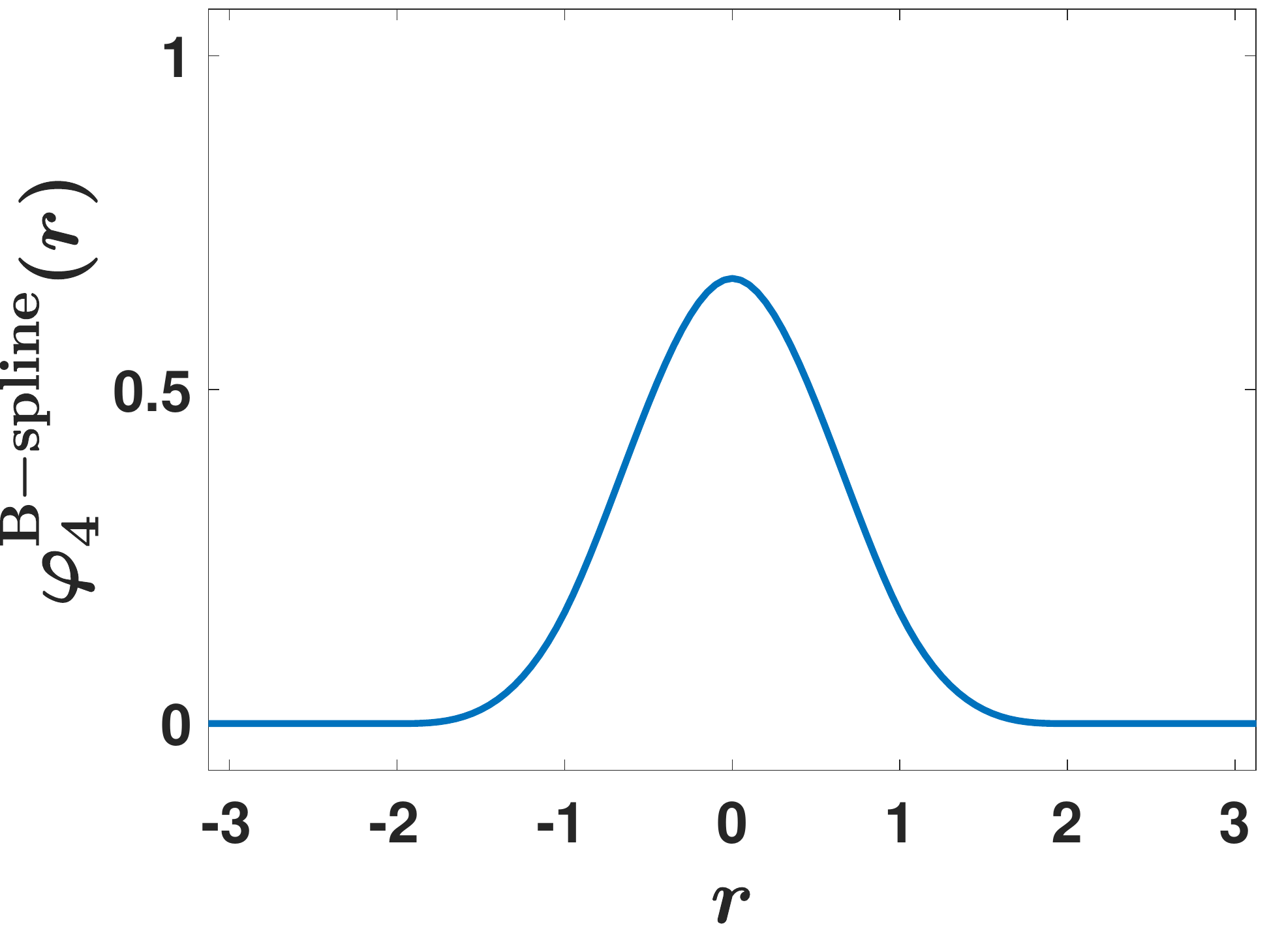}
			\caption{\parbox[t]{1.6cm}{\centering four-point\\B-spline}} \label{fig:bspline_4}
		\end{subfigure}
		\begin{subfigure}[t]{0.19\textwidth}
			\centering
			\includegraphics[scale = 0.1575]{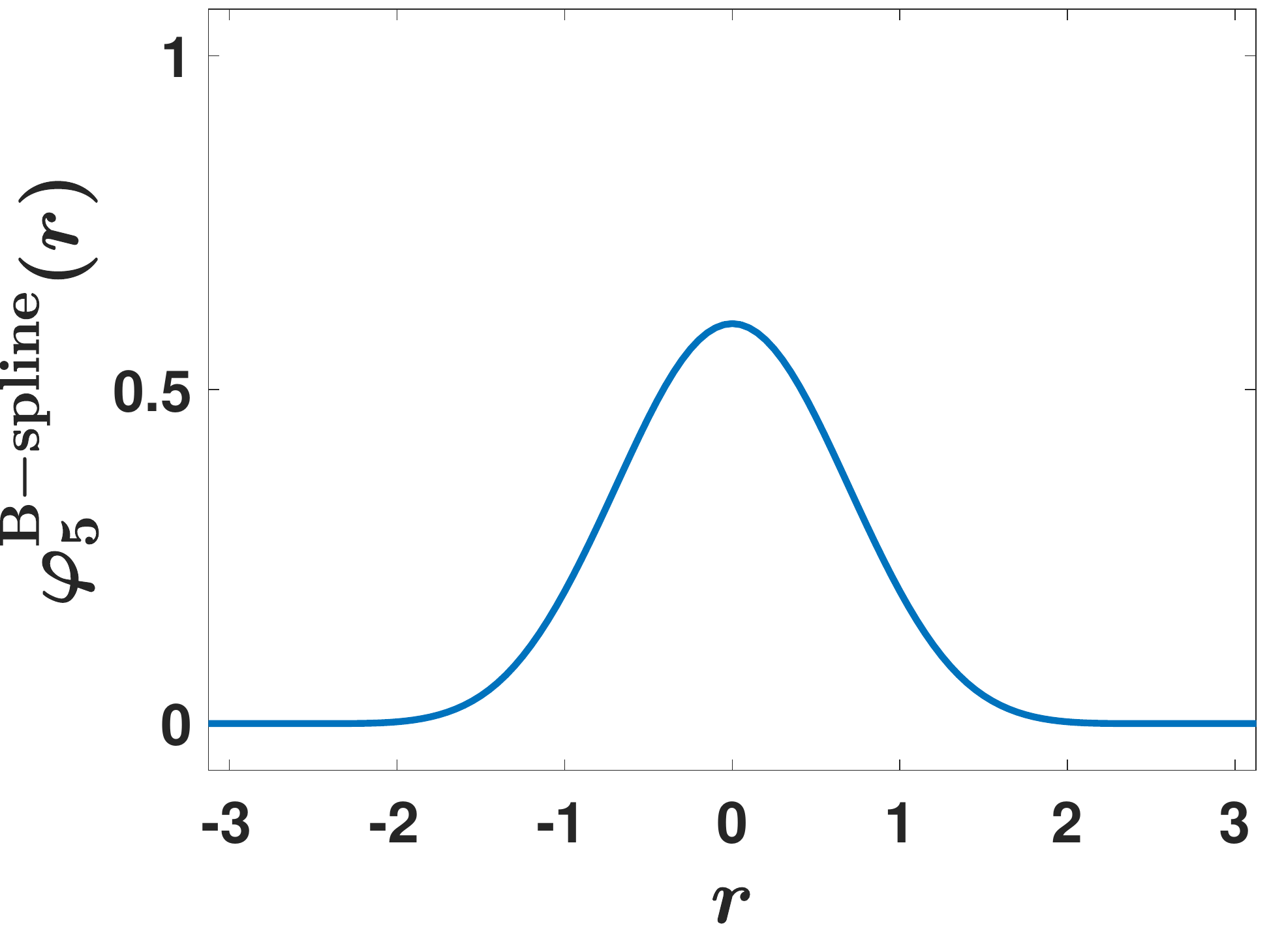}
			\caption{\parbox[t]{1.6cm}{\centering five-point\\B-spline}} \label{fig:bspline_5}
		\end{subfigure}
		\begin{subfigure}[t]{0.19\textwidth}
			\centering
			\includegraphics[scale = 0.1575]{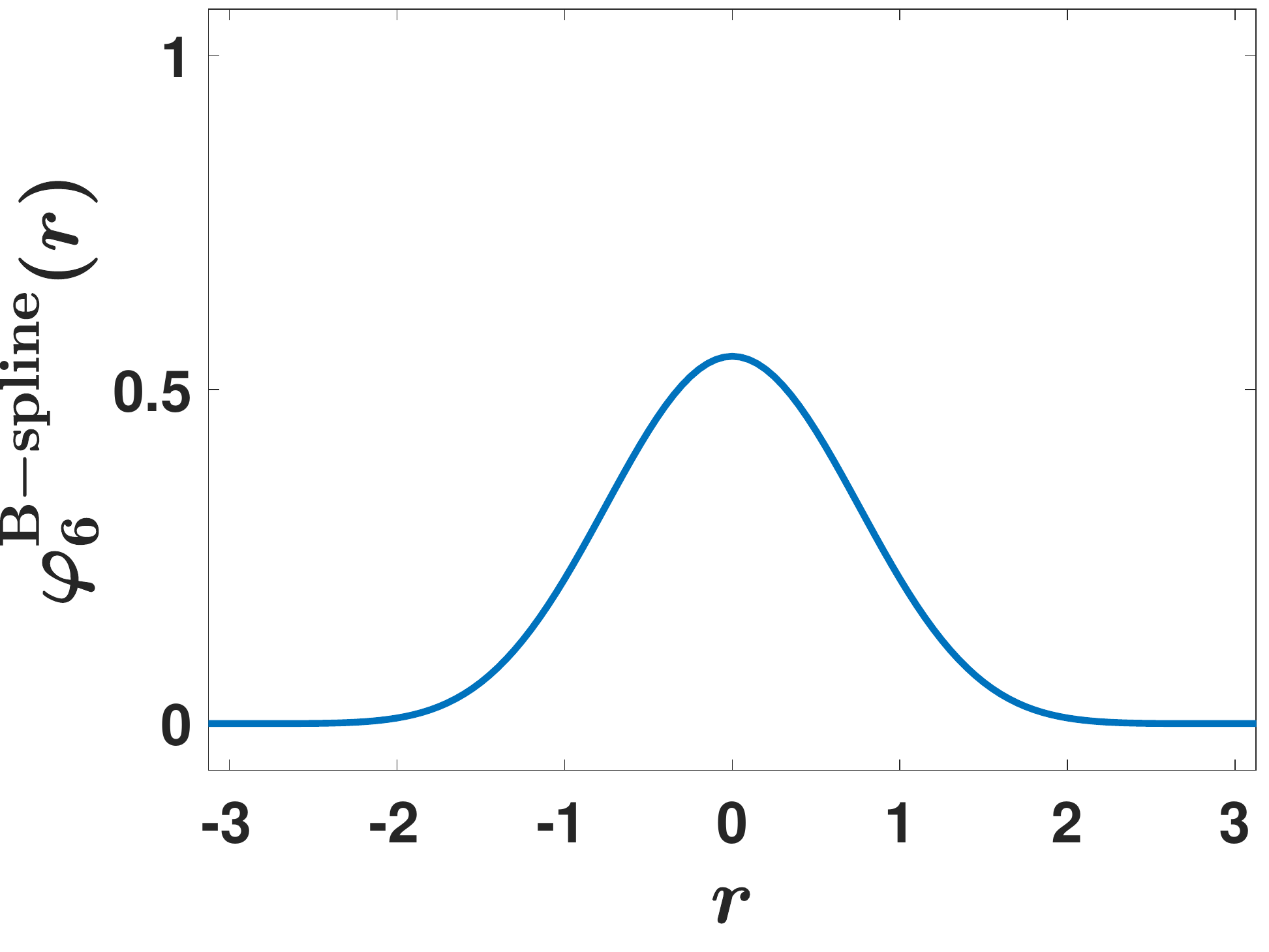}
			\caption{\parbox[t]{1.6cm}{\centering six-point\\B-spline}} \label{fig:bspline_6}
		\end{subfigure}
		\caption{Selected choices of regularized delta functions. One family of kernel functions is determined by imposing some or all of the conditions described by Peskin~\cite{Peskin2002}, and they present different properties depending on which of the moment conditions are satisfied. We can also consider B-spline kernels that are constructed by recursive convolution against piecewise-constant kernels.}
		\label{fig:delta_kernels}		
	\end{center}
\end{figure}
Following Peskin \cite{Peskin2002}, we construct the three-dimensional regularized delta function as the tensor product of one-dimensional delta functions, $\delta_h(\x) = \prod_{i=1}^3 \delta_h(x_i)$, and the one-dimensional regularized delta function is defined in terms of a basic kernel function via $\delta_h(x) = \frac{1}{h} \varphi\left(\frac{x}{h}\right)$. 
Note that $\varphi$ is different from $\phi$ used earlier to denote finite element basis functions.
Here, $\varphi\left(r\right)$ is continuous for all $r$ and zero outside of the radius of support.
Figure~\ref{fig:delta_kernels} shows different regularized delta functions considered in this study. 
One-dimensional kernel functions introduced by Peskin impose some or all of the following conditions~\cite{Peskin2002}:
\begin{align}
\label{eq:zeroth_moment} \text{zeroth moment: }& \sum_j \varphi(r-j) = 1;\\
\label{eq:even_odd} \text{ even--odd: }& \sum_{j \text{ even}} \varphi(r-j) = \sum_{j \text{ odd}} \varphi(r-j) = \frac{1}{2};\\
\label{eq:first_moment} \text{first moment: }& \sum_j(r-j)\, \varphi(r-j) = 0;\\
\label{eq:second_moment} \text{second moment: }& \sum_j(r-j)^2\, \varphi(r-j) = K, \text{ for some constant } K.
\end{align}
The zeroth moment condition implies total forces are equivalent in discretized Lagrangian or Eulerian form when $\delta_h$ is used for force spreading~\cite{Peskin2002}. 
The even--odd condition is designed to avoid the ``checkerboard'' instability in a collocated-grid fluid solver and thereby to suppress spurious high-frequency modes~\cite{Roma1999, Peskin2002, GriffithPPM2012, Bao2016, Bao2017}. 
Note that the even--odd condition implies the zeroth moment condition. 
The first moment condition implies the conservation of total torque. 
Along with the zeroth moment condition, it guarantees second-order accuracy in interpolating smooth functions~\cite{Peskin2002}. 
If a kernel function satisfies Eq~\eqref{eq:second_moment} with $K = 0$, then the second moment condition implies that the kernel achieves higher order accuracy in interpolating smooth fields. 
It is also possible to use the higher-order moment condition with $K \ne 0$, which can be used to impose higher continuity order on the kernel function~\cite{Bao2016}. 
Peskin also postulated a sum-of-squares condition,
\begin{equation}
\label{eq:sum-of-squares} \sum_j(\varphi(r-j))^2 = C, \text{ for some constant } C,
\end{equation}
which is a weak version of a grid translational invariance property~\cite{Peskin2002}.

At present, the kernel functions most commonly used with the IB method appear to be what we refer to as the \textit{IB kernels}, which satisfy some or all of the properties proposed by Peskin~\cite{Peskin2002} (Figures~\ref{fig:ib_3} and~\ref{fig:ib_4}). 
The three-point IB kernel is constructed by satisfying the zeroth and first moment conditions as well as the sum-of-squares condition, but not the even--odd condition~\cite{Peskin2002}.
The five-point IB kernel satisfies the same conditions as the three-point function along with second and third moment conditions with $K \neq 0$ chosen to yield higher continuity order~\cite{Bao2016, Bao2017}.
The four-point IB kernel is constructed by satisfying the even--odd condition (which also implies the zeroth moment condition) and first moment conditions as well as the sum-of-squares condition~\cite{Peskin2002}.
The six-point IB kernel satisfies the same conditions as the four-point function along with second and third moment conditions with $K \neq 0$ chosen to yield higher continuity order~\cite{Bao2016, Bao2017}.
We emphasize that the three-~and four-point IB kernels satisfy the same properties, except that the three-point kernel does not satisfy the even-odd condition. 
Likewise, the five-~and six-point IB kernels satisfy the same properties, except that the five-point kernel does not satisfy the even-odd condition.

This study also considers the performance of B-spline kernels (Figures~\ref{fig:bspline_3}--\ref{fig:bspline_6}), which are recursively constructed by convolution against the zeroth-order B-spline kernel (which is a piecewise-constant function):
\begin{align}
\label{eq:bspline_kernels}  \varphi^{\text{B-spline}}_{n}(r) &= \varphi^{\text{B-spline}}_{n-1}(r) \ast \varphi^{\text{B-spline}}_0(r) = \int_{-\infty}^\infty \varphi^{\text{B-spline}}_{n-1}(r-s) \, \varphi^{\text{B-spline}}_0(s) \, \mathrm{d}s.
\end{align} 
An $n^\text{th}$-order B-spline satisfies up to $n^\text{th}$-order moment conditions but does not satisfy the even--odd condition or the approximate grid translational invariance property.
Both the radius of support and the smoothness of the B-spline kernel increases with order, and the limiting function is a Gaussian~\cite{Unser1992, Bao2017}. One advantage of using B-spline kernels is that they are piecewise polynomial and can be evaluated efficiently. Table~\ref{table:kernels} shows a summary of properties and moment conditions satisfied by the kernels that are considered in this paper.

\begin{table}[t!!]
    \setlength\tabcolsep{4.5pt}
    \scriptsize
	\centering	
	\caption{Selected choices of regularized delta functions with properties and moment conditions that are satisfied. In the higher moment columns (second--fifth moment), the value of $K$ that satisfies the given moment condition is given.}
\begin{tabular}{| c | c | c | c | c | c | c | c | c|l}
\cline{1-9}
Kernel & Even--Odd & \makecell{Zeroth\\Moment} & \makecell{First\\Moment} & \makecell{Second\\Moment} & \makecell{Third\\Moment} & \makecell{Fourth\\Moment} & \makecell{Fifth\\Moment} & \makecell{Sum of\\Squares} & \\ 
\cline{1-9}
\multicolumn{1}{ |c| }{Piecewise-linear} & $\Cross$ & $\checkmark$ & $\checkmark$ & $\Cross$ & $\Cross$ & $\Cross$ & $\Cross$ & $\Cross$ &    \\ 
\cline{1-9}
\multicolumn{1}{ |c| }{IB (3-point)} & $\Cross$ & $\checkmark$ & $\checkmark$ & $\Cross$ & $\Cross$ & $\Cross$ & $\Cross$ & $\frac{1}{2}$ &    \\ 
\cline{1-9}
\multicolumn{1}{ |c| }{IB (4-point)} & $\checkmark$ & $\checkmark$ & $\checkmark$ & $\Cross$ & $\Cross$ & $\Cross$ & $\Cross$ & $\frac{3}{8}$ &    \\ 
\cline{1-9}
\multicolumn{1}{ |c| }{IB (5-point)} & $\Cross$ & $\checkmark$ & $\checkmark$ & $\frac{38}{60}-\frac{\sqrt{69}}{60}$ & $\approx 0$ & $\Cross$ & $\Cross$ & $\approx 0.393$ &    \\ 
\cline{1-9}
\multicolumn{1}{ |c| }{IB (6-point)} & $\checkmark$ & $\checkmark$ & $\checkmark$ & $\frac{59}{60}-\frac{\sqrt{29}}{20}$ & $\approx 0$ & $\Cross$ & $\Cross$ & $\approx 0.326$ &    \\ 
\cline{1-9}
\multicolumn{1}{ |c| }{B-spline (3-point)} & $\Cross$ & $\checkmark$ & $\checkmark$ & $\frac{1}{4}$ & $\Cross$ & $\Cross$ & $\Cross$ & $\Cross$ &    \\ 
\cline{1-9}
\multicolumn{1}{ |c| }{B-spline (4-point)} & $\Cross$ & $\checkmark$ & $\checkmark$ & $\frac{1}{3}$ & $0$ & $\Cross$ & $\Cross$ & $\Cross$ &    \\ 
\cline{1-9}
\multicolumn{1}{ |c| }{B-spline (5-point)} & $\Cross$ & $\checkmark$ & $\checkmark$ & $\approx 0.417$ & $\approx 0$ & $\approx 0.479$ & $\Cross$ & $\Cross$ &    \\ 
\cline{1-9}
\multicolumn{1}{ |c| }{B-spline (6-point)} & $\Cross$ & $\checkmark$ & $\checkmark$ & $\approx 0.500$ & $\approx 0$ & $\approx 0.700$ & $0$ & $\Cross$ &    \\ 
\cline{1-9}
\end{tabular}
	\label{table:kernels}	
\end{table} 

\subsection{Lagrangian mesh spacing}
\label{subsec:lag_spacing}
\begin{figure}[t!!]
 \centering
 \includegraphics[scale = 0.225]{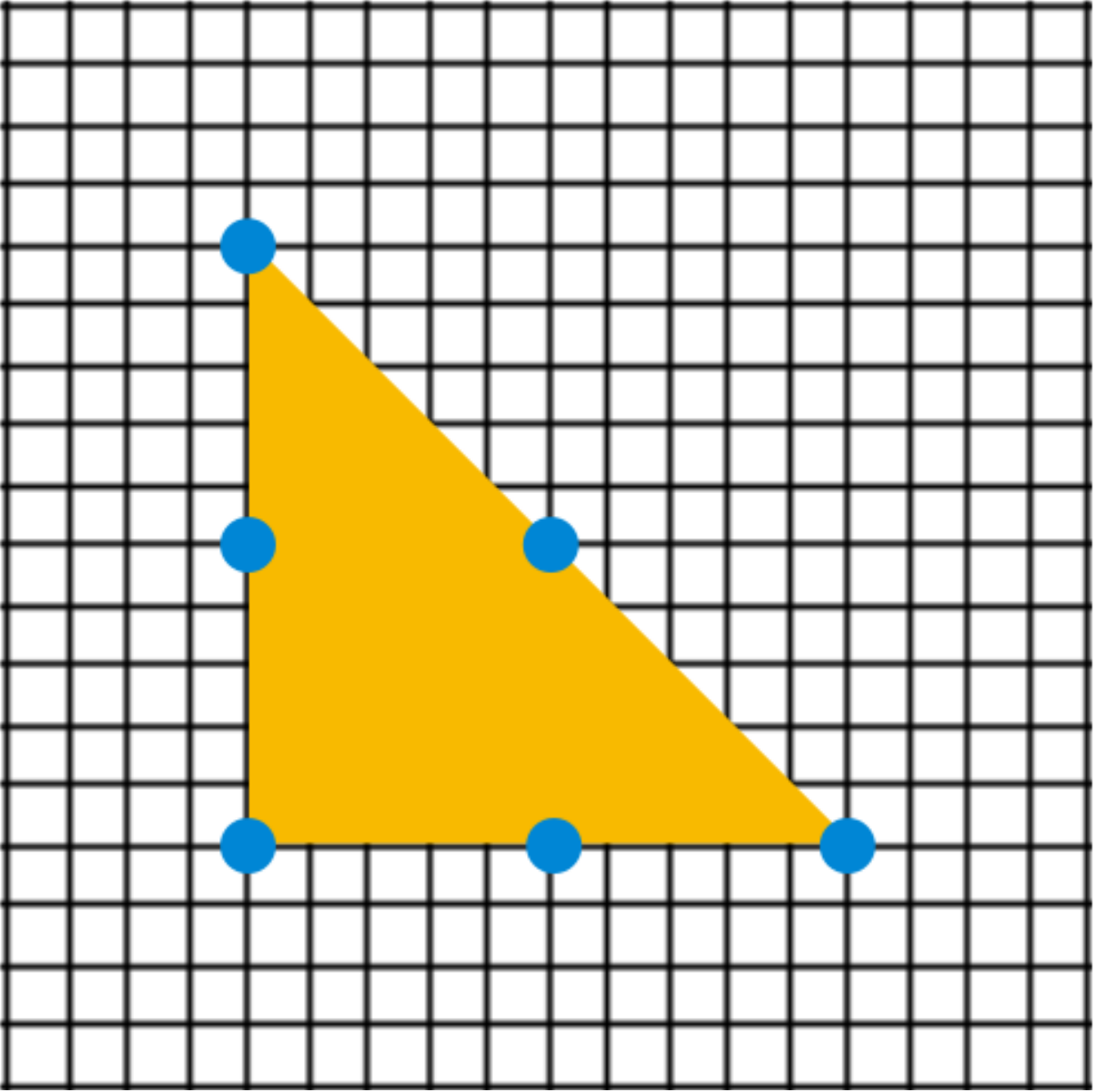}
 \caption{Description of $\mfac$, which is the ratio between background Cartesian grid spacing and finite element node spacing. In the case shown here, there are about five Cartesian grid cells between two finite element nodes for a second-order triangular ($\Ptwo$) element, so we say that $\mfac \approx 5$.}
 \label{fig:mfac}
\end{figure}

In addition to the choice of the regularized delta function kernel, the coupling strategy used in the IFED method allows us to study the impact of the interaction between the Lagrangian mesh and the Eulerian grid. 
We use \emph{mesh factor}, $\mfac$, to indicate the approximate ratio of Lagrangian element node spacing to the Eulerian grid spacing.
$\mfac$ is defined here as $\mfac \approx \frac{\lagrangiandx}{\efac \euleriandx}$, in which $\lagrangiandx$ is the Lagrangian element size, $\euleriandx$ is the Eulerian grid spacing in each coordinate direction, and the \emph{element factor} $\efac$ is $1$ for linear elements and $2$ for quadratic elements, and reflects the fact that, e.g., nodes are approximately $\lagrangiandx/2$ apart for quadratic elements.
See Figure~\ref{fig:mfac}.
For example, the usual ``rule of thumb'' described by Peskin~\cite{Peskin2002} restricts the structural mesh to be approximately twice as fine as the Eulerian grid, which corresponds to $\mfac=0.5$, independent of $\efac$.
We investigate the effect of the choice of $\mfac$, along with the choice of the regularized delta function kernel, in the accuracy of our solutions through our FSI benchmarks.

\subsection{Time discretization}
\label{subsec:time_discret}
We use an explicit midpoint rule for the structural deformation, a Crank-Nicolson scheme for the viscous term, and an Adams-Bashforth scheme for the convective term, as detailed previously~\cite{Griffith2017}. Each time step involves solving the time-dependent incompressible Stokes equations, one force evaluation and force spreading operation, and two velocity interpolation operations.

\subsection{Stabilization method for hyperelastic material models}
\label{subsec:volume_stabilization}
In the continuous IFED formulation, the immersed structure is automatically treated as incompressible because $\D{\vchi}{t}(\X,t) = \u(\vchi(\X,t),t)$ and $\grad \cdot \u(\x,t) = 0$.
In the spatially discretized equations, exact incompressibility can be lost in the solid.
We use a stabilization approach~\cite{Vadala-Roth2020} that effectively reinforces the incompressibility constraint.
This approach uses a splitting of the strain energy functional into isochoric and volumetric parts,
\begin{equation}
\Psi(\FF) = W(\overline{\FF}) + U(J),
\end{equation}
in which $\overline{\FF} = J^{-1/3}\FF$, as is commonly done in nearly incompressible elasticity models~\cite{Sansour2008}. 
We use the volumetric part of the strain energy as a stabilization term used to enforce the incompressibility of the elastic structures, and here we choose it to be~\cite{Liu1994}
\begin{equation}
    \label{eq:volumetric} U(J) = \beta(J\ln{J} - J + 1),
\end{equation}
in which $\beta$ is a numerical bulk modulus~\cite{Vadala-Roth2020}.
In this study, we empirically determine approximately the largest value of the bulk modulus that allow the scheme to remain stable for a given time step size $\Delta t$ to penalize any volume change in the structural mesh elements for each kernel and grid spacing.

\subsection{IBAMR}
\label{subsec:software}
FSI simulations use the IBAMR software infrastructure, which is a distributed-memory parallel implementation of the IB method with adaptive mesh refinement (AMR)~\cite{Griffith2007, IBAMR}. 
IBAMR uses SAMRAI~\cite{SAMRAIPaper} for Cartesian grid discretization management, libMesh~\cite{libMeshPaper} for finite element discretization management, and PETSc~\cite{petsc-user-ref} for linear solver infrastructure.

\section{Fluid-Structure Interaction Benchmarks and Results}
\label{sec:benchmarks}

This section systematically investigates the impact of the choice of regularized delta function as well as the relative spacings of the Lagrangian and Eulerian discretization on the IFED method using a series of FSI benchmarks.
The benchmarks are organized into shear-dominated and pressure-loaded cases, and we apply our findings from them to a large-scale FSI model.

\subsection{Two-dimensional flow past cylinder}
\label{subsec:cylinder}
\begin{figure}[t!!]
	\begin{center}
		\begin{subfigure}[t]{0.375\textwidth}
			\centering
			\includegraphics[scale = 0.65]{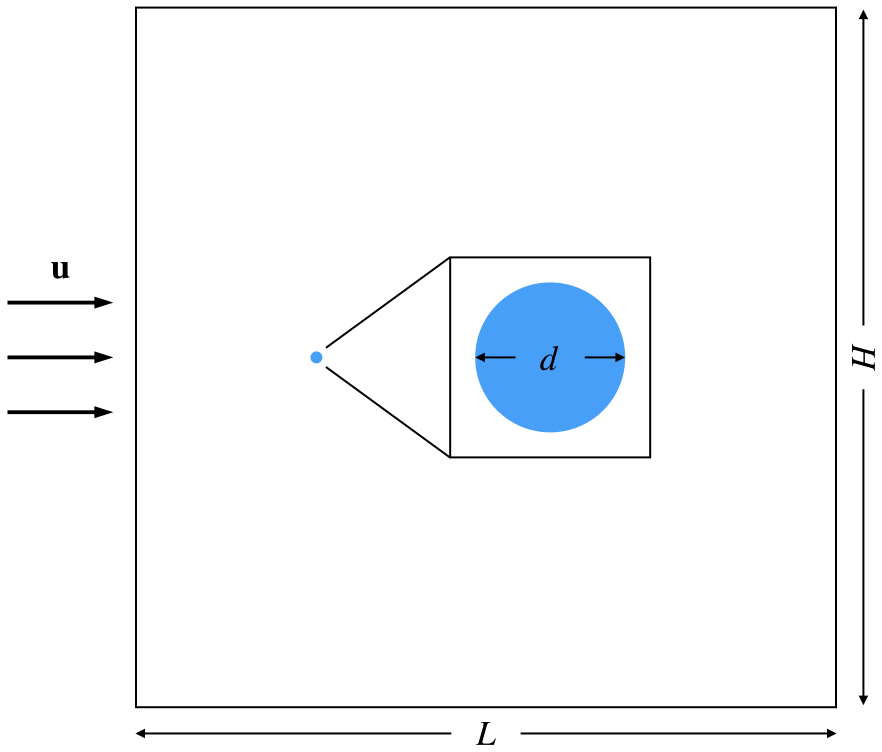}
			\caption{} \label{fig:cylinder_schematic}
		\end{subfigure}\hspace{0.05in}
		\begin{subfigure}[t]{0.6\textwidth}
			\centering
			\includegraphics[scale = 0.25]{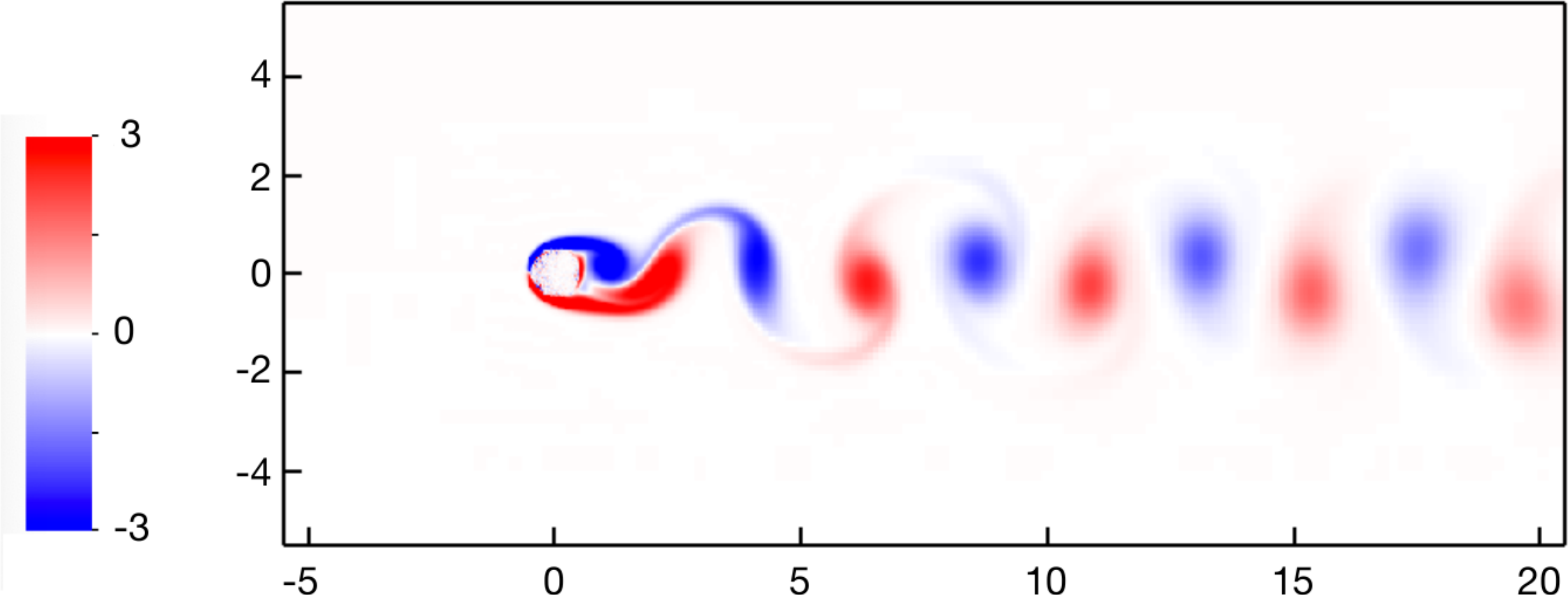}
			\caption{} \label{fig:cylinder_simulation}
		\end{subfigure}
		\caption{(a) Schematic of two-dimensional flow past a cylinder benchmark. Arrows represent the inflow boundary, where a uniform velocity boundary condition, $\u = (1,0)$, is applied. Zero normal traction and tangential velocity at the outflow boundary. For the top and bottom boundaries, we use zero normal velocity and tangential traction. We choose $Re = 200$ in our tests. (b) A magnified view of the vortices shed from a stationary circular cylinder from our simulation.}
		\label{fig:flow_past_cylinder}		
	\end{center}
\end{figure}

We begin by considering the widely used test of viscous incompressible flow past a stationary circular cylinder~\cite{Taira2007, Griffith2017}.
We use the penalty formulation, Eq.~\eqref{eq:tether_force}, to model the cylinder. 
The penalty parameters $\kappa$ and $\eta$ are determined to be approximately the largest stable values for a given time step size and Lagrangian and Eulerian mesh spacings.
The cylinder has diameter $D = 1$ and is embedded in a computational domain $\Omega$ with side lengths of $L = H = 60$. 
Figure~\ref{fig:cylinder_schematic} provides a schematic diagram. 
We use a uniform inflow velocity boundary condition, $\u = (1,0)$, on the left boundary of the computational domain and specify zero normal traction and tangential velocity conditions on the right boundary. 
For the top and bottom boundaries of the computational domain, we specify zero normal velocity and tangential traction conditions. 
The fluid has density $\rho = 1$ and viscosity $\mu = 0.005$, and the Reynolds number is $Re = \frac{\rho u_\infty D}{\mu} = 200$. 
We use the drag $\left(\CD = \frac{F_x}{\rho u_\infty^2D/2}\right)$ and lift $\left(\CL = \frac{F_y}{\rho u_\infty^2D/2}\right)$ coefficients to evaluate the effect of the choice of regularized delta function or mesh factor on the computed dynamics, in which $\F = (F_x, F_y)$ is the net force on the cylinder and $u_\infty$ is the characteristic flow speed (which we take to be $x$-component of the inflow velocity).
The computational domain is discretized using a six-level locally refined grid with a refinement ratio of two between levels and an $N\times N$ coarse grid.
The fine-grid Cartesian cell size is $\euleriandx = H/(32N)$, and the time step size is $\Delta t = 0.1875/N$.
Griffith and Luo~\cite{Griffith2017} previously conducted an initial study using this benchmark with the three-, four-, and six-point IB delta function kernels.
\begin{figure}[t!!]
 \centering
 \includegraphics[scale = 0.2375]{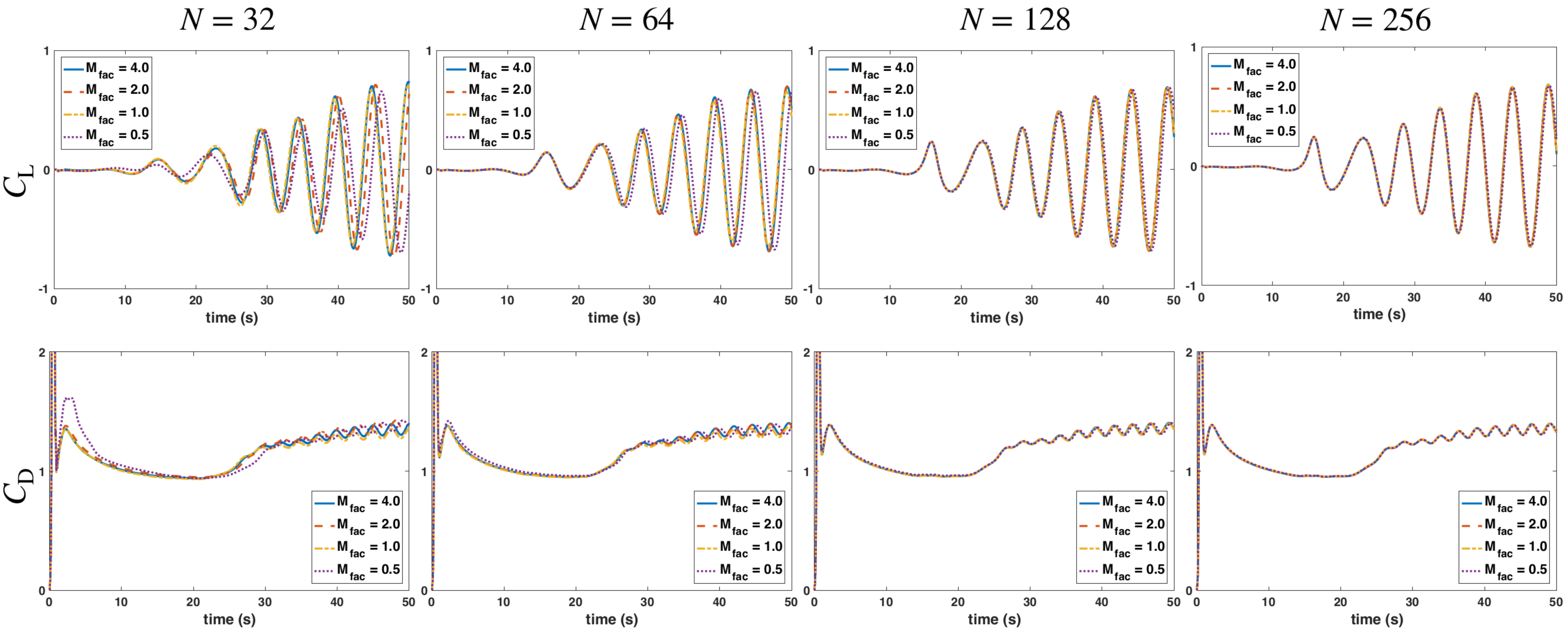}
 \caption{Representative lift ($\CL$) and drag ($\CD$) coefficients for flow past a stationary cylinder at $Re = 200$ using the three-point B-spline kernel. 
The computational domain $\Omega$ is discretized using a six-level locally refined grid with a refinement ratio of two between levels and an $N\times N$ coarse grid. 
We observe that the lift and drag coefficients show converging behavior under simultaneous Lagrangian and Eulerian grid refinement (from $N=32$ to $N=256$) for all values of $\mfac$.
Similar accuracy is observed with the other kernels except for the six-point IB kernel, which requires higher resolution to yield comparable accuracy.}
 \label{fig:cyl_grid_conv}
\end{figure}
Figure~\ref{fig:cyl_grid_conv} shows representative results of lift and drag coefficients using the three-point B-spline kernel, which shows converging behavior under simultaneous Lagrangian and Eulerian grid refinement (from $N=32$ to $N=256$) for $\mfac=0.5,1,2,$ and 4.
The method yields similar convergence behavior under grid refinement for the other kernel functions that we consider and for the chosen values of $\mfac$.

Although the scheme converges under grid refinement for all choices of kernels and for all values of $\mfac$, we observe that four- and six-point IB kernels clearly require high grid resolution to yield converged solutions for $\mfac = 0.5$ and 1.
Considering specifically the intermediate Cartesian resolution corresponding to $N=64$, we find that some kernels show markedly lower accuracy for some $\mfac$ values.
At the same moderate resolution, the three-point IB kernel (along with the piecewise linear and B-spline kernels) are less sensitive for $\mfac \geq 1$ compared to the four- and six-point IB kernels, which agrees with the results reported by Griffith and Luo~\cite{Griffith2017}.

\begin{figure}[t!!]
 \centering
 \includegraphics[scale = 0.2375]{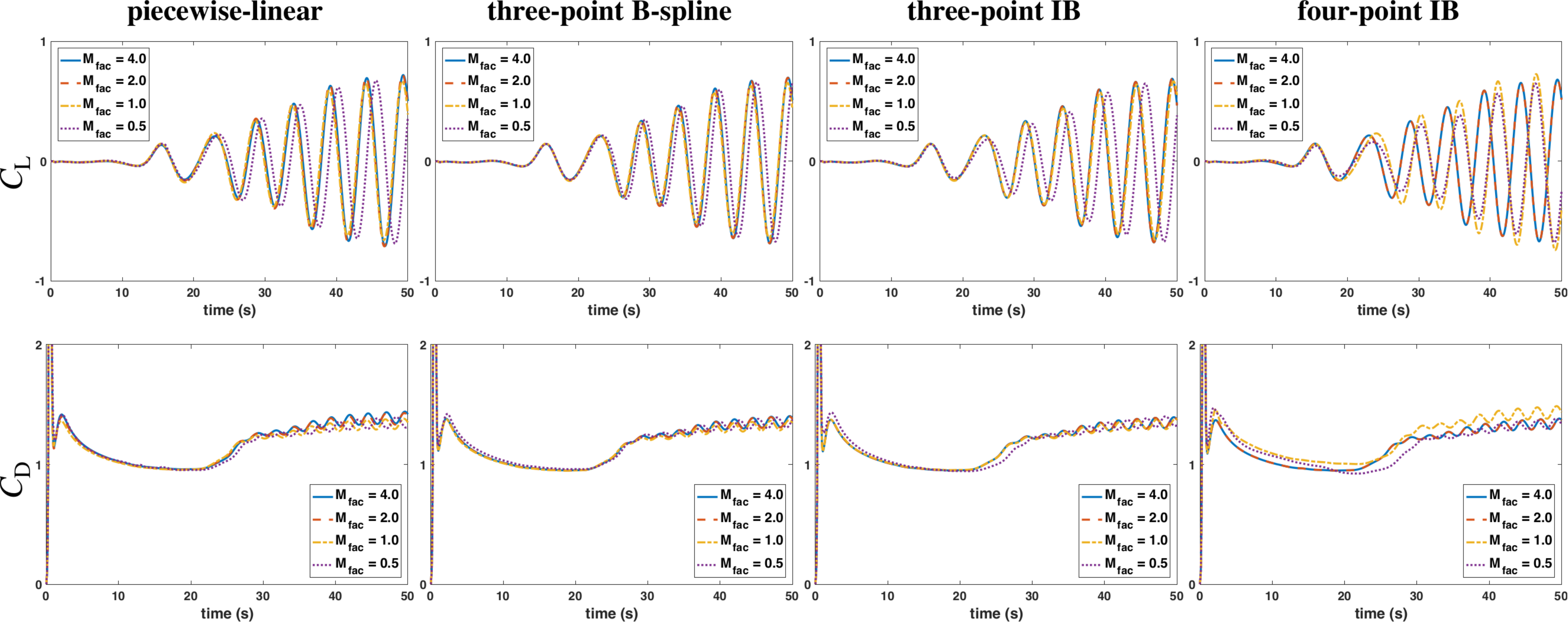}
 \caption{Representative results of lift ($\CL$) and drag ($\CD$) coefficients for flow past a stationary cylinder at $N = 64$ using four representative regularized delta functions and different relative structural mesh spacing for a fixed Eulerian grid ($\mfac$ = 0.5, 1, 2, and 4).}
 \label{fig:cyl_comparison}
\end{figure}

Next, we compare results at the same resolution with a broader selection of kernel functions to identify which kernels give more consistent results over different values of $\mfac$ at intermediate resolution.
Table~\ref{table:flow_past_cylinder_comparison} compares lift and drag coefficients and Strouhal numbers at $N=64$ using different kernel functions for $\mfac$ = 0.5, 1, 2, and 4.
These quantities converge to $\CL = \pm0.67$, $\CD = 1.361\pm0.041$, and $St=0.200$ under grid refinement, and we observe that the three-point IB and three- and four-point B-spline kernels with $\mfac > 1$ result in the best agreement with the converged values at $N=64$. 
For $\mfac = 0.5$, lift amplitudes differ up to 25\% from the converged value, compared to up to 9\% for $\mfac \geq 1$ for some kernels.
These three kernels also give consistent Strouhal numbers ($St=0.200$) for the values of $\mfac$ considered.
Figure~\ref{fig:cyl_comparison} compares lift and drag coefficients as functions of time for four representative kernels.
Although Table~\ref{table:flow_past_cylinder_comparison} suggests that the three-point IB and three-point B-spline kernels yield similar values for lift and drag coefficients with similar root-mean-square error with respect to the converged results, Figure~\ref{fig:cyl_comparison} shows that the lift and drag coefficients for $\mfac = 0.5$ using the three-point IB kernel clearly deviates from the results for other values of $\mfac$.
This suggests that the three-point B-spline kernel yields more consistent results as we vary $\mfac$.
These results also are concordant with previous work by Griffith and Luo~\cite{Griffith2017} in that refining the Lagrangian mesh while keeping the Eulerian grid resolution fixed generally lowers the accuracy.
We find that the three-point B-spline kernel shows the least sensitivity at the coarsest grid spacings amongst the kernel functions considered in this study.
Possible explanations for relatively lower accuracy and consistency from other kernels could be that the piecewise linear kernel is not sufficiently smoothing out the high-frequency errors and the larger IB and B-spline kernels are generating unphysically large numerical boundary layers near fluid-structure interfaces.
We also note that the four-and six-point IB kernels satisfy the even--odd condition, and they yield lower accuracy than the corresponding three- and five-point IB kernels that do not satisfy the even-odd condition.

\begin{table}[t!!]
    \setlength\tabcolsep{4.5pt}
    \scriptsize
	\centering	
	\caption{Comparison of lift ($\CL$) and drag ($\CD$) coefficients for flow past a stationary cylinder at an intermediate Cartesian resolution of $N=64$ using different regularized delta functions and relative structural mesh spacing ($\mfac$). These values converge to $\CL = \pm 0.67$, $\CD = 1.361\pm0.041$, $St = 0.200$ under further grid refinement. We observe that the three-point IB and three- and four-point B-spline kernels result in the highest accuracies and least variation across $\mfac$ values at $N=64$.}
\begin{tabular}{c c c c |c c c | c c c| c c c|l}
\cline{2-13}
& \multicolumn{3}{ |c| }{$\mfac = 0.5$} & \multicolumn{3}{ |c| }{$\mfac = 1.0$} & \multicolumn{3}{ |c| }{$\mfac = 2.0$} & \multicolumn{3}{ |c| }{$\mfac = 4.0$} \\  
\cline{1-13}
\multicolumn{1}{ |c| }{Kernel} & $\CL$ & $\CD$ & $St$ & $\CL$ & $\CD$ & $St$ & $\CL$ & $\CD$ & $St$ & $\CL$ & $\CD$ & $St$ & \\ 
\cline{1-13}
\multicolumn{1}{ |c| }{Piecewise-linear} & \parbox[b]{0.8cm}{$\\ \pm 0.45$} &   \parbox[b]{0.9cm}{$\\1.350\\ \pm 0.045$} & 0.180  &  \parbox[b]{0.8cm}{$\\ \pm 0.61$} & \parbox[b]{0.9cm}{$\\1.346\\ \pm 0.030$}   &  0.200 & \parbox[b]{0.8cm}{$\\ \pm 0.69$}  &  \parbox[b]{0.9cm}{$\\1.389\\ \pm 0.039$}  & 0.200 & \parbox[b]{0.8cm}{$\\ \pm 0.70$} & \parbox[b]{0.9cm}{$\\1.400\\ \pm 0.042$}  & 0.200    \\ 
\cline{1-13}
\multicolumn{1}{ |c| }{IB (3-point)} & \parbox[b]{0.8cm}{$\\ \pm 0.53$} &   \parbox[b]{0.9cm}{$\\1.375\\ \pm 0.045$} & 0.200  &  \parbox[b]{0.8cm}{$\\ \pm 0.61$} & \parbox[b]{0.9cm}{$\\1.347\\ \pm 0.028$}   &  0.200 & \parbox[b]{0.8cm}{$\\ \pm 0.66$}  &  \parbox[b]{0.9cm}{$\\1.357\\ \pm 0.036$}  & 0.200 & \parbox[b]{0.8cm}{$\\ \pm 0.66$} & \parbox[b]{0.9cm}{$\\1.358\\ \pm 0.036$}  & 0.200   \\ 
\cline{1-13}
\multicolumn{1}{ |c| }{IB (4-point)} & \parbox[b]{0.8cm}{$\\ \pm 0.44$} &   \parbox[b]{0.9cm}{$\\1.359\\ \pm 0.042$} & 0.200  &  \parbox[b]{0.8cm}{$\\ \pm 0.62$} & \parbox[b]{0.9cm}{$\\1.446\\ \pm 0.047$}   &  0.180 & \parbox[b]{0.8cm}{$\\ \pm 0.64$}  &  \parbox[b]{0.9cm}{$\\1.347\\ \pm 0.031$}  & 0.200 & \parbox[b]{0.8cm}{$\\ \pm 0.64$} & \parbox[b]{0.9cm}{$\\1.348\\ \pm 0.031$}  & 0.200   \\ 
\cline{1-13}
\multicolumn{1}{ |c| }{IB (5-point)} & \parbox[b]{0.8cm}{$\\ \pm 0.46$} &   \parbox[b]{0.9cm}{$\\1.360\\ \pm 0.044$} & 0.180  &  \parbox[b]{0.8cm}{$\\ \pm 0.55$} & \parbox[b]{0.9cm}{$\\1.432\\ \pm 0.053$}   &  0.200 & \parbox[b]{0.8cm}{$\\ \pm 0.64$}  &  \parbox[b]{0.9cm}{$\\1.346\\ \pm 0.032$}  & 0.200 & \parbox[b]{0.8cm}{$\\ \pm 0.64$} & \parbox[b]{0.9cm}{$\\1.343\\ \pm 0.033$}  & 0.200  \\ 
\cline{1-13}
\multicolumn{1}{ |c| }{IB (6-point)} & \parbox[b]{0.8cm}{$\\ \pm 0.51$} &   \parbox[b]{0.9cm}{$\\1.366\\ \pm 0.042$} & 0.180  &  \parbox[b]{0.8cm}{$\\ \pm 0.70$} & \parbox[b]{0.9cm}{$\\1.467\\ \pm 0.043$}   &  0.180 & \parbox[b]{0.8cm}{$\\ \pm 0.63$}  &  \parbox[b]{0.9cm}{$\\1.332\\ \pm 0.029$}  & 0.200 & \parbox[b]{0.8cm}{$\\ \pm 0.63$} & \parbox[b]{0.9cm}{$\\1.332\\ \pm 0.029$}  & 0.200    \\ 
\cline{1-13}
\multicolumn{1}{ |c| }{B-spline (3-point)} & \parbox[b]{0.8cm}{$\\ \pm 0.51$} &   \parbox[b]{0.9cm}{$\\1.354\\ \pm 0.042$} & 0.200  &  \parbox[b]{0.8cm}{$\\ \pm 0.62$} & \parbox[b]{0.9cm}{$\\1.336\\ \pm 0.032$}   &  0.200 & \parbox[b]{0.8cm}{$\\ \pm 0.67$}  &  \parbox[b]{0.9cm}{$\\1.363\\ \pm 0.037$}  & 0.200 & \parbox[b]{0.8cm}{$\\ \pm 0.67$} & \parbox[b]{0.9cm}{$\\1.366\\ \pm 0.037$}  & 0.200   \\ 
\cline{1-13}
\multicolumn{1}{ |c| }{B-spline (4-point)} & \parbox[b]{0.8cm}{$\\ \pm 0.50$} &   \parbox[b]{0.9cm}{$\\1.350\\ \pm 0.043$} & 0.200  &  \parbox[b]{0.8cm}{$\\ \pm 0.61$} & \parbox[b]{0.9cm}{$\\1.355\\ \pm 0.031$}   &  0.200 & \parbox[b]{0.8cm}{$\\ \pm 0.66$}  &  \parbox[b]{0.9cm}{$\\1.357\\ \pm 0.035$}  & 0.200 & \parbox[b]{0.8cm}{$\\ \pm 0.66$} & \parbox[b]{0.9cm}{$\\1.356\\ \pm 0.035$}  & 0.200   \\ 
\cline{1-13}
\multicolumn{1}{ |c| }{B-spline (5-point)} & \parbox[b]{0.8cm}{$\\ \pm 0.49$} &   \parbox[b]{0.9cm}{$\\1.358\\ \pm 0.043$} & 0.200  &  \parbox[b]{0.8cm}{$\\ \pm 0.60$} & \parbox[b]{0.9cm}{$\\1.389\\ \pm 0.040$}   &  0.200 & \parbox[b]{0.8cm}{$\\ \pm 0.65$}  &  \parbox[b]{0.9cm}{$\\1.351\\ \pm 0.034$}  & 0.200 & \parbox[b]{0.8cm}{$\\ \pm 0.65$} & \parbox[b]{0.9cm}{$\\1.349\\ \pm 0.034$}  & 0.200   \\ 
\cline{1-13}
\multicolumn{1}{ |c| }{B-spline (6-point)} & \parbox[b]{0.8cm}{$\\ \pm 0.49$} &   \parbox[b]{0.9cm}{$\\1.368\\ \pm 0.044$} & 0.180  &  \parbox[b]{0.8cm}{$\\ \pm 0.56$} & \parbox[b]{0.9cm}{$\\1.422\\ \pm 0.051$}   &  0.200 & \parbox[b]{0.8cm}{$\\ \pm 0.64$}  &  \parbox[b]{0.9cm}{$\\1.346\\ \pm 0.032$}  & 0.200 & \parbox[b]{0.8cm}{$\\ \pm 0.64$} & \parbox[b]{0.9cm}{$\\1.344\\ \pm 0.032$}  & 0.200   \\ 
\cline{1-13}
\end{tabular}
	\label{table:flow_past_cylinder_comparison}	
\end{table} 

% SLANTED CHANNEL FLOW BENCHMARK
\subsection{Two-dimensional channel flow}
\label{subsec:channel}
\begin{figure}[t!!]
	\begin{center}
		\hspace{-0.55in}
		\begin{subfigure}[t]{0.4\textwidth}
			\centering
			\includegraphics[scale = 0.175]{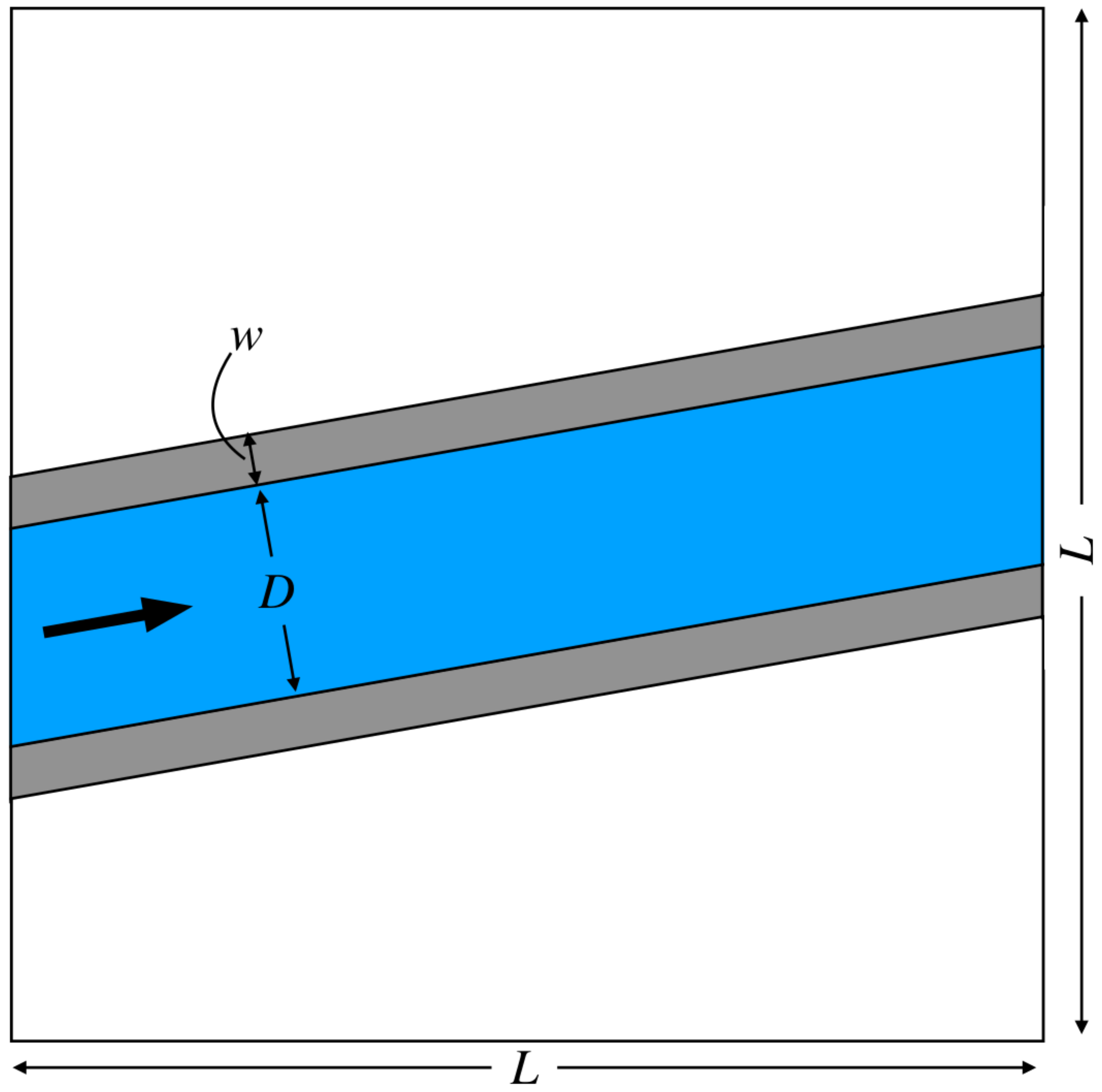}
			\caption{} \label{fig:slanted_channel}
		\end{subfigure}\hspace{0.1in}
		\begin{subfigure}[t]{0.4\textwidth}
			\centering
			\includegraphics[scale = 0.175]{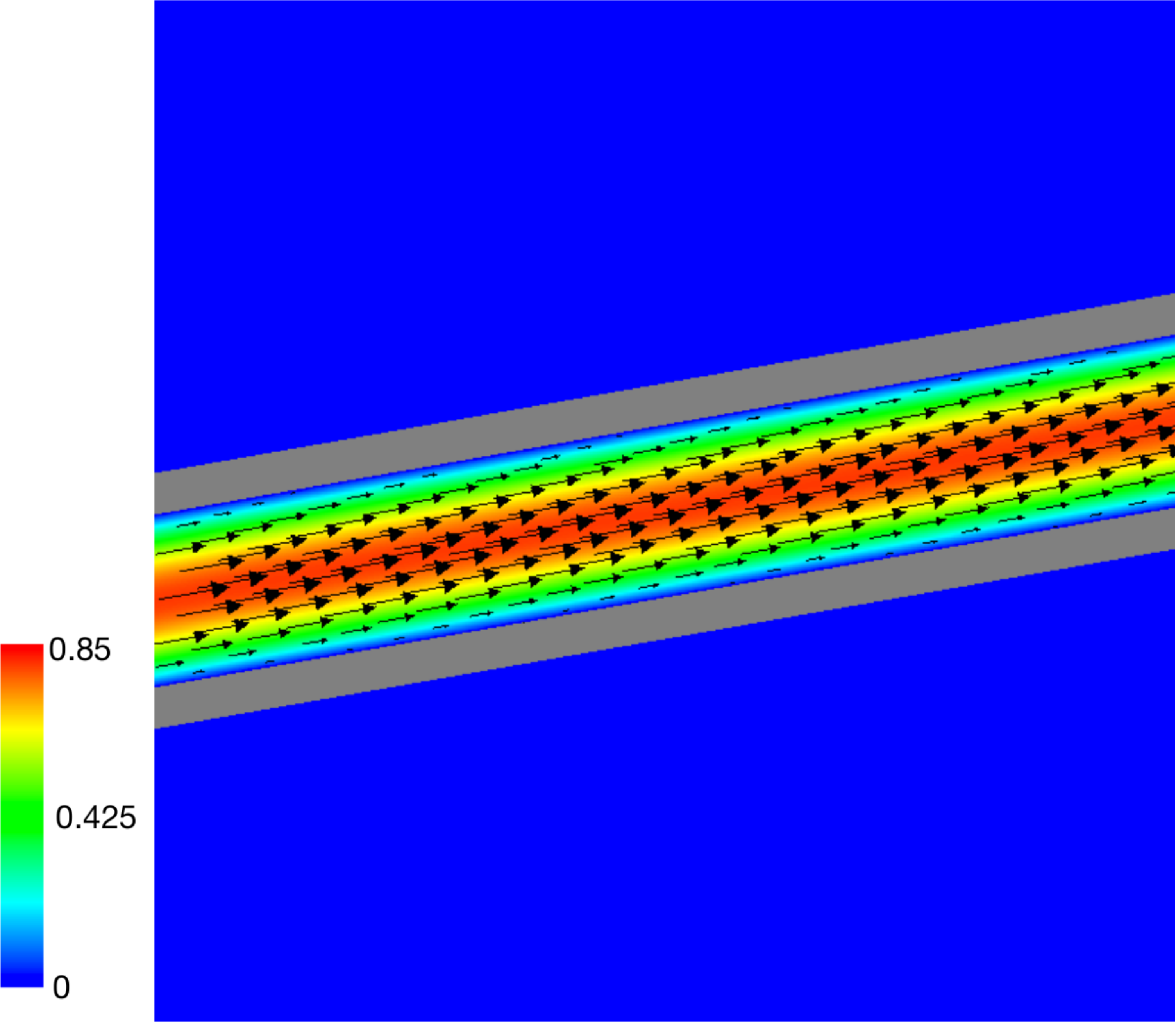}
			\caption{} \label{fig-slanted_channel}
		\end{subfigure}
		\caption{(a) Schematic of two-dimensional flow through a slanted channel. (b) Representative steady-state velocity solution vector field for the two-dimensional slanted channel flow benchmark. This simulation uses a three-level locally refined grid with a refinement ratio of two between levels and an $N \times N$ coarse grid with $N = 256$. The computation uses a piecewise linear kernel and $\mfac = 2.0$.}
		\label{fig:channel_schematic}		
	\end{center}
\end{figure}

This section considers the benchmark problem of two-dimensional channel flow test adopted from Kolahdouz et al.~\cite{Kolahdouz2020}.
We consider a domain $\Omega = [0, L]^2$ with two parallel plates, with channel width $D$ and wall width $w = 0.24D$.
The exact steady-state solution is described by the plane Poiseuille equation,
\begin{align} 
\label{eqn_plane_poiseuille_straight} u(y) &= \frac{\chi D}{2\mu } (y-y_0)\left(1 -\frac{y-y_0}{D}\right),
\end{align}
in which $y_0$ is the height of inner wall of the lower channel plate and $\chi = \frac{2p_0}{L}$ is the pressure gradient between the inflow and the outflow. 
To avoid a purely grid aligned test, we consider a slanted channel. 
Figure~\ref{fig:slanted_channel} provides a schematic.
This is done by rotating the channel walls by an angle $\theta$, so that for every point on the walls $(x,y)$, we transform the $y$-coordinate to $y' = y+\left(x-\frac{L}{2}\right)\tan\theta$ and let $(x,y')$ be the new coordinates for the walls.
The steady-state solution is then transformed to
\begin{align}
\label{eqn_plane_poiseuille} u(\eta) &= \frac{\chi D}{2\mu } (\eta-\eta_0)\left(1 -\frac{\eta-\eta_0}{D}\right),
\end{align}
in which $\eta = -x\sin\theta+(y-y_0)\cos\theta$ and $\chi=\frac{2p_0}{L/\cos\theta+D\tan\theta}$.
In our computations, we use $D = 1$, $\mu = 0.01$, $\rho = 1.0$, $L = 6D$, $p_0 = 0.2$, and $\theta=\pi/18$.
The maximum velocity is $U_\text{max} = 1$, and the average velocity is $\overline{U} = 2/3$, which implies that the Reynolds number is $Re = \frac{\rho D \overline{U}}{\mu} \approx 66.67$. 
The fine-grid Cartesian cell size is $\euleriandx = D/(4N)$, and the time step size is $\Delta t = 0.0375/N$.
At the inlet and outlet of the channel, the rotated analytical solution of the steady-state velocity (Eq.~\eqref{eqn_plane_poiseuille}) provides velocity boundary conditions.
This benchmark assesses which choices of kernel and $\mfac$ give the best accuracy for the flow within a confined, stationary geometry.

The channel walls are modeled as a stiff neo-Hookean material with $W_{\text{wall}} = \frac{c_{\text{wall}}}{2}(\bar{I}_1 - 3)$, in which $\bar{I}_1$ is the modified first invariant of the right Cauchy-Green tensor $\bar{\CC} = \bar{\FF}^T\bar{\FF} = J^{-\frac{2}{3}} \FF^T\FF$, and $c_{\text{wall}}\propto \frac{\Delta x}{\Delta t^2}$ is a penalty stiffness parameter, so that the body becomes infinitely rigid as $\Delta t \to 0$. 
In addition to the penalty stiffness, we use penalty body and damping forces described in Eq.~\eqref{eq:tether_force} to enforce rigidity of the structure as well as to keep it stationary. 
We empirically determine approximately the largest values of the penalty parameters that allow the scheme to remain stable for a given time step size $\Delta t$ for each kernel and grid spacing.

\begin{figure}[t!!]
 \centering
 \includegraphics[scale = 0.25]{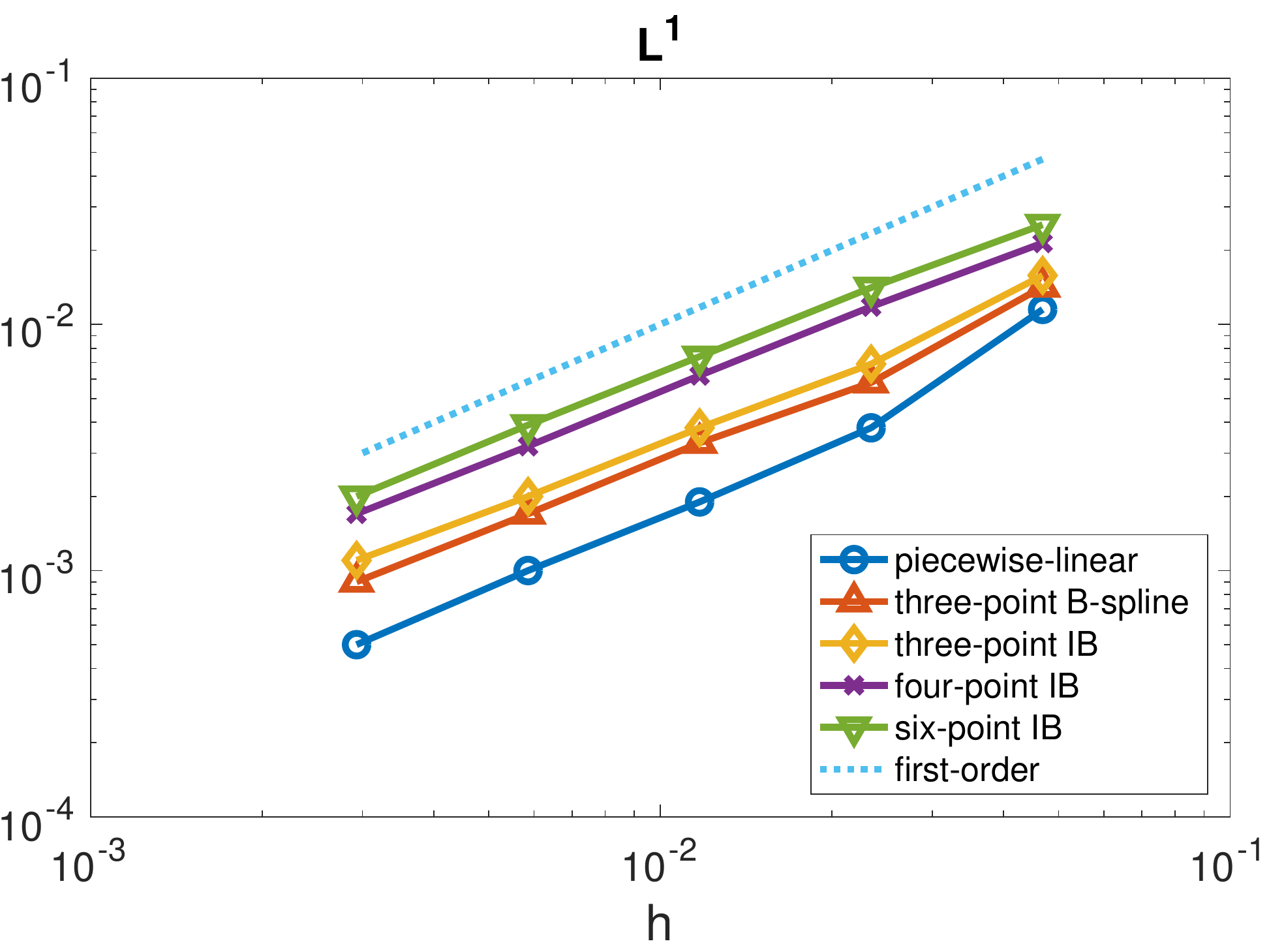}
 \includegraphics[scale = 0.25]{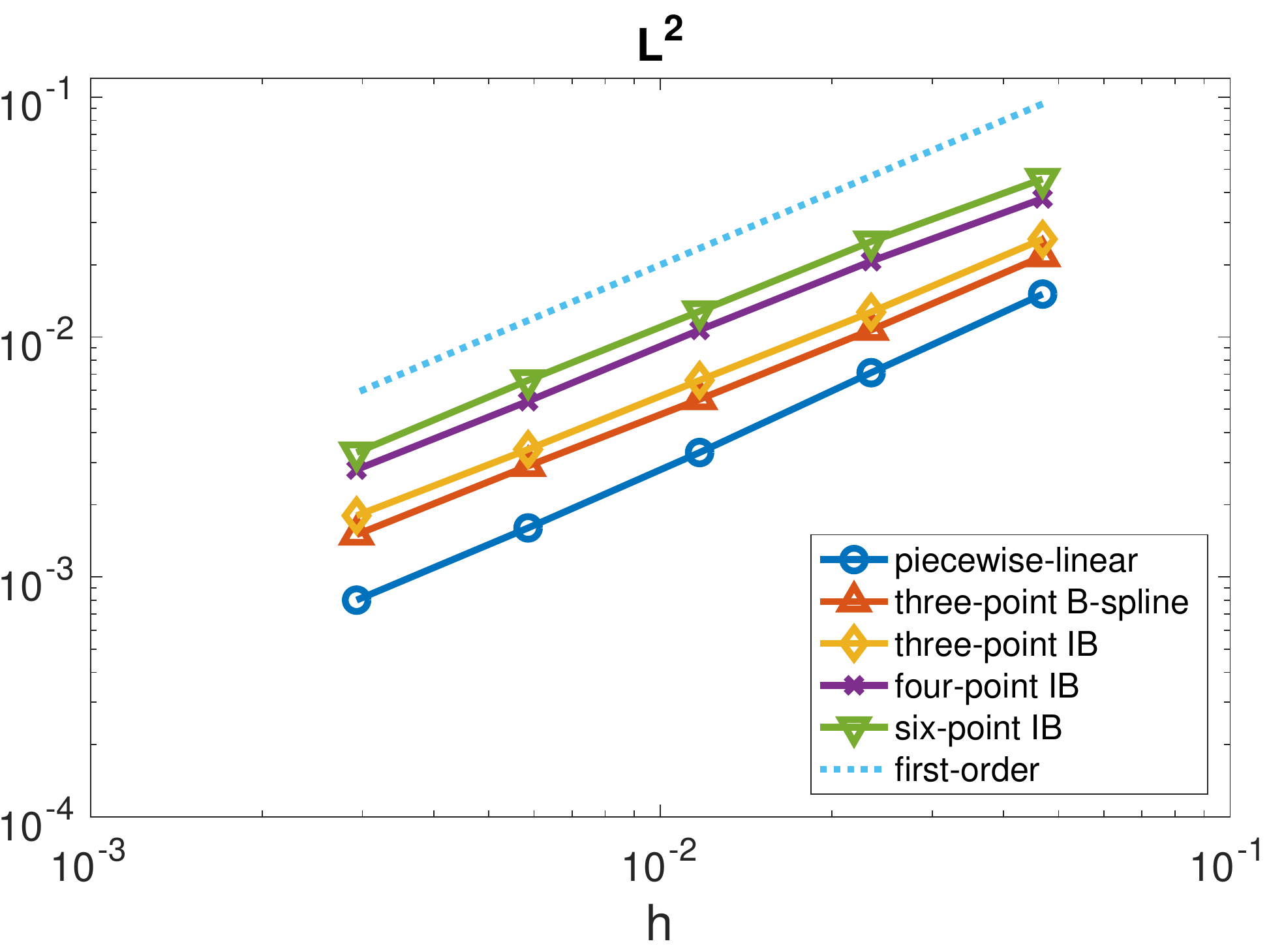}
 \includegraphics[scale = 0.25]{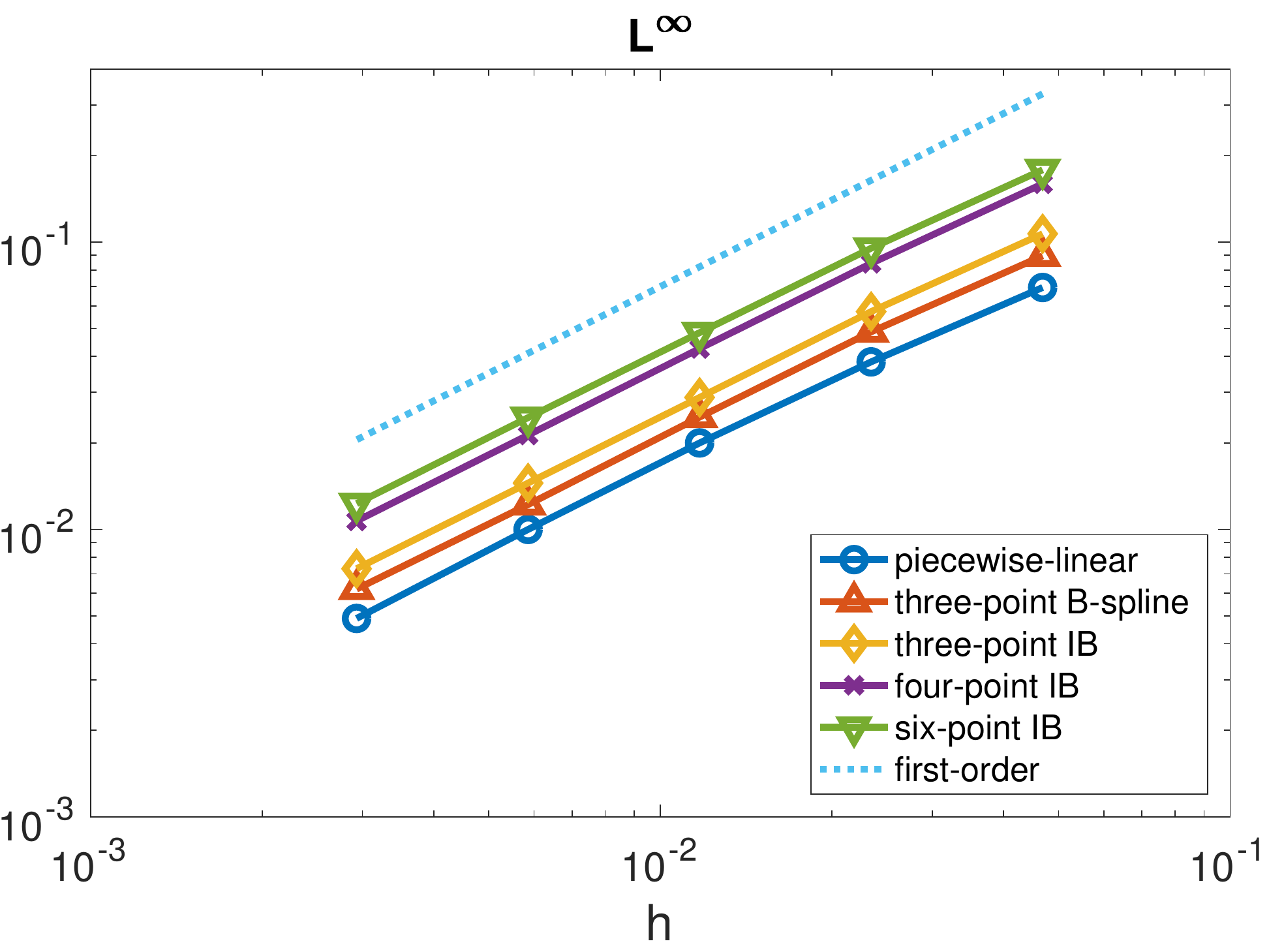}
 \caption{Representative $\log$-$\log$ plot of different error norms in velocity with respect to the finest Eulerian mesh width $h$ for various kernels with $\mfac = 2.0$. The piecewise linear kernel shows the smallest errors. Notice that first-order convergence is obtained with all choices of kernels and for all values of $\mfac$.}
 \label{fig-channel_error}
 \end{figure}

\begin{figure}[t!!]
 \centering
 \includegraphics[scale = 0.25]{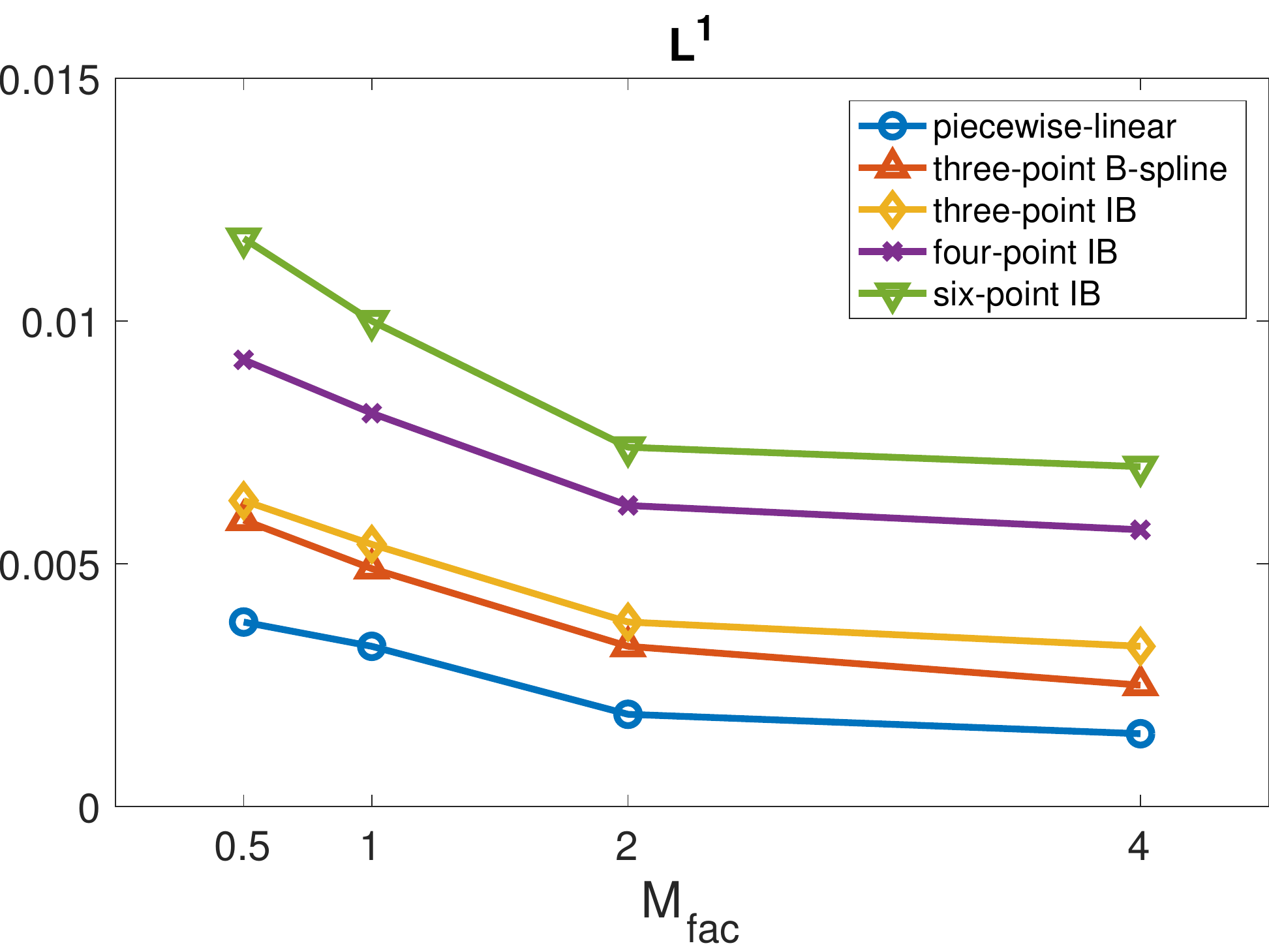}
 \includegraphics[scale = 0.25]{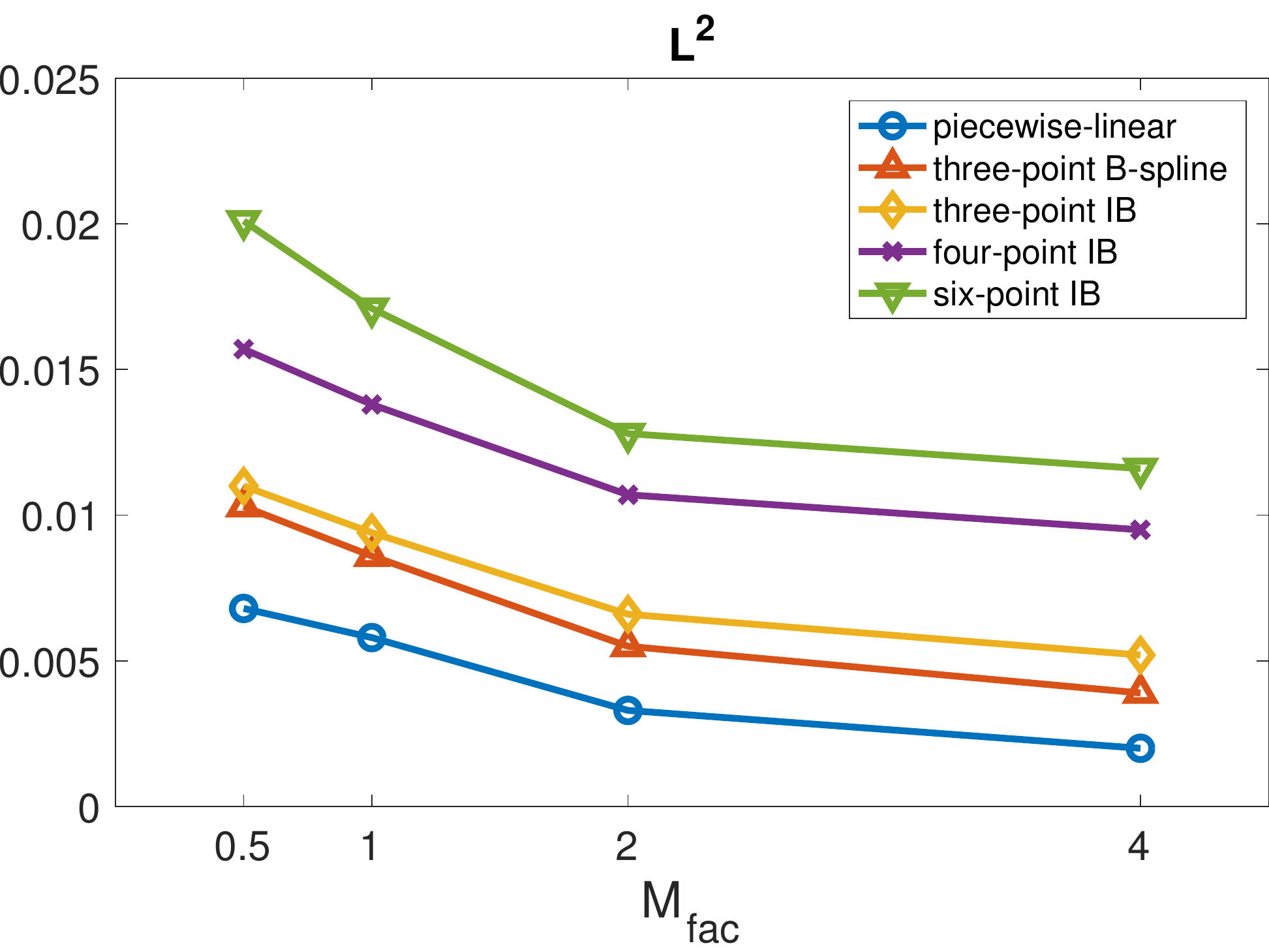}
 \includegraphics[scale = 0.25]{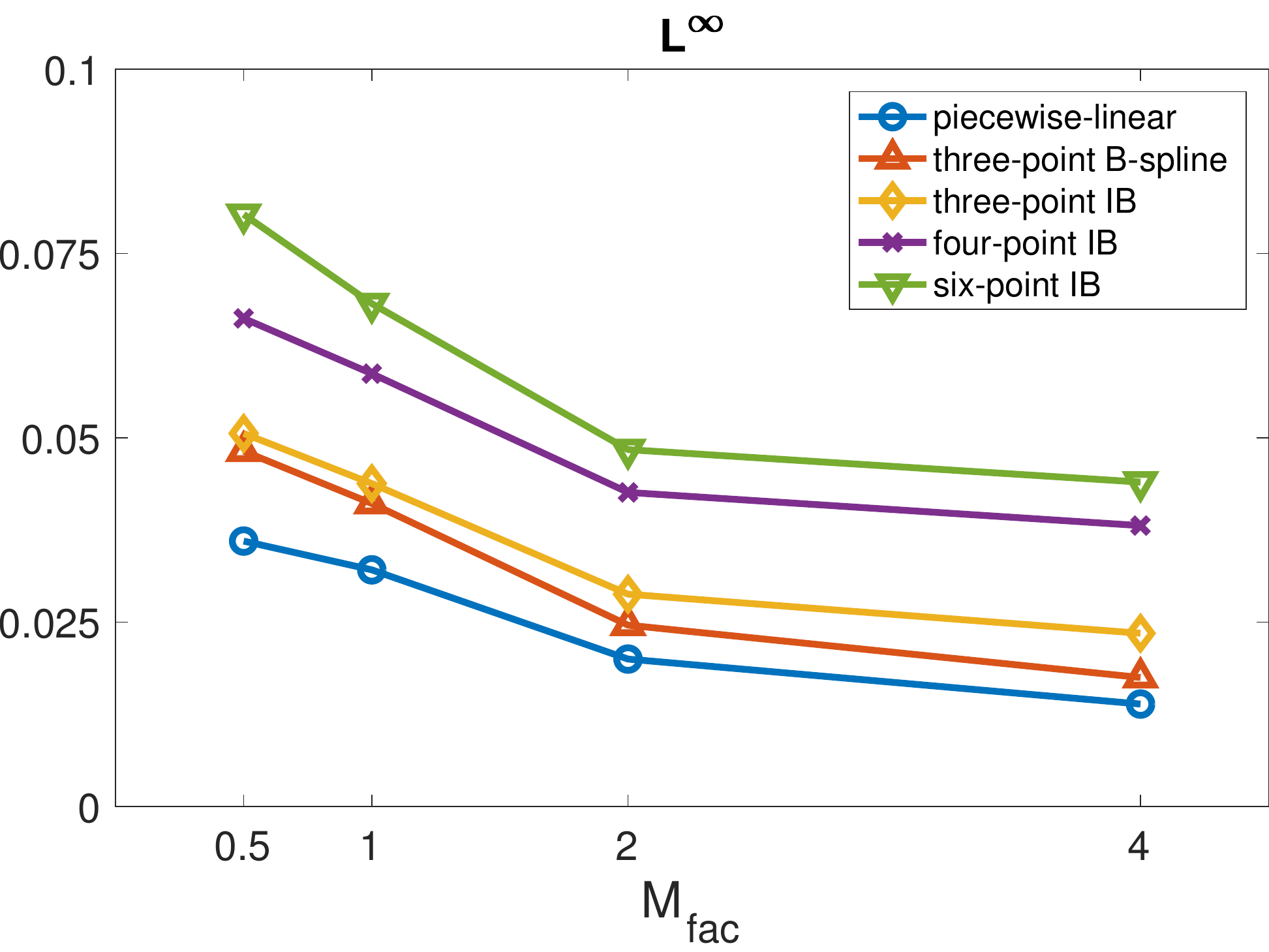}
 \caption{Plot of the error norms in velocity for values of $\mfac$ = 0.5, 1, 2, and 4 at $N = 128$ for various kernels. It shows that the structural mesh that is relatively coarser than the finest background Cartesian grids yields the lowest error. The piecewise linear kernel shows the smallest errors. Note that similar results are obtained at all resolutions with all choices of kernels.}
 \label{fig-channel_error_vs_mfac}
 \end{figure}

Figure~\ref{fig-channel_error} shows a convergence study using different error norms for representative kernels with $\mfac = 2.0$.
These results indicate that the velocity converges at first order for all kernels.
Similar convergence rates are observed for $\mfac$ = 0.5, 1.0, 2.0, and 4.0.
We find that using the piecewise linear kernel leads to the best accuracy for this test.
Figure~\ref{fig-channel_error_vs_mfac} shows the error plots in velocity for representative kernels for $\mfac$ = 0.5, 1, 2, and 4 at $N = 128$. 
In all cases, we observe the general trend that the cases with $\mfac \geq 1$, i.e., cases in which the structural mesh is coarser than the background Cartesian grids, result in better accuracy.
This finding is in agreement with the results of Section~\ref{subsec:cylinder}.
Similar results are obtained at all resolutions with all choices of kernels.
These results demonstrate that the kernels with relatively narrower support (piecewise linear, three-point IB, and three-point B-spline) and using a structural mesh that is coarser than the background Cartesian grid yield the best accuracy for simulating internal flow within a stationary geometry.
As in the tests reported in Section~\ref{subsec:cylinder}, we observe that the scheme converges under grid refinement for all choices of kernels and for all values of $\mfac$.

\subsection{Modified Turek-Hron benchmark}
\label{subsec:turek-hron}
\begin{figure}[t!!]
	\begin{center}
		\begin{subfigure}[t]{0.5\textwidth}
			\centering
			\includegraphics[scale = 0.9]{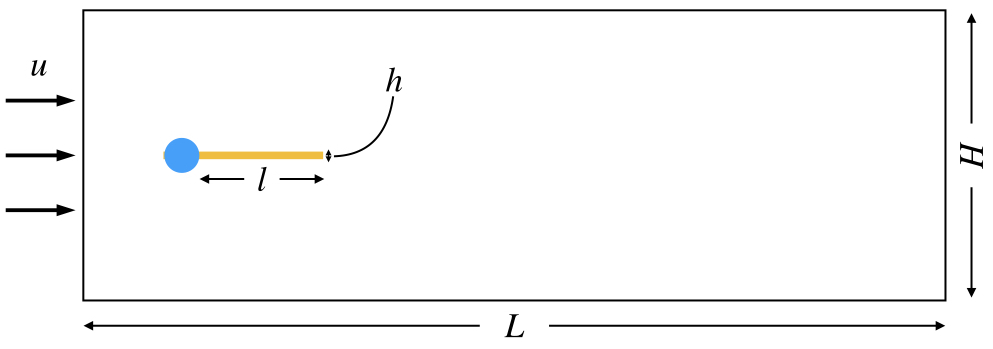}
			\caption{} \label{fig:turek_hron_schematic_full}
		\end{subfigure}
		\begin{subfigure}[t]{0.45\textwidth}
			\centering
			\includegraphics[scale = 1]{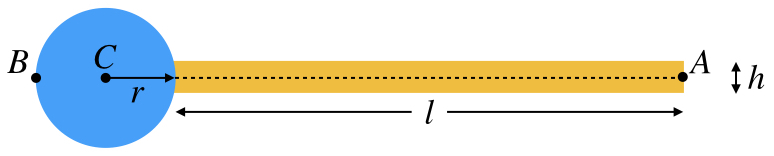}
			\caption{} \label{fig:turek_hron_schematic_zoom}
		\end{subfigure}
		\caption{(a) Schematic of the Turek-Hron benchmark~\cite{Turek2007}. (b) Detail of the immersed cylinder and flexible beam.}
		\label{fig:turek_hron_schematic}		
	\end{center}
\end{figure}

Next we consider a version of the Turek-Hron FSI benchmark of flow interacting with a flexible elastic beam mounted to a stationary circular cylinder~\cite{Turek2007}. 
Figure~\ref{fig:turek_hron_schematic} shows a schematic of the setup for this benchmark. 
In the original Turek-Hron benchmark, the domain length is $L = 2.5$ and height is $H = 0.41$.
Our modification to this benchmark uses $L = 2.46 = 6.0H$ for the domain to obtain square Cartesian grid cells, but this change is small enough that it does not affect the results substantially.
The fine-grid Cartesian cell size is $\euleriandx = H/(4N)$, and the time step size is $\Delta t = 0.001025/N$.
The circular cylinder is centered at $(0.2, 0.2)$ with radius $r = 0.05$. 
The elastic beam has length $l = 0.35$ and height $h = 0.02$. 
The left end of the beam is fixed at the back of cylinder. 
We track the position of the control point $A$ highlighted in Figure~\ref{fig:turek_hron_schematic_zoom}, whose initial position is $A(0) = (0.6, 0.2)$. 
The boundary conditions are $u(0,y) = 1.5\overline{U}\frac{y(H-y)}{\left(H/2\right)^2}$ for $x = 0$, in which $\overline{U}=2$ is the average velocity, zero normal traction and zero tangential velocity conditions for $x = L$, and zero velocity condition for $y = 0$ and $y = H$. 
The Reynolds number is $Re=\frac{\rho\overline{U}d}{\mu}=200$, in which $d=2r=0.1$ is the diameter of the cylinder, $\rho=1000$, and $\mu=1$.
Notice that without the elastic beam, this problem reduces to a version of the flow past a cylinder benchmark already considered in Section~\ref{subsec:cylinder}. 
Results reported in Section~\ref{subsec:cylinder} indicate that the three-point B-spline kernel provides the best accuracy at a given spatial resolution among the kernel functions considered in this study. 
Consequently, here we use the three-point B-spline kernel with $\mfac = 2.0$ for the cylinder in all cases, so that we can isolate the effects of the choice of kernel function and $\Mfac$ on the dynamics of the elastic beam. 
Also note that, following the specification of the test by Turek and Hron, the immersed body is positioned asymmetrically in the $y$-direction to ensure a consistent onset of beam motion across discretization and solver approaches.

In this benchmark, we use an incompressible neo-Hookean material for the elastic beam, whose strain energy functional is defined as
\begin{equation}
W_{\text{NH}} = \frac{1}{2}G_\text{s}(\bar{I}_1 - 3),
\end{equation}
in which $G_\text{s}$ is the shear modulus. 
This differs from the problem specification in Turek and Hron's original paper, which uses a compressible St.~Venant-Kirchhoff model for the elastic beam, but our numerical framework enforces incompressibility on both solid and fluid by formulation, and so we cannot readily model the elastic beam as a compressible material.
However, we also note that our results with an incompressible material model still fall within the range of results reported by Turek and Hron using a compressible material model~\cite{Turek2011}.
We do not expand on those results because that comparison is not the main focus of this study.

\begin{table}[t!!]
    \setlength\tabcolsep{4.5pt}
    \scriptsize
	\centering	
	\caption{Results for the modified Turek-Hron benchmark using the three-point B-spline kernel with various values of $\mfac$ under different grid resolutions. $N$ is the number of grid cells on coarsest grid level, $A_x$ and $A_y$ are $x$-, $y$-displacements of the point $A$, and $St_x$ and $St_y$ are Strouhal numbers for the oscillations of $A_x$ and $A_y$.}
\begin{tabular}{c c c |c c | c c| c c|l}
\cline{2-9}
& \multicolumn{2}{ |c| }{$\mfac = 0.5$} & \multicolumn{2}{ |c| }{$\mfac = 1.0$} & \multicolumn{2}{ |c| }{$\mfac = 2.0$} & \multicolumn{2}{ |c| }{$\mfac = 4.0$} \\  
\cline{1-9}
\multicolumn{1}{ |c| }{$N$} & $A_x$ $(\times 10^{-3})$ & $St_x$ & $A_x$ $(\times 10^{-3})$ & $St_x$ & $A_x$ $(\times 10^{-3})$ & $St_x$ & $A_x$ $(\times 10^{-3})$ & $St_x$ & \\ \cline{1-9}
\multicolumn{1}{ |c| }{32} & $-2.30 \pm 2.16$ & 10.4 & $-2.63 \pm 2.47$ & 10.4 & $-2.89 \pm 2.85$ & 10.4 & $-3.16 \pm 3.00$ & 10.8 &    \\ 
\cline{1-9}
\multicolumn{1}{ |c| }{64} & $-2.57 \pm 2.44$ & 10.8 & $-2.69 \pm 2.55$ & 10.8 & $-2.76 \pm 2.66$ & 10.8 & $-3.03 \pm 2.89$ & 10.8 &    \\ 
\cline{1-9}
\multicolumn{1}{ |c| }{128} & $-2.75 \pm 2.61$ & 10.8 & $-2.77 \pm 2.63$ & 10.8 & $-2.83 \pm 2.70$ & 10.8 & $-2.88 \pm 2.73$ & 10.8 &    \\ 
\cline{1-9}
\multicolumn{1}{ |c| }{256} & $-2.79 \pm 2.64$ & 10.8 & $-2.82 \pm 2.67$ & 10.8 & $-2.83 \pm 2.69$ & 10.8 & $-2.85 \pm 2.70$ & 10.8 &    \\ 
\hhline{=========}
\multicolumn{1}{ |c| }{$N$} & $A_y$ $(\times 10^{-3})$ & $St_y$ & $A_y$ $(\times 10^{-3})$ & $St_y$ & $A_y$ $(\times 10^{-3})$ & $St_y$ & $A_y$ $(\times 10^{-3})$ & $St_y$ &\\
\cline{1-9}
\multicolumn{1}{ |c| }{32} &  $1.47 \pm 30.5$ & 5.00 & $1.67 \pm 32.4$ & 5.00 & $1.56 \pm 34.9$ & 5.00 & $1.23 \pm 36.1$ & 5.42 & \\ 
\cline{1-9}
\multicolumn{1}{ |c| }{64} &  $1.41 \pm 32.7$ & 5.00 &  $1.44 \pm 33.4$ & 5.00 & $1.41 \pm 34.2$ & 5.00 & $1.49 \pm 35.3$ & 5.00 &\\ 
\cline{1-9}
\multicolumn{1}{ |c| }{128} &  $1.42 \pm 33.9$ & 5.00 & $1.43 \pm 34.0$ & 5.00 & $1.44 \pm 34.5$ & 5.00 & $1.42 \pm 34.7$ & 5.00 &    \\ 
\cline{1-9}
\multicolumn{1}{ |c| }{256} &  $1.42 \pm 34.2$ & 5.00 & $1.43 \pm 34.3$ & 5.00 & $1.42 \pm 34.4$ & 5.00 & $1.42 \pm 34.5$ & 5.00 &    \\ 
\cline{1-9}
\end{tabular}
	\label{table:nh_kinematics_grid_refinement}	
\end{table} 

\begin{figure}[t!!]
\centering
\includegraphics[scale = 0.365,trim={10 30 0 10},clip]{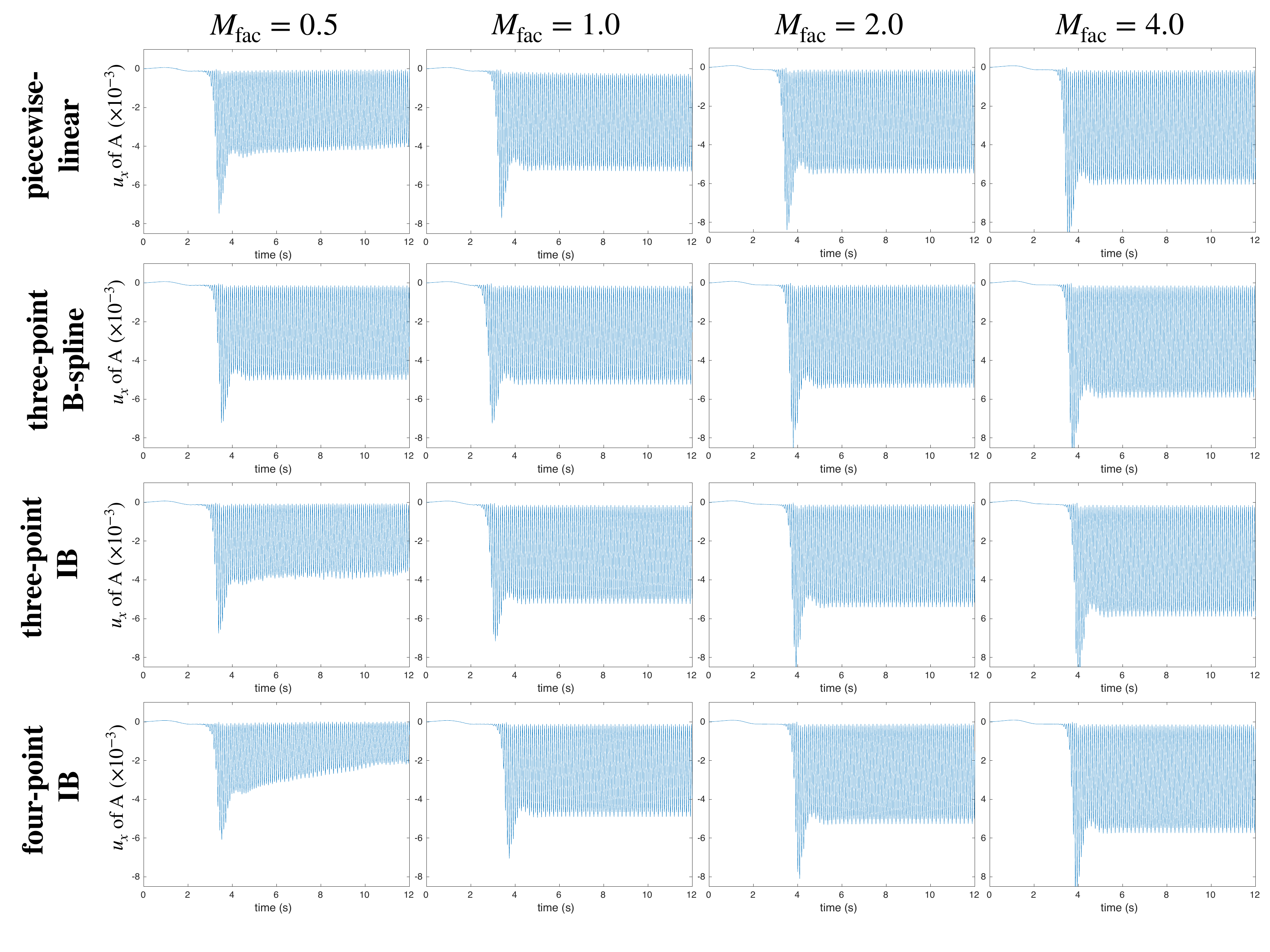} 
\includegraphics[scale = 0.365,trim={445 30 340 10},clip]{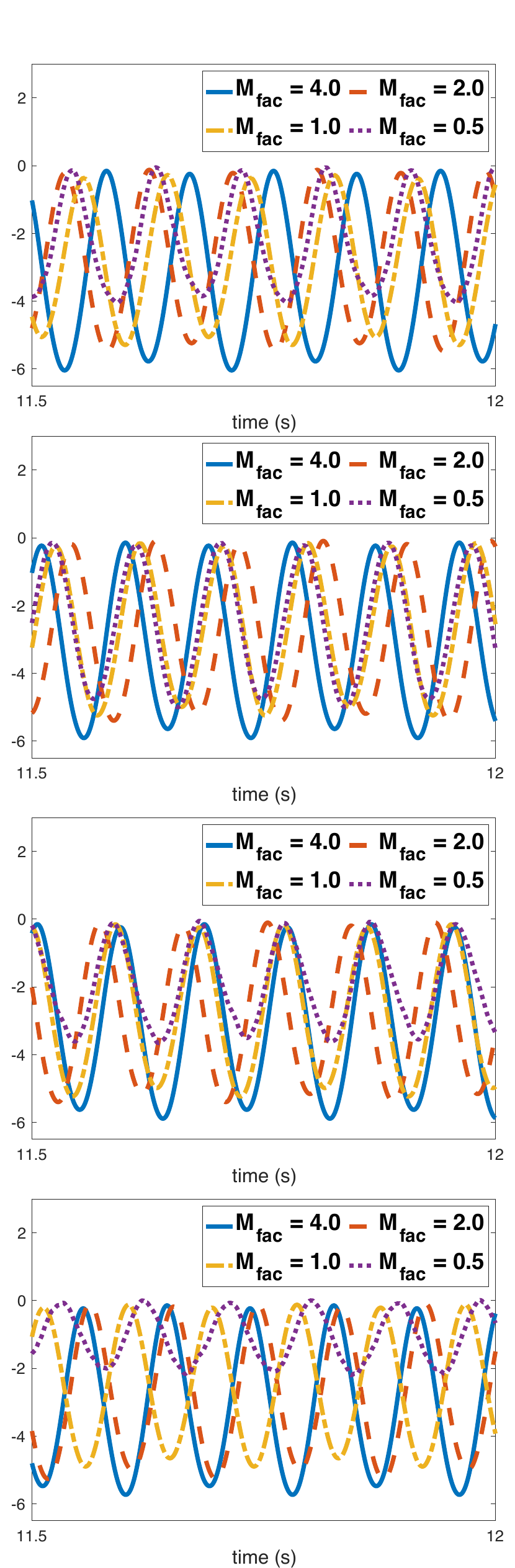} 
\caption{$x$-displacement ($A_x$) of the point $A$ for different values of $\mfac$ for the modified Turek-Hron benchmark using different kernels at a Cartesian resolution of $N = 64$. Panels in the rightmost column show the periodic oscillations between $t=11.5$ and $t=12$.}
\label{fig:A_x_comparison_N=64}
\end{figure}

\begin{figure}[t!!]
\centering
\includegraphics[scale = 0.365,trim={10 30 0 10},clip]{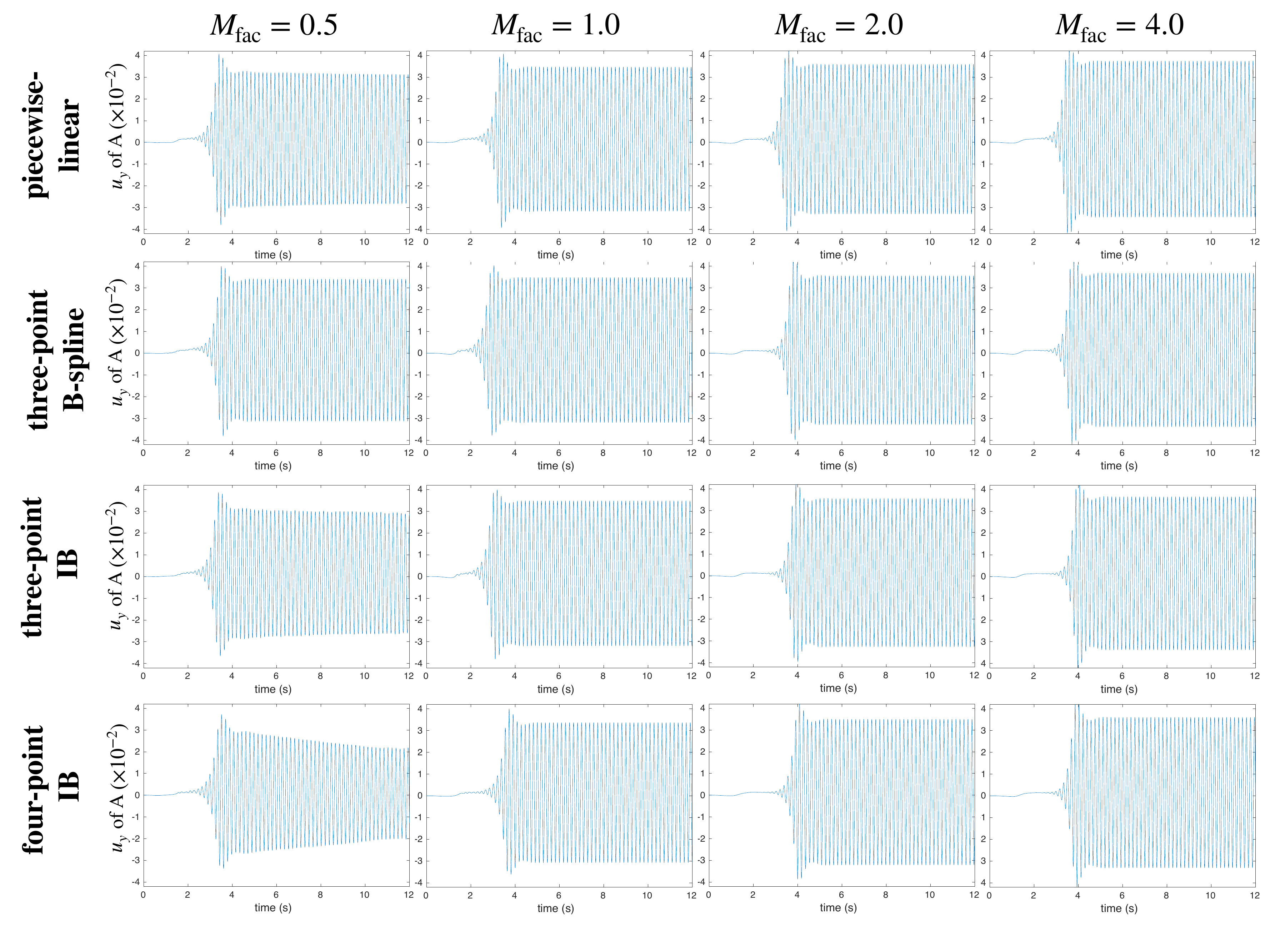} 
\includegraphics[scale = 0.365,trim={445 30 340 10},clip]{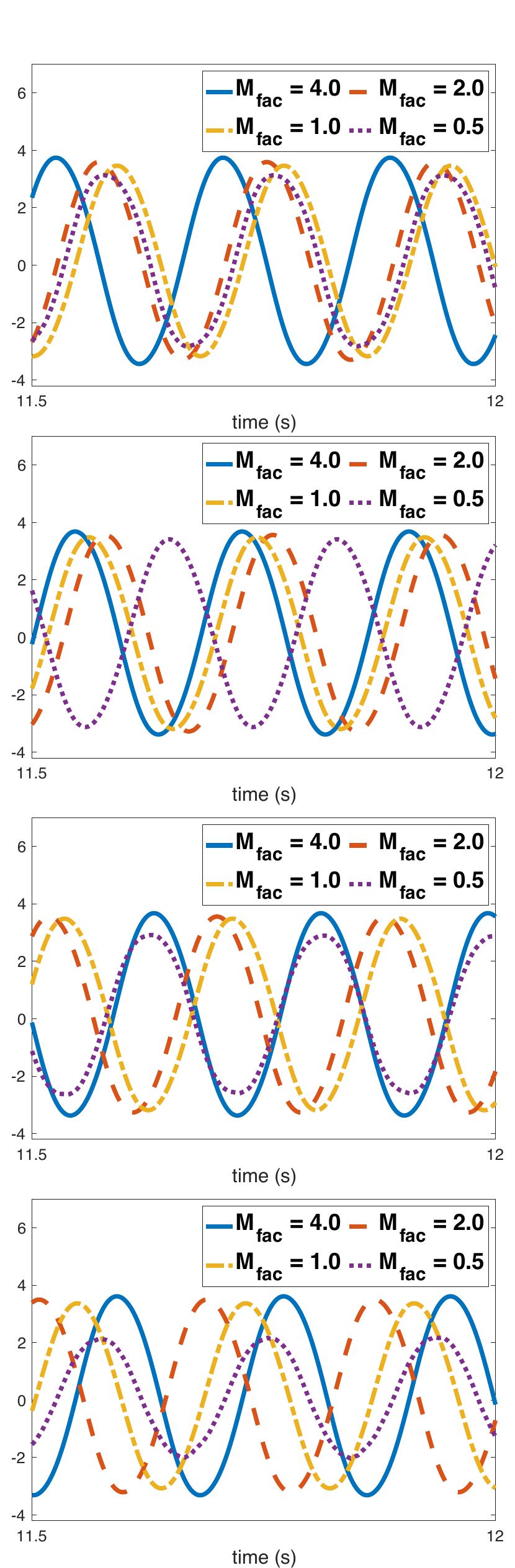} 
\caption{$y$-displacement ($A_y$) of the point $A$ for different values of $\mfac$ for the modified Turek-Hron benchmark using different kernels at a Cartesian resolution of $N = 64$. Figures in the rightmost column show the periodic oscillations between $t=11.5$ and $t=12$.}
\label{fig:A_y_comparison_N=64}
\end{figure}

Table~\ref{table:nh_kinematics_grid_refinement} shows comparisons for the three-point B-spline kernel for $\mfac = 0.5, 1, 2,$ and $4$ under grid refinement.
It reports the average and the amplitude of the $x$- and $y$-displacements ($A_x$ and $A_y$) of the point $A$, as well as the Strouhal numbers ($St_x$ and $St_y$) to quantify the oscillations of $A_x$ and $A_y$.
We obtain comparable results under grid refinement, which are more consistent between different $\mfac$ values as we refine the resolution. 
Similar to the results from other benchmarks, this benchmark also indicates that under grid refinement, the results become independent of $\mfac$ and the type of kernel. 
We again focus on the effect of $\mfac$ and the choice of kernel function at an intermediate spatial resolution. 
Figures~\ref{fig:A_x_comparison_N=64} and~\ref{fig:A_y_comparison_N=64} compare representative kernels at $N=64$ for $\mfac$ = 0.5, 1, 2, and 4.
The three-point B-spline kernel clearly yields more consistent results for different values of $\mfac$ at this resolution. 
Table~\ref{table:nh_kinematics_comparison} summarizes the differences between selected kernels in Figures~\ref{fig:A_x_comparison_N=64} and~\ref{fig:A_y_comparison_N=64}.
Appendix~\ref{sec:TH_appendix} provides the results for the remaining IB and B-spline kernels (see Figures~\ref{fig:A_x_comparison_N=64_appendix} and~\ref{fig:A_y_comparison_N=64_appendix}), and Tables~\ref{table:nh_kinematics_comparison} and~\ref{table:nh_kinematics_N=64} show that the piecewise linear and three-~and four-point IB kernels yield displacements that show large discrepancies from the converged displacements if we set $\mfac = 0.5$. 
The five-~and six-point IB and three-, four-, five-, and six-point B-spline kernels produce displacements that are relatively consistent.
However, Table~\ref{table:nh_kinematics_grid_refinement} yields that the Strouhal numbers converge to 10.8, and only the three-point B-spline kernel shows the converged value for Strouhal number consistently for all values of $\mfac$ = 0.5, 1, 2, and 4 at the intermediate Cartesian resolution of $N=64$.
These results indicate that the three-point B-spline kernel is less sensitive to changes in $\mfac$, whereas other kernels show clear loss of accuracy as we refine the Lagrangian mesh for a fixed Eulerian grid that is of intermediate spatial resolution.

\begin{table}[t!!]
    \setlength\tabcolsep{4.5pt}
    \scriptsize
	\centering	
	\caption{Results for the modified Turek-Hron benchmark using kernels in Figures~\ref{fig:A_x_comparison_N=64} and~\ref{fig:A_y_comparison_N=64}. The number of Cartesian grid cells on the coarsest level is $N = 64$, $A_x$ and $A_y$ are $x$-, $y$-displacements of the point $A$, and $St_x$ and $St_y$ are Strouhal numbers for the oscillations of $A_x$ and $A_y$.}
\begin{tabular}{c c c |c c | c c| c c|l}
\cline{2-9}
& \multicolumn{2}{ |c| }{$\mfac = 0.5$} & \multicolumn{2}{ |c| }{$\mfac = 1.0$} & \multicolumn{2}{ |c| }{$\mfac = 2.0$} & \multicolumn{2}{ |c| }{$\mfac = 4.0$} \\  
\cline{1-9}
\multicolumn{1}{ |c| }{Kernel} & $A_x$ $(\times 10^{-3})$ & $St_x$ & $A_x$ $(\times 10^{-3})$ & $St_x$ & $A_x$ $(\times 10^{-3})$ & $St_x$ & $A_x$ $(\times 10^{-3})$ & $St_x$ & \\ \cline{1-9}
\multicolumn{1}{ |c| }{Piecewise-linear} & $-2.20 \pm 2.14$ & 10.4 & $-2.74 \pm 2.53$ & 10.8 & $-2.80 \pm 2.68$ & 10.4 & $-3.10 \pm 2.95$ & 10.8 &    \\ 
\cline{1-9}
\multicolumn{1}{ |c| }{B-spline (3-point)} & $-2.57 \pm 2.44$ & 10.8 & $-2.69 \pm 2.55$ & 10.8 & $-2.76 \pm 2.66$ & 10.8 & $-3.03 \pm 2.89$ & 10.8 &    \\ 
\cline{1-9}
\multicolumn{1}{ |c| }{IB (3-point)} & $-2.03 \pm 1.97$ & 10.4 & $-2.69 \pm 2.54$ & 10.4 & $-2.76 \pm 2.65$ & 10.8 & $-3.02 \pm 2.87$ & 10.8 &    \\ 
\cline{1-9}
\multicolumn{1}{ |c| }{IB (4-point)} & $-1.55 \pm 1.56$ & 10.4 & $-2.51 \pm 2.39$ & 10.4 & $-2.69 \pm 2.58$ & 10.8 & $-2.94 \pm 2.80$ & 10.8 &    \\ 
\hhline{=========}
\multicolumn{1}{ |c| }{Kernel} & $A_y$ $(\times 10^{-3})$ & $St_y$ & $A_y$ $(\times 10^{-3})$ & $St_y$ & $A_y$ $(\times 10^{-3})$ & $St_y$ & $A_y$ $(\times 10^{-3})$ & $St_y$ &\\
\cline{1-9}
\multicolumn{1}{ |c| }{Piecewise-linear} & $1.38 \pm 30.8$ & 5.00 & $1.46 \pm 33.2$ & 5.00 & $1.45 \pm 34.4$ & 5.00 & $1.47 \pm 35.9$ & 5.00 &    \\ 
\cline{1-9}
\multicolumn{1}{ |c| }{B-spline (3-point)} & $1.41 \pm 32.7$ & 5.00 & $1.44 \pm 33.4$ & 5.00 & $1.41 \pm 34.2$ & 5.00 & $1.49 \pm 35.3$ & 5.00 &    \\ 
\cline{1-9}
\multicolumn{1}{ |c| }{IB (3-point)} & $1.37 \pm 29.2$ & 5.00 & $1.45 \pm 33.3$ & 5.00 & $1.42 \pm 34.1$ & 5.00 & $1.48 \pm 35.2$ & 5.00 &    \\ 
\cline{1-9}
\multicolumn{1}{ |c| }{IB (4-point)} & $1.04 \pm 25.9$ & 5.00 & $1.46 \pm 32.2$ & 5.00 & $1.41 \pm 33.5$ & 5.00 & $1.47 \pm 34.6$ & 5.00 &    \\ 
\cline{1-9}
\end{tabular}
	\label{table:nh_kinematics_comparison}	
\end{table} 

\subsection{Two-dimensional pressure-loaded elastic band}
\label{subsec:pressurized_band}
\begin{figure}[t!!]
	\begin{center}
	\includegraphics[scale = 0.15,trim={0 120 0 120},clip]{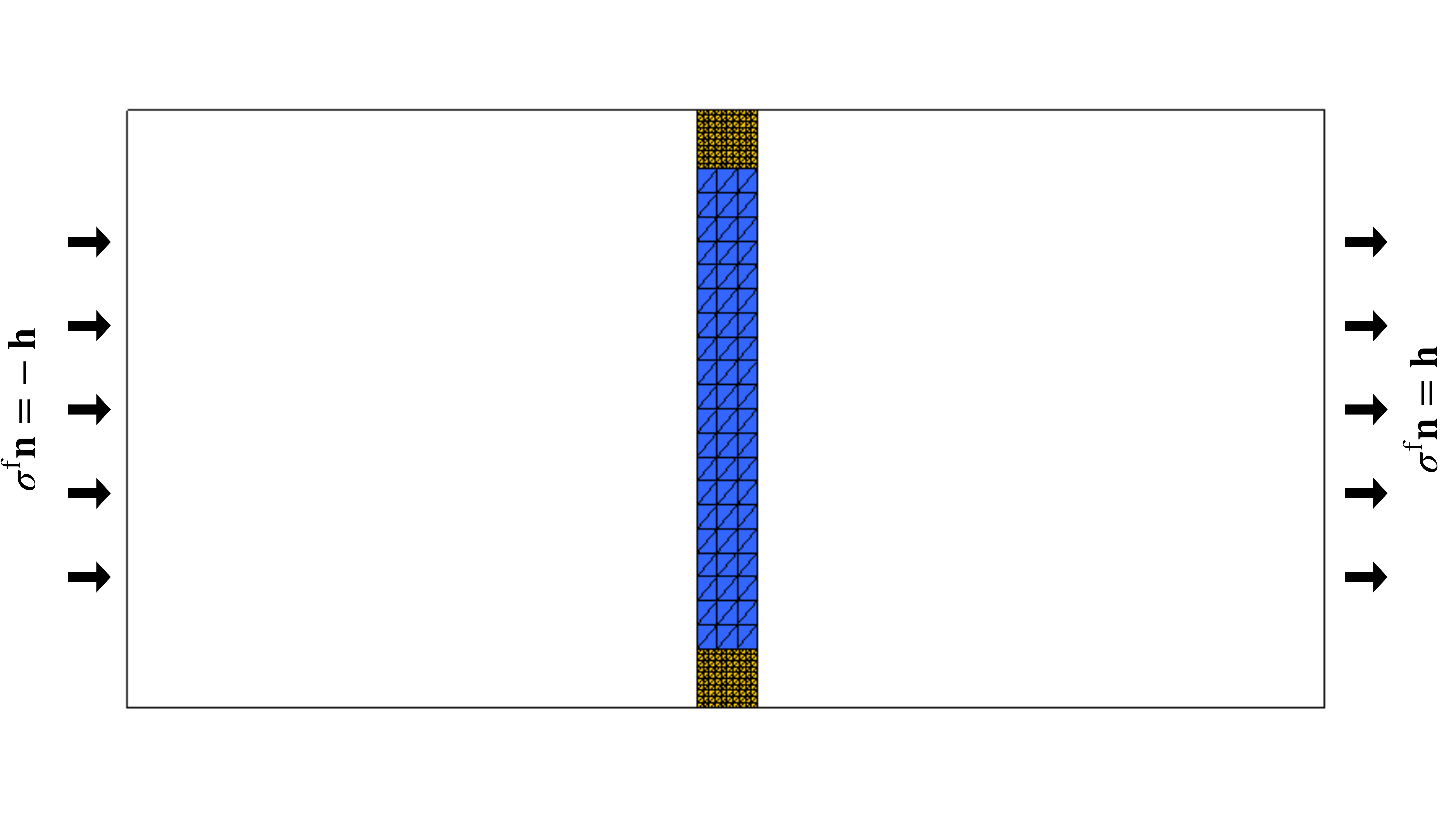}
	\caption{Schematic of two-dimensional pressure-loaded elastic band adopted from Vadala-Roth et al.~\cite{Vadala-Roth2020}. The loading on the band (blue) is driven by fluid forces induced by the pressure gradient between the left and right boundaries of the computational domain. The effective shape of the kernel function changes near the boundary of the computational domain. So we avoid issues that may arise from using a finer structural mesh ($\mfac = 0.5$) for the two rigid blocks (yellow) by which the top and bottom of the band are fixed in place. In this figure, $\mfac = 2$ for the band away from the boundary.}
	\label{fig:pressurized_band_schematic}		
	\end{center}
\end{figure}

\begin{figure}[t!!]
	\begin{center}
		\hspace{-0.1in}
		\begin{subfigure}[t]{0.33\textwidth}
			\centering
			\includegraphics[scale = 0.205]{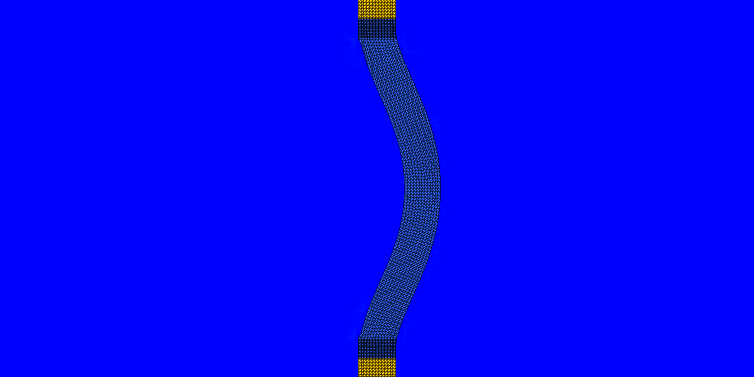} 
			\caption{$\mfac = 0.5$} \label{fig:pb_u_mfac_0_5}
		\end{subfigure}
		\begin{subfigure}[t]{0.33\textwidth}
			\centering
			\includegraphics[scale = 0.205]{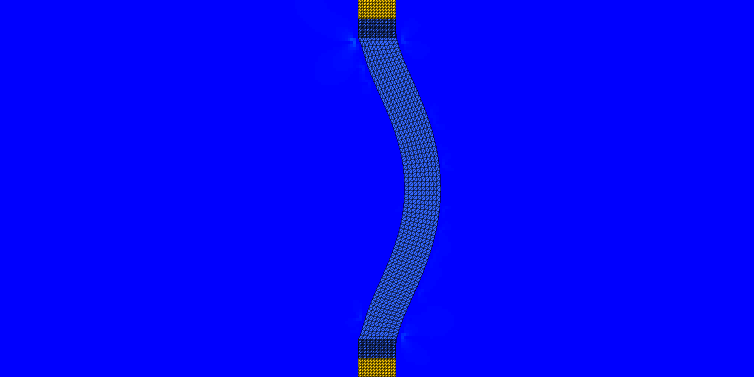} 
			\caption{$\mfac = 0.75$} \label{fig:pb_u_mfac_0_75}
		\end{subfigure} 
		\begin{subfigure}[t]{0.33\textwidth}
			\centering
			\includegraphics[scale = 0.205]{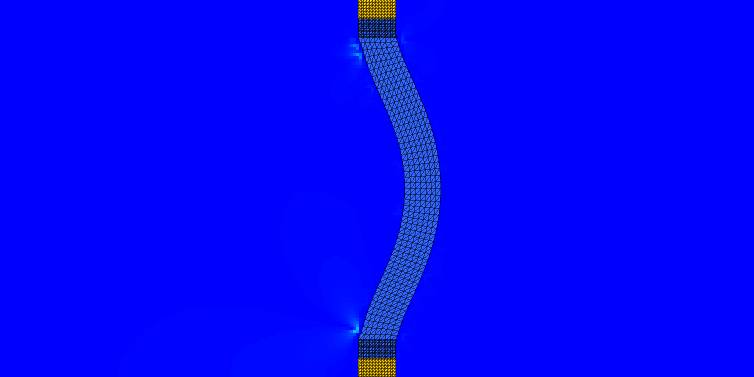} 
			\caption{$\mfac = 1$} \label{fig:pb_u_mfac_1}	
		\end{subfigure}\\
		\hspace{-0.2in}
		\includegraphics[scale = 0.11]{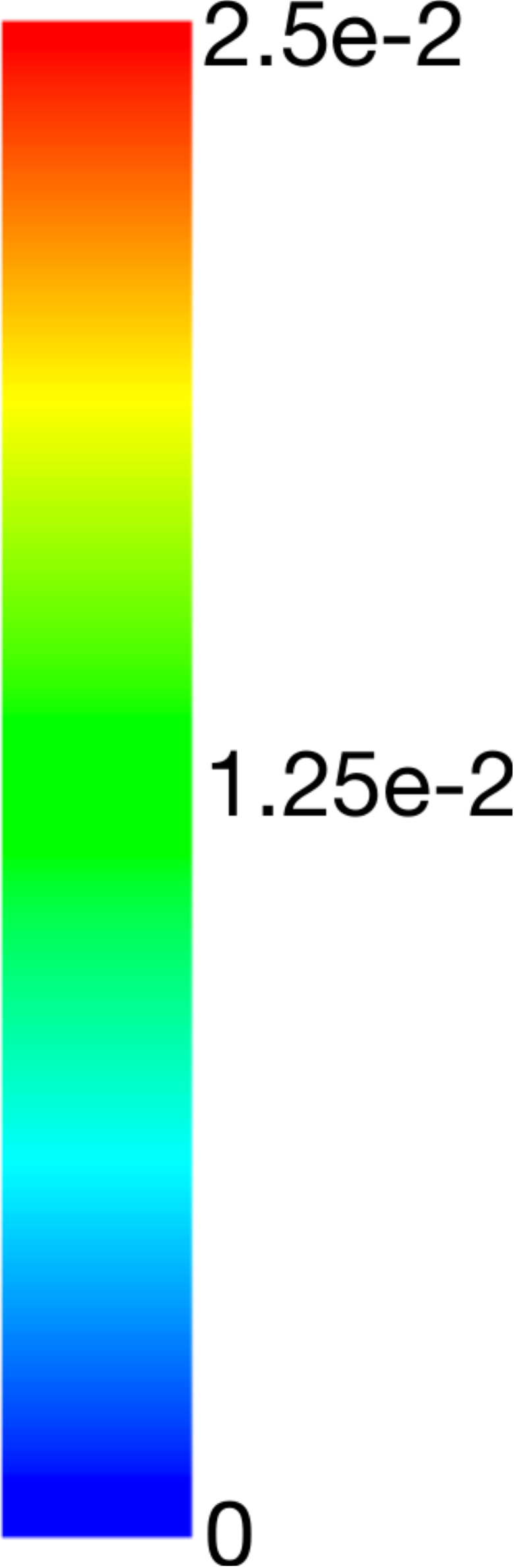}
		\hspace{0.25in}
		\begin{subfigure}[t]{0.33\textwidth}
			\centering
			\includegraphics[scale = 0.205]{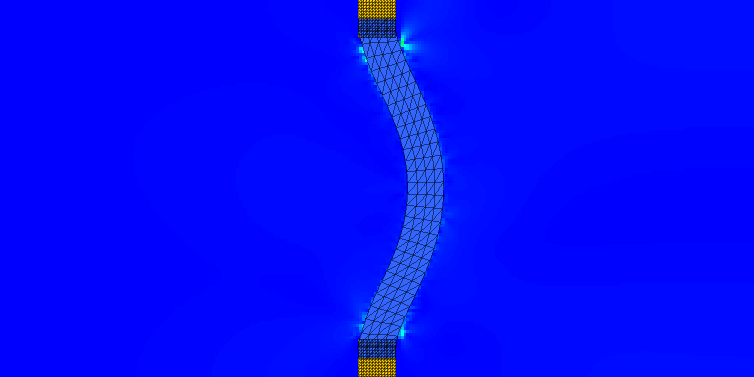} 
			\caption{$\mfac = 2$} \label{fig:pb_u_mfac_2}
		\end{subfigure}
		\begin{subfigure}[t]{0.33\textwidth}
			\centering
			\includegraphics[scale = 0.205]{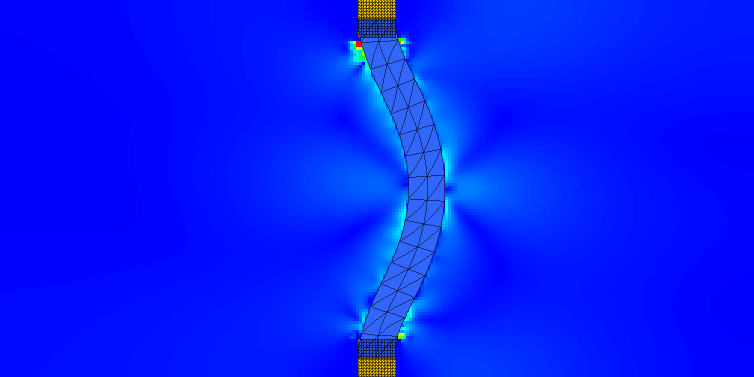} 
			\caption{$\mfac = 4$} \label{fig:pb_u_mfac_4}	
		\end{subfigure}
		\caption{Comparison of velocity fields from the pressure-loaded two-dimensional elastic band for $\mfac$ = 0.5, 0.75, 1, 2, and 4. The simulations use an $N \times N$ 	grid with $N = 128$. The three-point B-spline kernel is used for this figure, but we observe the similar results with other kernels. If the structural mesh is relatively coarser ($\mfac > 1$) than the finest background Cartesian grids, then we obtain low accuracy for simulating pressurized elastic band.}
		\label{fig:pressurized_band}		
	\end{center}
\end{figure}

Results reported in Sections~\ref{subsec:cylinder}, \ref{subsec:channel}, and \ref{subsec:turek-hron} suggest that larger $\mfac$ values generally give higher accuracy at a fixed Cartesian grid resolution, independent of the choice of kernel function.
The tests considered so far, however, are examples of shear-dominant flows.
Here we consider cases in which pressure loading dominates, as commonly encountered in biological and biomedical applications. 
To do so, we use a pressure-loaded ``elastic band'' model (Figure~\ref{fig:pressurized_band_schematic}) that is adopted from Vadala-Roth et al.~\cite{Vadala-Roth2020}. 
This uses an incompressible neo-Hookean material model, as described in Section~\ref{subsec:turek-hron}, with the shear modulus $G_\text{s} = 200$.
We set $\rho = 1.0$ and $\mu = 0.01$. 
The computational domain is $2L \times L$, in which $L = 1$. 
The simulations use a uniform grid with an $2N \times N$ grid with $N = 128$.
The fine-grid Cartesian cell size is $\euleriandx = L/N$, and the time step size is $\Delta t = 0.001/N$.
Fluid tractions $\vec{\tau}(\x,t) = \vec{\bbsigma} (\x, t)\n(\x) = -\vec{h}$ and $\vec{\tau}(\x,t) = \vec{\bbsigma} (\x, t)\n(\x) = \vec{h}$ are imposed on the left and right boundaries of the computational domain, in which $\bbsigma = -p {\mathbb I} + \mu\left(\grad \vec{u} + \grad \vec{u}^T\right)$ is the fluid stress tensor and $\vec{h}=(5,0)$, and zero velocity is enforced along the top and bottom boundaries.
The elastic band deforms and ultimately reaches a steady-state configuration determined by the pressure difference across the band. 
We use a grid resolution ($N=128$) that is fine enough so that the elastic bands are well-resolved for all cases to isolate the effect of $\mfac$ on the Lagrangian-Eulerian coupling and eliminate the effect of elastic response from the band. 
Moreover, the effective shape of the kernel function changes near the boundary of the computational domain.
We attempt to isolate the elastic model from issues that may arise at or near the physical boundaries by attaching it to the boundary through rigid blocks (of height $h=0.1$) that are discretized using relatively fine structural meshes ($\mfac = 0.5$).
In the continuous problem, there is no flow at equilibrium, but Figure~\ref{fig:pressurized_band} demonstrates that if the structural mesh is coarser than the finest background Cartesian grid spacing ($\mfac > 1$), spurious flows clearly develop near the fluid-structure interface.
Table~\ref{table:elastic_band} confirms that the error is an order of magnitude larger when $\mfac$ is increased from 0.5 to $\mfac > 1$, and up to 65 times larger when increased to $\mfac = 4$.
Similar results are obtained for all of the kernels considered here.

\begin{table}[t!!]
    \setlength\tabcolsep{4.5pt}
    \scriptsize
	\centering	
	\caption{Quantification of errors in velocity fields from the pressure-loaded two-dimensional elastic band using different regularized delta functions and relative structural grid refinement ($\mfac$ = 0.5, 0.75, 1, 2, and 4). If the structural mesh is relatively coarser ($\mfac > 1$) than the finest background Cartesian grids, then we obtain low accuracy for simulating pressurized elastic band.}
\begin{tabular}{c c c |c c | c c| c c| c c|l}
\cline{2-11}
& \multicolumn{2}{ |c| }{$\mfac = 0.5$} & \multicolumn{2}{ |c| }{$\mfac = 0.75$} & \multicolumn{2}{ |c| }{$\mfac = 1.0$} & \multicolumn{2}{ |c| }{$\mfac = 2.0$} & \multicolumn{2}{ |c| }{$\mfac = 4.0$} \\  
\cline{1-11}
\multicolumn{1}{ |c| }{Kernel} & \parbox[t]{1.0cm}{\centering $L^2$\\$(\times 10^{-4})$} & \parbox[t]{1.0cm}{\centering $L^\infty$\\$(\times 10^{-3})$} & \parbox[t]{1.0cm}{\centering $L^2$\\$(\times 10^{-4})$} & \parbox[t]{1.0cm}{\centering $L^\infty$\\$(\times 10^{-3})$} & \parbox[t]{1.0cm}{\centering $L^2$\\$(\times 10^{-4})$} & \parbox[t]{1.0cm}{\centering $L^\infty$\\$(\times 10^{-3})$} & \parbox[t]{1.0cm}{\centering $L^2$\\$(\times 10^{-4})$} & \parbox[t]{1.0cm}{\centering $L^\infty$\\$(\times 10^{-3})$} & \parbox[t]{1.0cm}{\centering $L^2$\\$(\times 10^{-4})$} & \parbox[t]{1.0cm}{\centering $L^\infty$\\$(\times 10^{-3})$} & \\ 
\cline{1-11}
\multicolumn{1}{ |c| }{IB (3-point)} & 0.36 & 1.26 & 0.88 & 2.48 & 1.49 & 3.88 & 6.27 & 9.77 & 17.01 & 42.25 &    \\ 
\cline{1-11}
\multicolumn{1}{ |c| }{IB (4-point)} & 0.52 & 1.26 & 1.20 & 4.44 & 0.95 & 2.41 & 5.99 & 12.27 & 21.66 & 19.19 &    \\ 
\cline{1-11}
\multicolumn{1}{ |c| }{IB (5-point)} & 1.06 & 2.12 & 1.92 & 4.18 & 1.80 & 5.73 & 5.91 & 10.8 & 20.45 & 18.63 &    \\ 
\cline{1-11}
\multicolumn{1}{ |c| }{IB (6-point)} & 2.33 & 2.81 & 3.19 & 3.76 & 3.94 & 5.96 & 6.19 & 13.69 & 23.92 & 19.79 &    \\ 
\cline{1-11}
\multicolumn{1}{ |c| }{B-spline (3-point)} & 0.33 & 0.61 & 0.65 & 2.05 & 1.61 & 5.41 & 6.53 & 10.62 & 16.93 & 39.40 &    \\
\cline{1-11}
\multicolumn{1}{ |c| }{B-spline (4-point)} & 1.26 & 2.86 & 1.85 & 4.62 & 2.10 & 8.08 & 6.42 & 8.44 & 15.83 & 9.59 &    \\ 
\cline{1-11}
\multicolumn{1}{ |c| }{B-spline (5-point)} & 1.81 & 4.47 & 2.50 & 4.63 & 2.56 & 9.06 & 6.90 & 11.04 & 18.15 & 10.51 &    \\ 
\cline{1-11}
\multicolumn{1}{ |c| }{B-spline (6-point)} & 1.74 & 3.09 & 2.48 & 4.46 & 2.29 & 7.77 & 5.97 & 9.37 & 20.21 & 17.39 &    \\ 
\cline{1-11}
\end{tabular}
	\label{table:elastic_band}	
\end{table} 

\subsection{Bioprosthetic heart valve dynamics in a pulse duplicator}
\label{subsec:BHV}
We aim to verify our key findings in a large-scale FSI model.
To do so, we consider a dynamic model of a bovine pericaridal bioprosthetic heart valve (BHV) in a pulse duplicator, as described in detail by Lee et al.~\cite{Lee2020, LeeJTCVS}. 
The simulation setup includes a detailed IFED model of the aortic test section of an experimental pulse duplicator, and the simulation includes both substantial pressure-loads (during diastole, when the valve is closed) and shear-dominant flows (during systole, when the valve is open). 
The bovine pericardial valve leaflets are described by a modified version~\cite{Lee2020} of the Holzapfel--Gasser--Ogden (HGO) model~\cite{Gasser2006},
\begin{equation}
	W_{\text{BHV}} = C_{10}\{\exp{\left[C_{01}(\bar{I}_1-3)\right]}-1\} + \frac{k_1}{2k_2}\{\exp{\left[k_2(\kappa\bar{I}_1 + (1-3\kappa)\bar{I}_{4}^{\star}-1)^2\right]-1}\},
\end{equation}
in which $\bar{I}_{4}^{\star} = \max(\bar{I}_{4}, 1) = \max(\e_{0}^T \bar{\CC}\e_{0}, 1)$, and $\e_{0}$ is a unit vector aligned with the mean fiber direction in the reference configuration. The parameter $\kappa \in [0,\frac{1}{3}]$ describes collagen fiber angle dispersion. In the our simulations, we use $C_{10} = 0.119$~kPa, $C_{01} = 22.59$, $k_1 = 2.38$~MPa, $k_2 = 149.8$, and $\kappa = 0.292$~\cite{Lee2020}.
We use $\rho = 1.0$ g/cm$^3$ and $\mu = 1.0$ cP, and we can calculate the peak Reynolds number~\cite{Lee2020, LeeJTCVS}, $Re_\text{peak}=\frac{\rho Q_\text{peak}D}{\mu A}\approx 14800$, in which $D=25$~mm and $A$ are the geometrical diameter and cross-sectional area of the valve, respectively.
The computational domain is 5.05 cm $\times$ 10.1 cm $\times$ 5.05 cm. 
The simulations use a three-level locally refined grid with a refinement ratio of two between levels and an $N/2 \times N \times N/2$ coarse grid with $N = 64$, which yields a fine-grid Cartesian resolution of 0.4~mm. 
Here, we use the piecewise linear kernel for the test section and consider the effects of different choices of kernel functions for the valve leaflets.
Three-element Windkessel (R--C--R) models establish the upstream driving and downstream loading conditions for the aortic test section. 
A combination of normal traction and zero tangential velocity boundary conditions are used at the inlet and outlet to couple the reduced-order models to the detailed description of the flow within the test section. 
Solid wall boundary conditions are imposed on the remaining boundaries of the computational domain.
See Lee et al.~\cite{Lee2020} for further details.

\begin{figure}[t!!]
	\begin{center}
		\begin{subfigure}[t]{0.45\textwidth}
			\centering
			\includegraphics[scale = 0.125]{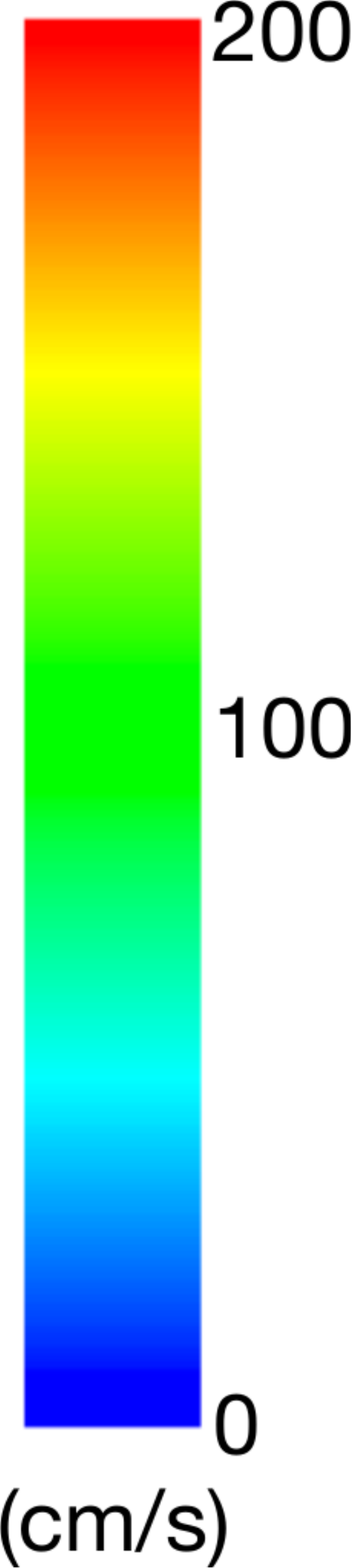}
			\includegraphics[scale = 0.2]{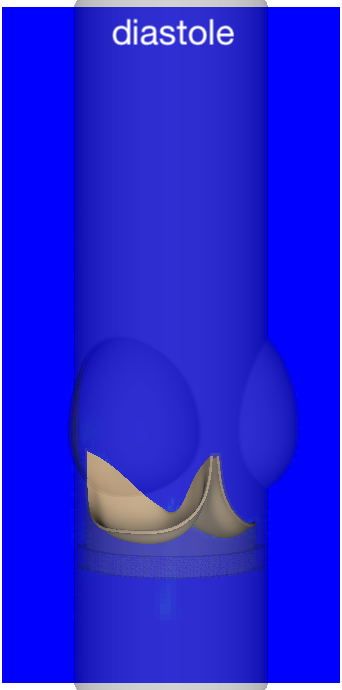} 
			\includegraphics[scale = 0.2]{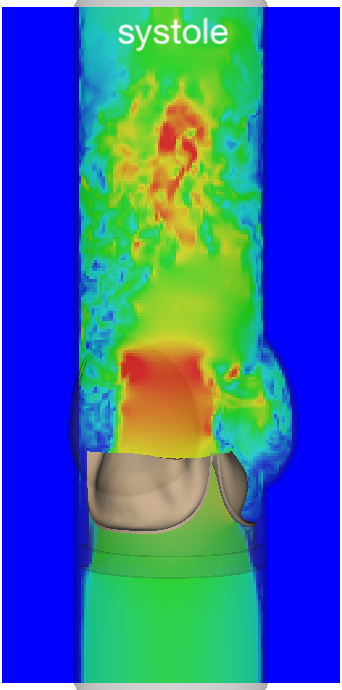} 
			\caption{$\|\u\|$ $(\mfac = 0.75)$} \label{fig:bhv_mfac_0_75}
		\end{subfigure}\hspace{-0.4in}
		\begin{subfigure}[t]{0.45\textwidth}
			\centering
			\includegraphics[scale = 0.2]{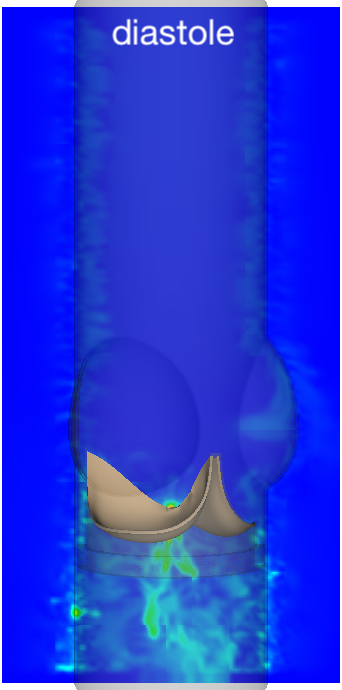} 
			\includegraphics[scale = 0.2]{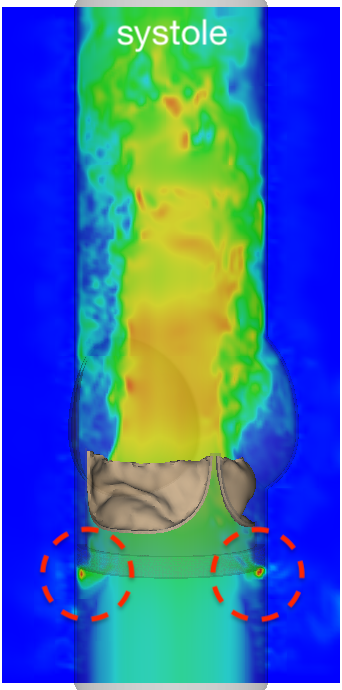} 
			\caption{$\|\u\|$ $(\mfac = 1.5)$} \label{fig:bhv_mfac_1_5}
		\end{subfigure}\\
		\begin{subfigure}[t]{0.33\textwidth}
			\centering
			\includegraphics[scale = 0.265]{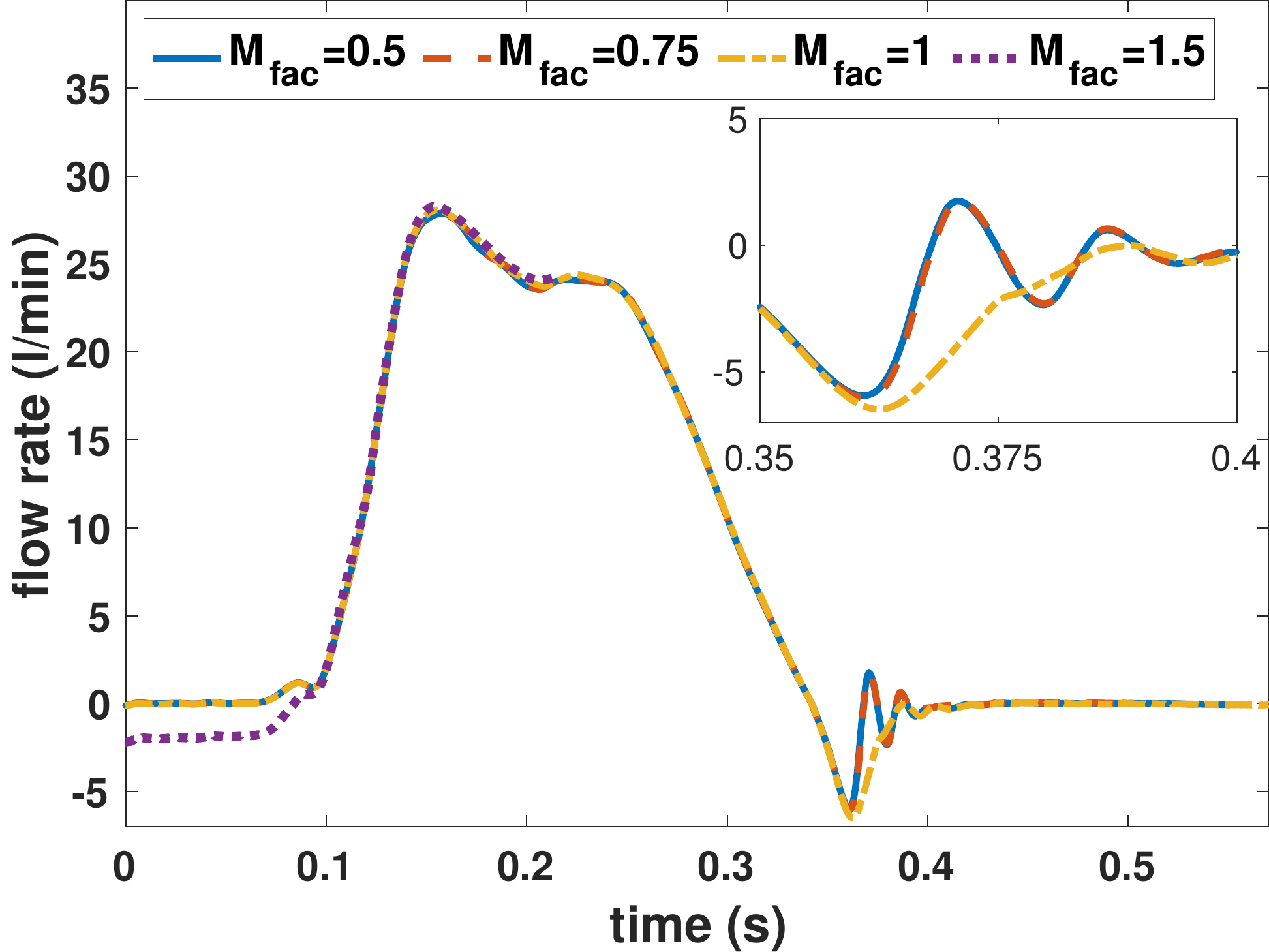} 
			\caption{$Q_\text{Ao}$} \label{fig:q_comparison}
		\end{subfigure}\hspace{-0.05in}
		\begin{subfigure}[t]{0.33\textwidth}
			\centering
			\includegraphics[scale = 0.265]{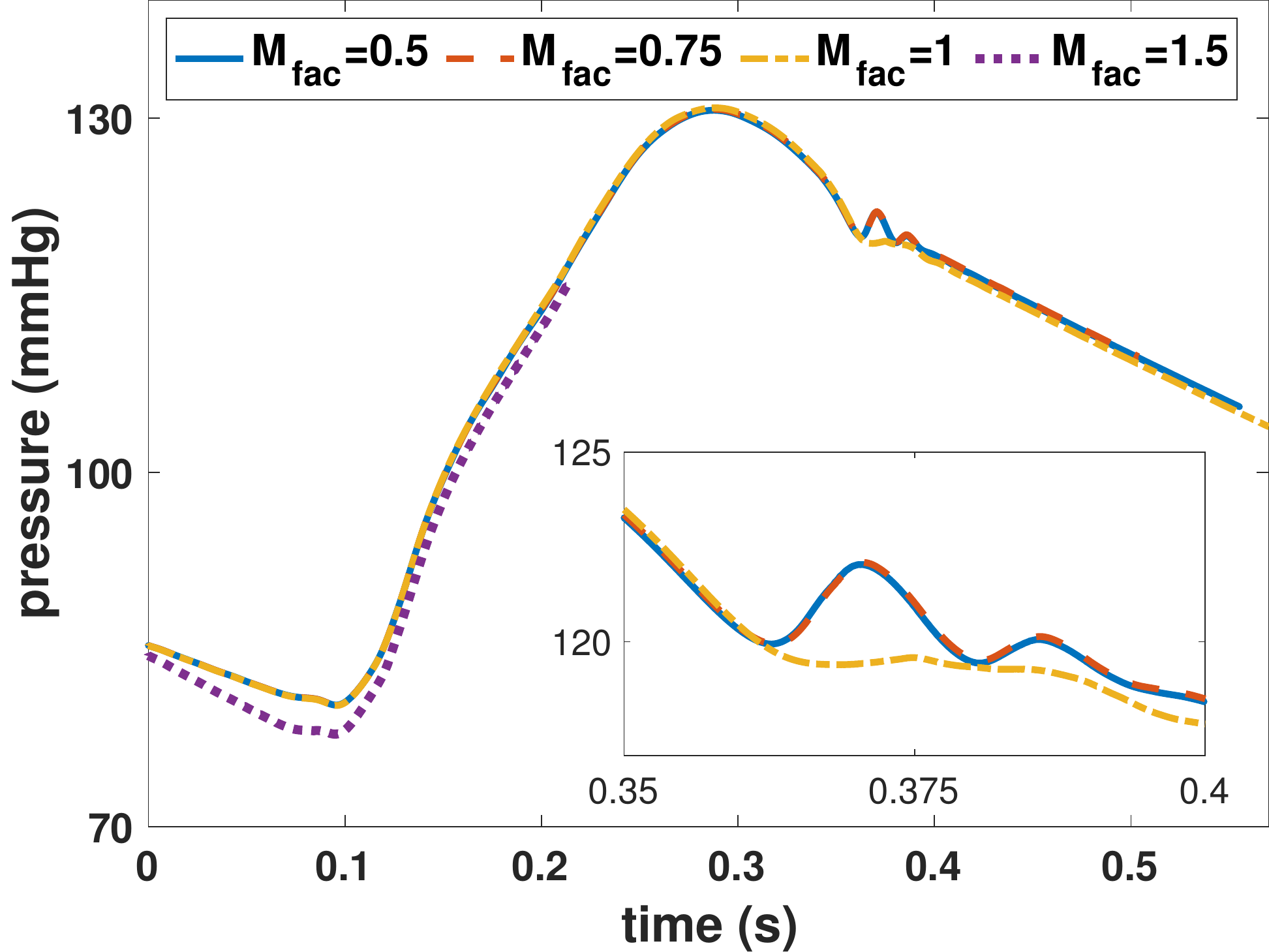} 
			\caption{$P_\text{Ao}$} \label{fig:p_down_comparison}
		\end{subfigure}\hspace{-0.05in}
		\begin{subfigure}[t]{0.33\textwidth}
			\centering
			\includegraphics[scale = 0.265]{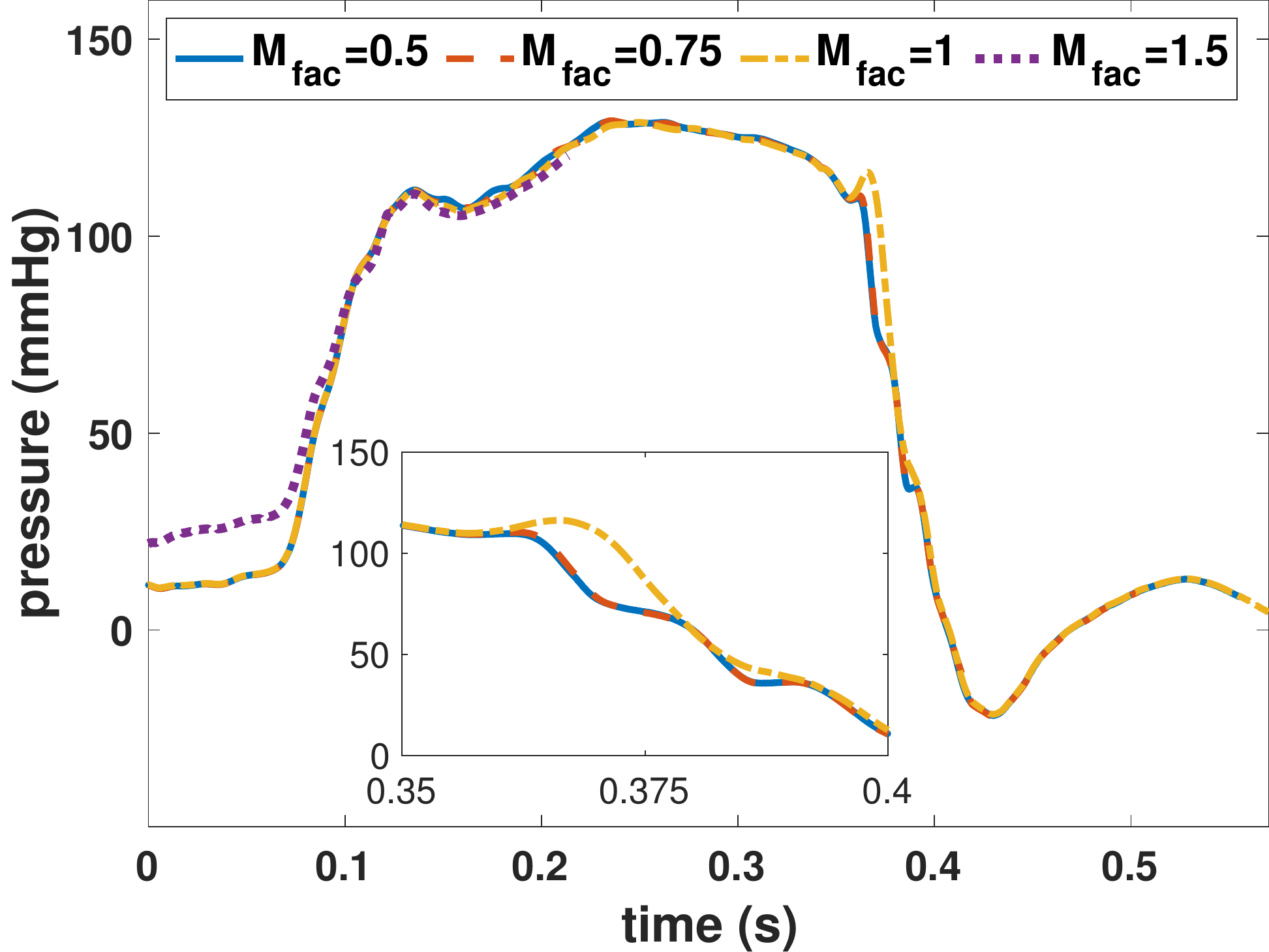} 
			\caption{$P_\text{LV}$} \label{fig:p_up_comparison}
		\end{subfigure}
		\caption{Representative comparison of cross-section views of simulated velocity magnitudes for the bovine pericardial valve models for (a) $\mfac$ = 0.75 and (b) 1.5 during diastole (pressure-loaded when the valve is closed) and systole (shear-dominant flow when the valve is open). The simulations use a three-level locally refined grid with a refinement ratio of two between levels and an $N/2 \times N \times N/2$ coarse grid with $N = 64$, which corresponds to $N = 256$ at the finest level. We also look at comparisons of (c) flow rates ($Q_\text{Ao}$), (d) downstream pressure ($P_\text{Ao}$), and (e) upstream pressure ($P_\text{LV}$) waveforms measured from simulations with $\mfac =$ 0.5, 0.75, 1, 1.5. In panel (b) we see spurious velocities through the structure during diastole as well as local regions with unphysical velocity concentrations during systole (red dashed circles). All comparisons in panels (c)--(e) also indicate that $\mfac =$ 0.5 and 0.75 are in excellent agreement and $\mfac = 1$ shows minor discrepancy during closure as shown in the magnified views. This is because for the case of $\mfac = 1$, not all elements have $\mfac$ exactly equal to 1, but some element have $\mfac > 1$. However, we clearly observe spurious velocities in (b) with $\mfac = 1.5$, which are reflected as negative flow rate measurement as shown in (c). As a result, we also observe discrepancies in both downstream and upstream pressure for $\mfac = 1.5$.}
		\label{fig:BHV}		
	\end{center}
\end{figure}

We first consider the effect of $\mfac$ when using the three-point B-spline kernel for the valve leaflets.
Figures~\ref{fig:bhv_mfac_0_75} and~\ref{fig:bhv_mfac_1_5} compare cross-section views of velocity magnitude for the bovine BHV models for $\mfac = 0.75$ and $1.5$. 
It is clear in Figure~\ref{fig:bhv_mfac_1_5} that there is significant spurious flows through the structure during diastole, but not in Figure~\ref{fig:bhv_mfac_0_75}. 
Figures~\ref{fig:q_comparison},~\ref{fig:p_down_comparison}, and~\ref{fig:p_up_comparison} compare simulated and experimental flow rates and pressure waveforms. 
The spurious velocities that are evident in Figure~\ref{fig:bhv_mfac_1_5} are also reflected in these measurements.
Using $\mfac = 0.75$, the flow rate and pressure data are in excellent agreement with those from the corresponding experiment, whereas we clearly observe the effect of spurious velocities through the structure that are reflected as negative flow rates when the valve is closed and supporting a physiological pressure load (Figure~\ref{fig:q_comparison}). 
We also remark that the simulation using $\mfac = 1.5$ is not able to proceed beyond $t \approx 0.25$ s without significantly reducing the time step size because of the high spurious velocities in the regions highlighted in Figure~\ref{fig:bhv_mfac_1_5}. 
The BHV leaflets in the $\mfac = 1.5$ case also experience unphysical deformations at the free edges of the leaflets during systole.

\begin{figure}[t!!]
	\begin{center}
		\includegraphics[scale = 0.125]{bhv_scale_labeled.pdf}
		\begin{subfigure}[t]{0.25\textwidth}
			\centering
			\includegraphics[scale = 0.275]{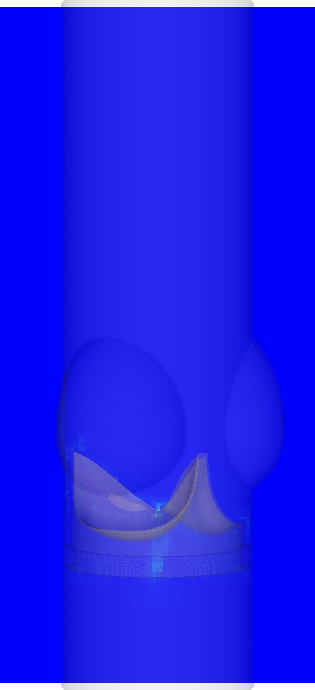} 
			\caption{three-point B-spline} \label{fig:bhv_bsp3}
		\end{subfigure}\hspace{-0.05in}
		\begin{subfigure}[t]{0.25\textwidth}
			\centering
			\includegraphics[scale = 0.275]{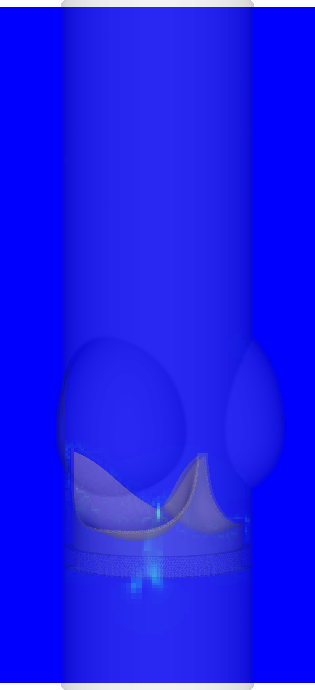} 
			\caption{three-point IB} \label{fig:bhv_ib3}
		\end{subfigure}\hspace{-0.05in}
		\begin{subfigure}[t]{0.25\textwidth}
			\centering
			\includegraphics[scale = 0.275]{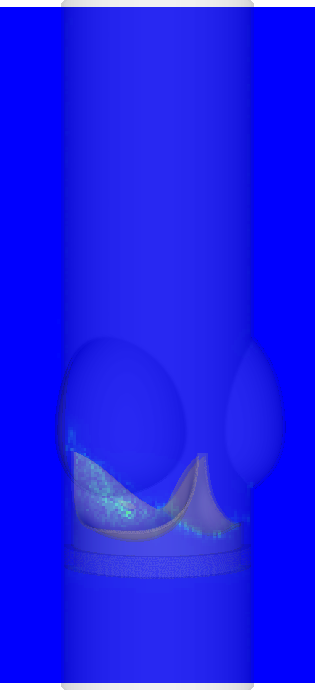} 
			\caption{four-point IB} \label{fig:bhv_ib4}
		\end{subfigure}
		\caption{Representative comparison of cross-section views of simulated velocity magnitudes for the bovine pericardial valve models for (a) the three-point B-spline kernel, (b) the three-point IB kernel, and (c) the four-point IB kernel during diastole. In all of the cases, we use the piecewise linear kernel for the housing and set $\mfac=0.75$. We observe that with the four-point IB kernel, there are relatively large unphysical flows ``through'' the valve leaflets, whereas the only flow we see with the three-point B-spline kernel and the three-point IB kernel are small leakage flows through the middle gap between the valve leaflets, which appear to be physical.}
		\label{fig:BHV_kernel}		
	\end{center}
\end{figure}

We also look at the effect of different kernels under a fixed value of $\mfac$. 
In Figure~\ref{fig:BHV_kernel}, we compare three cases in which we use the three-point B-spline kernel, the three-point IB kernel, and the four-point IB kernel for the valve leaflets, and for all of them use the piecewise linear kernel for the aortic test section and set $\mfac=0.75$. 
We observe immediately during diastole that there is unphysical velocity through the valve leaflets when using the four-point IB kernel (Figure~\ref{fig:bhv_ib4}).

\section{Discussion}
\label{sec:discussion}
This study explores the impacts of various choices of regularized delta functions to approximate the integral transforms that connect the Lagrangian and Eulerian variables, Eqs.~\eqref{eq:fsiconstraint} and~\eqref{eq:noslip}, in the IFED method.
It also investigates the effect of variations in the structural mesh spacing relative to the background Cartesian grid spacing for different kernels on the accuracy using standard FSI benchmark studies.
Our results suggest that kernels satisfying the even--odd condition require higher resolution to achieve similar accuracy as kernels that do not satisfy this condition (e.g., the four-~and six-point IB kernels versus the three-~and five-point IB kernels). 
We also find that, at least for the tests considered herein, narrower kernels are more robust, and that structural meshes that are coarser than the background Cartesian grid can yield improved accuracy compared to structural meshes that are comparable to or finer than the background grid for shear-dominated cases, but not for cases with large normal forces along the fluid-structure interface.
This suggests that to handle both cases within a single model, one needs to use structural meshes with resolutions that are at least as fine as the background grid to avoid instabilities along with a narrower kernel.
The impact of the choice of regularized delta function or relative mesh spacings will likely depend on the many details of the Lagrangian and Eulerian spatial discretizations. 
For instance, different Lagrangian-Eulerian coupling strategies, such as node-based approximations to the integral transforms, may not be suitable for general use with $\mfac > 1$.
A possible explanation for why narrower kernels yield higher accuracies for a given resolution is that the kernels with smaller support result in smaller numerical boundary layers along the fluid-structure interface.
In cases in which the numerical boundary layer is comparable to or larger than the physical boundary layer, changes in numerical boundary layer thicknesses may substantially impact accuracy.
For instance, the channel flow benchmark is an interesting case in which the piecewise linear kernel leads to the best accuracy, as opposed to other tests in which the three-point B-spline kernel yields the best accuracy.
Unlike in other cases considered herein, however, in laminar channel flow, the flow field is completely tangential to the immersed structures, so the dominating factor that affects the accuracy here is indeed the numerical boundary layer.
The piecewise linear kernel leads to the smallest numerical boundary layer effect on the solution, and it provides the best accuracy for this test.
It is evident in Figures~\ref{fig-channel_error} and~\ref{fig-channel_error_vs_mfac} that the error increases as the radius of support of the kernel increases.
More broadly, these results suggest that smoother kernels do not necessarily yield improved accuracy. 
Indeed, based on these results, we speculate that there is a benefit in accuracy to using the minimal amount of smoothing required for a particular model. 
Specifically, kernel functions that smooth out spatial variations in the Lagrangian force when it is spread to the Eulerian grid effectively prevent those variations from impacting the dynamics of the fluid-structure system.
Similarly, kernel functions that smooth out spatial variations in the Eulerian velocity field prevent those variations from influencing the motion of the structure.
This can allow such modes to persist in the computed solution unless otherwise suppressed through physical or numerical smoothing mechanisms.
In the present models, viscous dissipation is the primary physical smoothing mechanism, and some additional smoothing is provided by numerical dissipation from the PPM-type discretization of the convective terms in the momentum equation.
The physical viscosity may have a limited impact on the computed dynamics at moderate-to-high Reynolds numbers at practical spatial resolutions, particularly in three spatial dimensions. 
This suggests that the overall methodology may benefit from additional stabilization that is tailored to the Lagrangian-Eulerian coupling operators, although the construction of such stabilization procedures is beyond the scope of the present study.
As with the higher-order kernel functions, kernels that satisfy the even--odd condition will be oblivious to grid-scale even--odd oscillations in the velocity field when interpolating the velocity to the structure, and instead will only see the mean velocity. 
Specifically, interpolating an alternating $+U/-U$ velocity pattern will yield a structural velocity that is identically 0 for any value of $U$ when using a kernel that satisfies the even--odd condition. 
Although such oscillations will be damped by viscosity, if viscosity is small, these modes may decay slowly.
We believe that this can allows oscillatory modes to persist near or inside the structure, like those that appear in Figure~\ref{fig:BHV_kernel} when the valve is closed for the four-point IB kernel but not for the three-point B-spline or IB kernels.

The results in Sections~\ref{subsec:cylinder},~\ref{subsec:channel}, and~\ref{subsec:turek-hron} indicate that we obtain improved accuracy with a given Cartesian grid resolution for these shear-dominant cases by using relatively coarser Lagrangian nodal spacing ($\mfac > 1$).
This means that the structural mesh in the IFED method can be coarser than that follows the ``rule of thumb'' ($\mfac = 0.5$) for the nodally interacting IB method by a factor of 8 ($\mfac = 4$), which results in a significant improvement in both accuracy and efficiency.
However, in the pressure-loaded case considered in Section~\ref{subsec:pressurized_band}, we observe that the Lagrangian mesh needs to have a resolution that is similar to or relatively finer than the Cartesian grid ($\mfac \leq 1$) to avoid spurious velocities through the structure.
In fact, it is common in simulations using complex geometries to have many mesh elements that are comparable to or finer than the background Cartesian grid to preserve fine-scale geometric features.
Understanding this transition in accuracy between shear-dominant and pressure-loaded cases is another possible future area of research.

The benchmarks suggest that the three-point B-spline kernel is the best overall choice considering both shear- and pressure-dominant flows because it is less sensitive to the relative structural mesh spacing.
We emphasize, however, that under sufficiently fine grid resolution, different kernels all appear ultimately to converge to the same results.
However, this study also suggests that optimal choices of numerical and discretization parameters can provide consistent solutions at the coarser grid resolutions that are needed to facilitate the deployment of the methodology to large-scale three-dimensional models.
Using these results, we also applied our findings from benchmark studies to an FSI model of bovine pericardial BHV in a pulse duplicator, which involves both pressure-loaded and shear-dominant flows in a rigid and stationary channel with immersed elastic structures inside. 
Results obtained using this large-scale model are consistent with the key findings of benchmark test cases, and we obtain accurate results only for $\mfac < 1$.
For the case in which $\mfac = 1$, the results are in excellent agreement with results using $\mfac = 0.5$ and $0.75$, except for a slight discrepancy during closure as shown in Figures~\ref{fig:q_comparison}--\ref{fig:p_up_comparison} because some elements have $\mfac > 1$.
A limitation of this study is that not all possible kernel function constructions are considered.
Another limitation is that it considers specific Lagrangian and Eulerian spatial discretizations. 
Although this study is done within the context of the IFED method, the effect of different kernels could be important not just for this method, but more generally for other IB-type methods that use regularized delta functions to mediate fluid-structure interaction.

%\clearpage
\section*{Acknowledgments}
B.E.G. acknowledges funding from the NIH (Awards R01HL117063, U01HL143336, and R01HL157631) and NSF (Awards OAC 1450327, OAC 1652541, OAC 1931516, and CBET 1757193).
J.H.L. acknowledges funding from the Integrative Vascular Biology Training Program (NIH Award 5T32HL069768-17) at the University of North Carolina School of Medicine.
Simulations were performed using facilities provided by the University of North Carolina at Chapel Hill through the Research Computing Division of UNC Information Technology Services.
We thank Aaron Barrett, Ebrahim M. Kolahdouz, Ben Vadala-Roth, David Wells, and the reviewers for their constructive comments that improved the manuscript.

\section*{Conflict of Interest} 
No conflict of interest.

%%%%%%%%%%%%%%%%%%%%%%%%%%%%%%%%%%%%%%%%%%%%%%%%%%%%%%%%%%%%%%%%%%%%%%%%%%%%%%%%%%%%%%%%%%%%%%%%%%%%%%%%%%%%%%%%%%%%%%%%%%%

\bibliography{reference}

\clearpage
\renewcommand\thefigure{S\arabic{figure}}
\renewcommand\thetable{S\arabic{table}}
\renewcommand\theequation{S\arabic{equation}}
\setcounter{equation}{0}
\setcounter{figure}{0}
\setcounter{table}{0}
\appendix
\section*{Appendix}
\section{Turek-Hron Benchmark Results for Various Choices of Kernel Function} 
\label{sec:TH_appendix}
\begin{figure}[t!!]
\centering
\includegraphics[scale = 0.365,trim={10 30 0 10},clip]{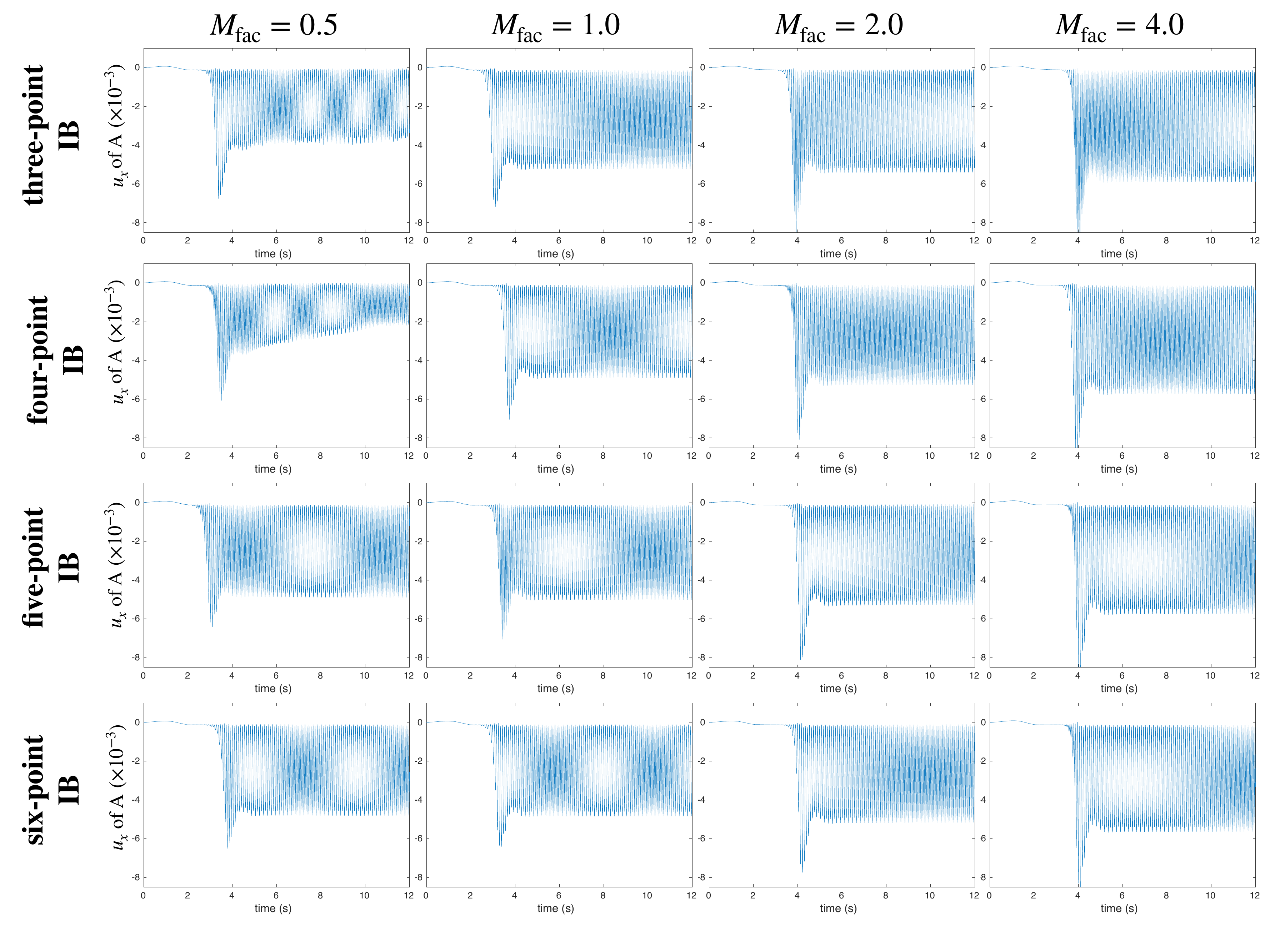}
\includegraphics[scale = 0.365,trim={445 30 340 10},clip]{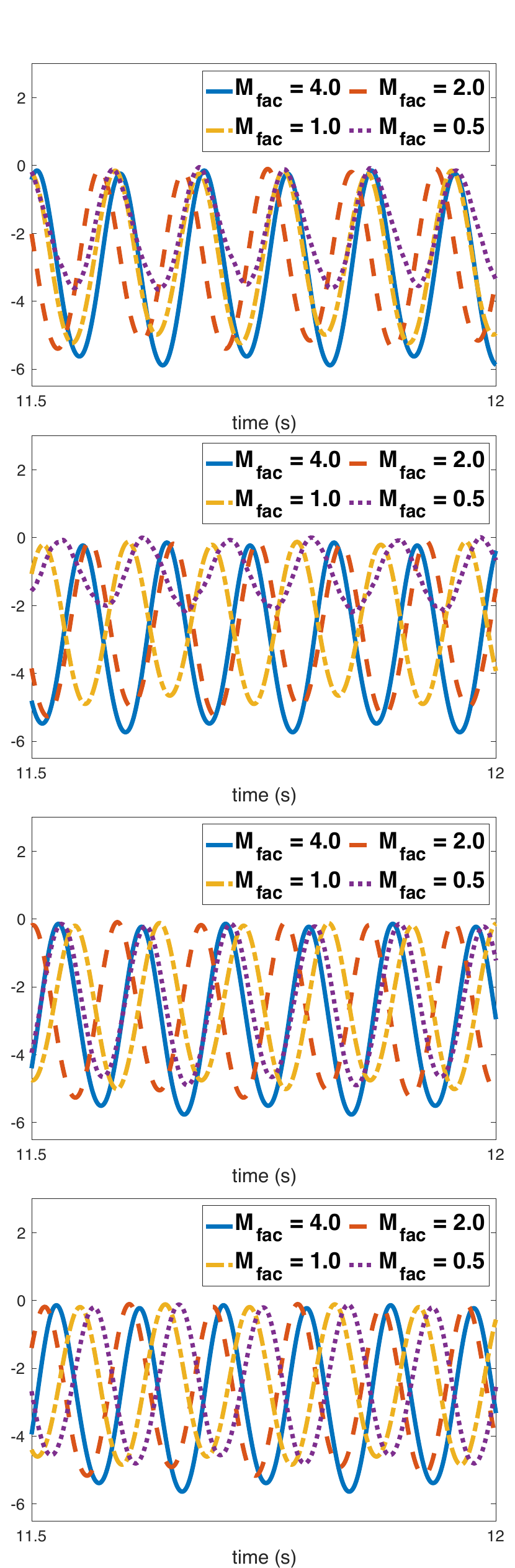}
\includegraphics[scale = 0.365,trim={10 30 0 35},clip]{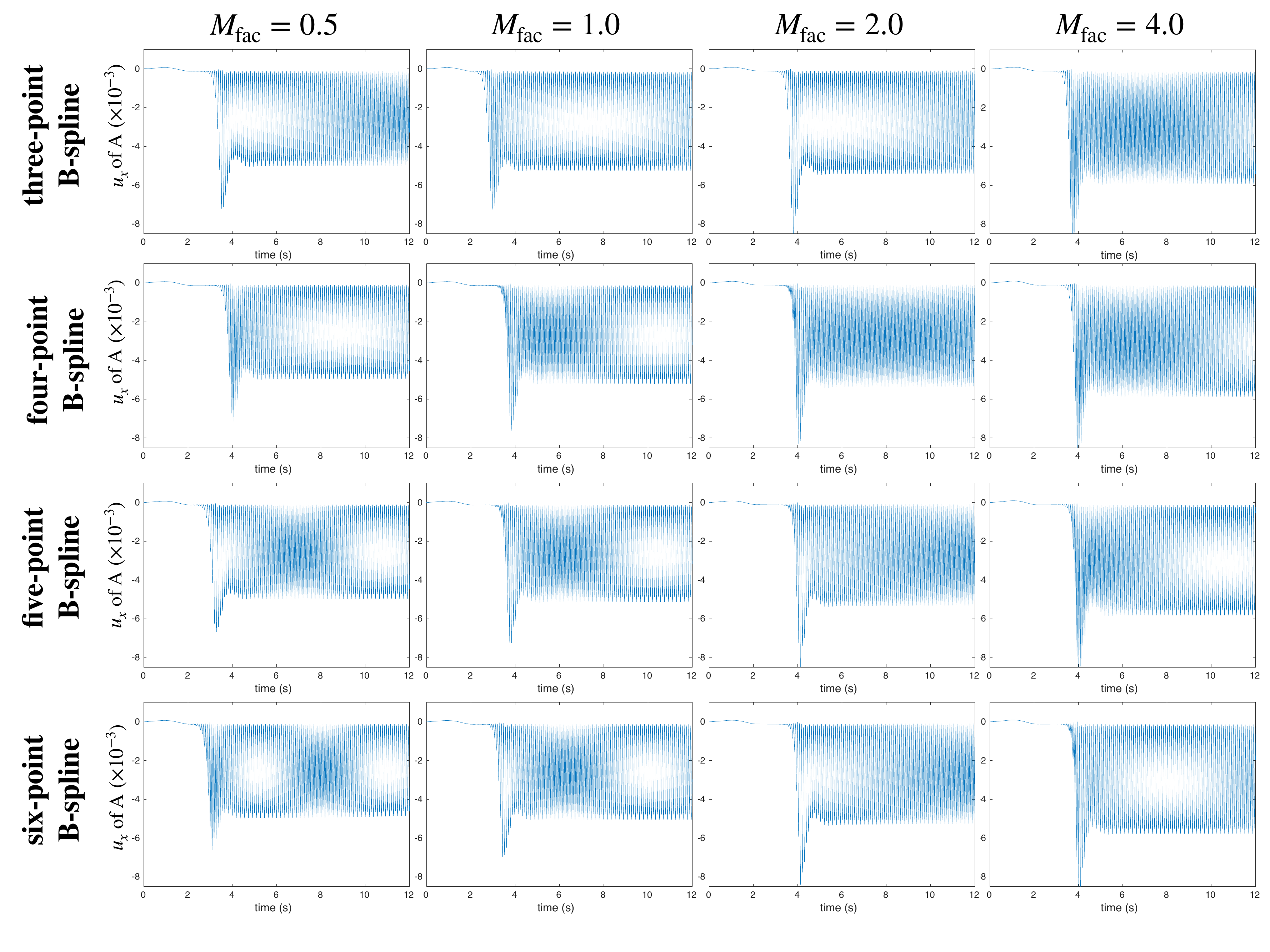} 
\includegraphics[scale = 0.365,trim={445 30 340 35},clip]{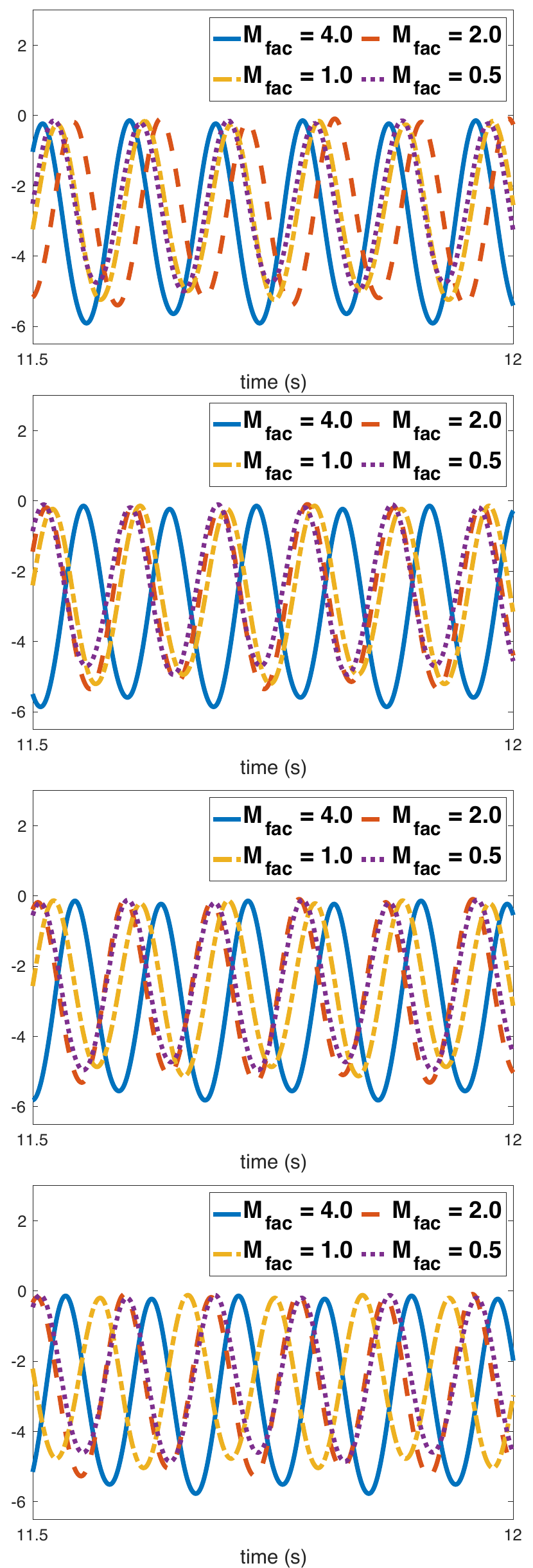}
\hspace{-0.09in}
\caption{$x$-displacement ($A_x$) of the point $A$ for different values of $\mfac$ for the modified Turek-Hron benchmark using different IB and B-spline kernels at a Cartesian resolution of $N = 64$. Figures in the rightmost panels show the periodic oscillations between $t=11.5$ and $t=12$.}
\label{fig:A_x_comparison_N=64_appendix}
\end{figure}

\begin{figure}[t!!]
\centering
\includegraphics[scale = 0.365,trim={10 30 0 10},clip]{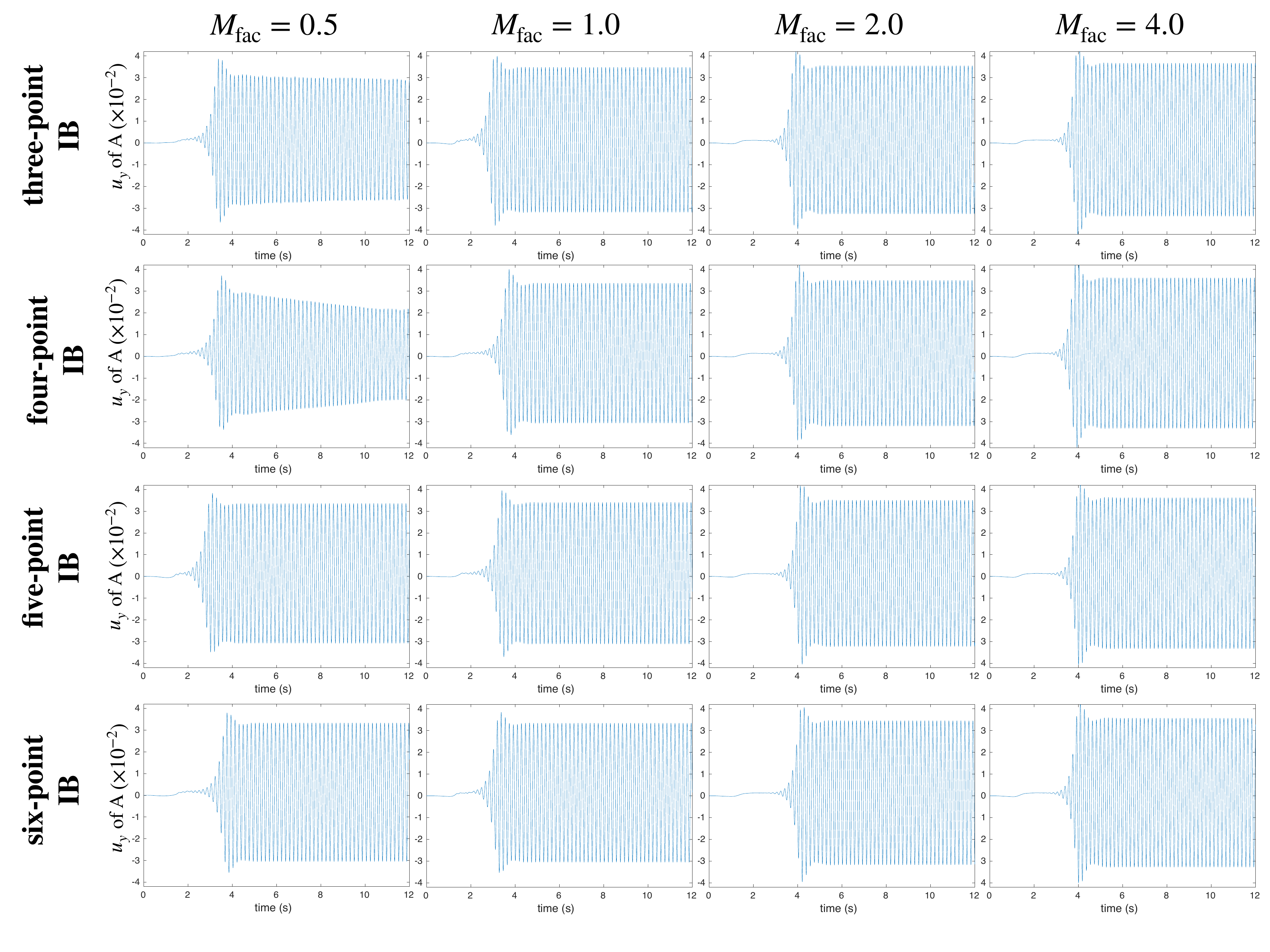}
\includegraphics[scale = 0.365,trim={445 30 340 10},clip]{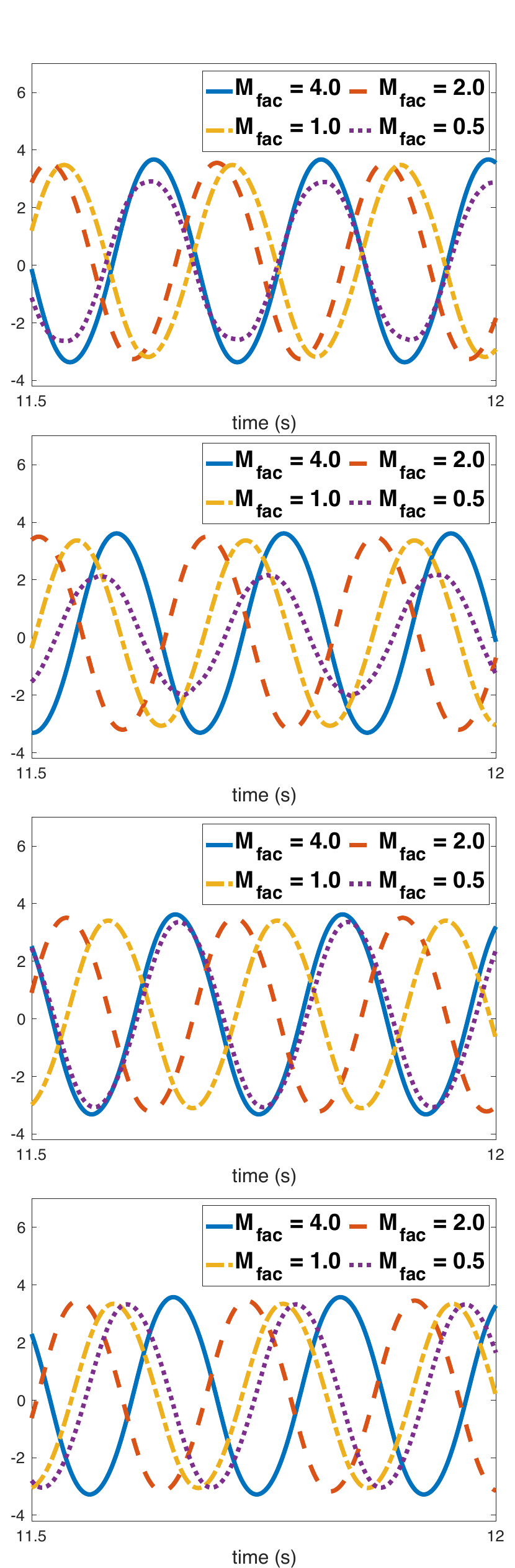}
\includegraphics[scale = 0.365,trim={10 30 0 35},clip]{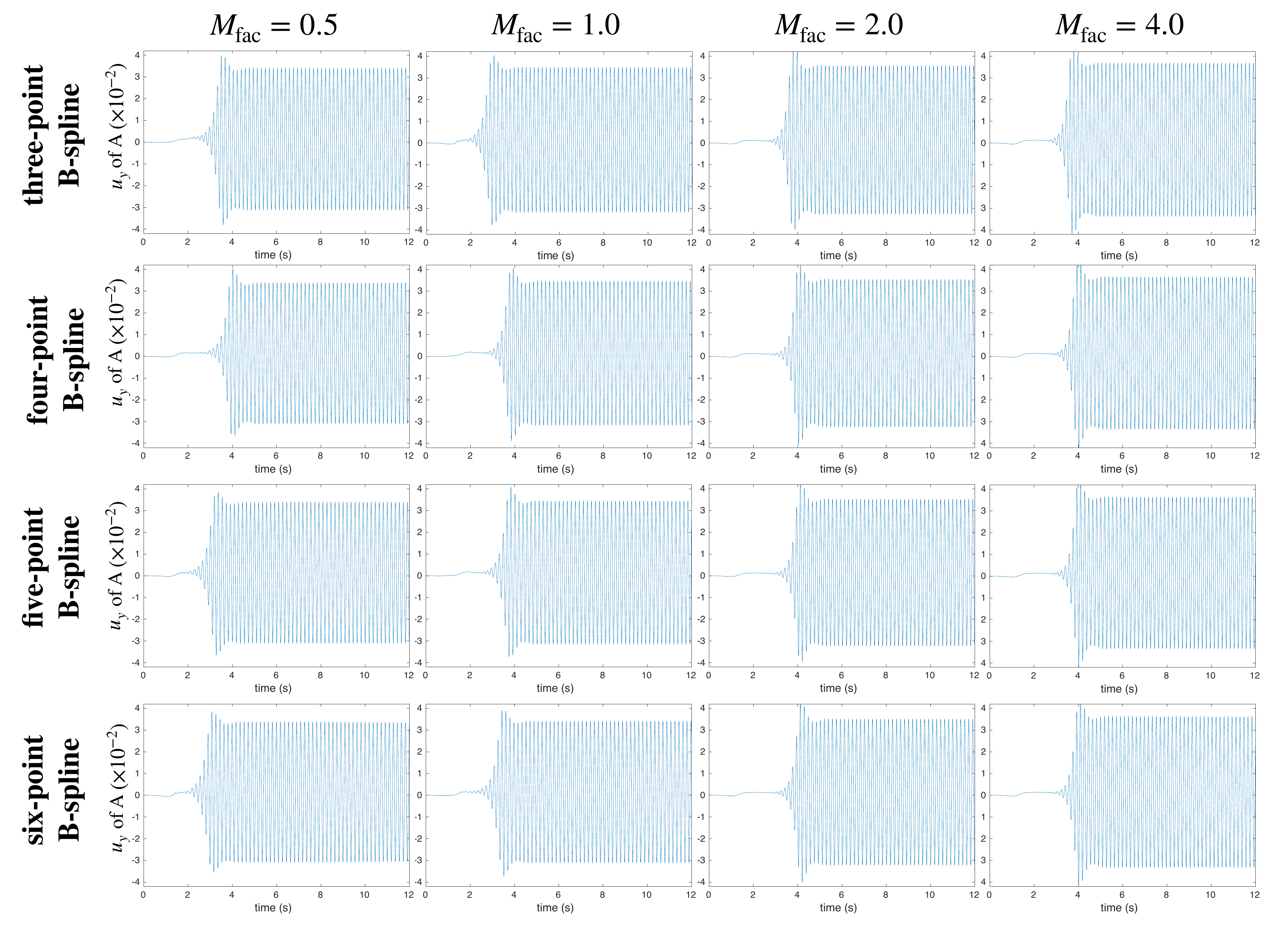} 
\includegraphics[scale = 0.365,trim={445 30 340 35},clip]{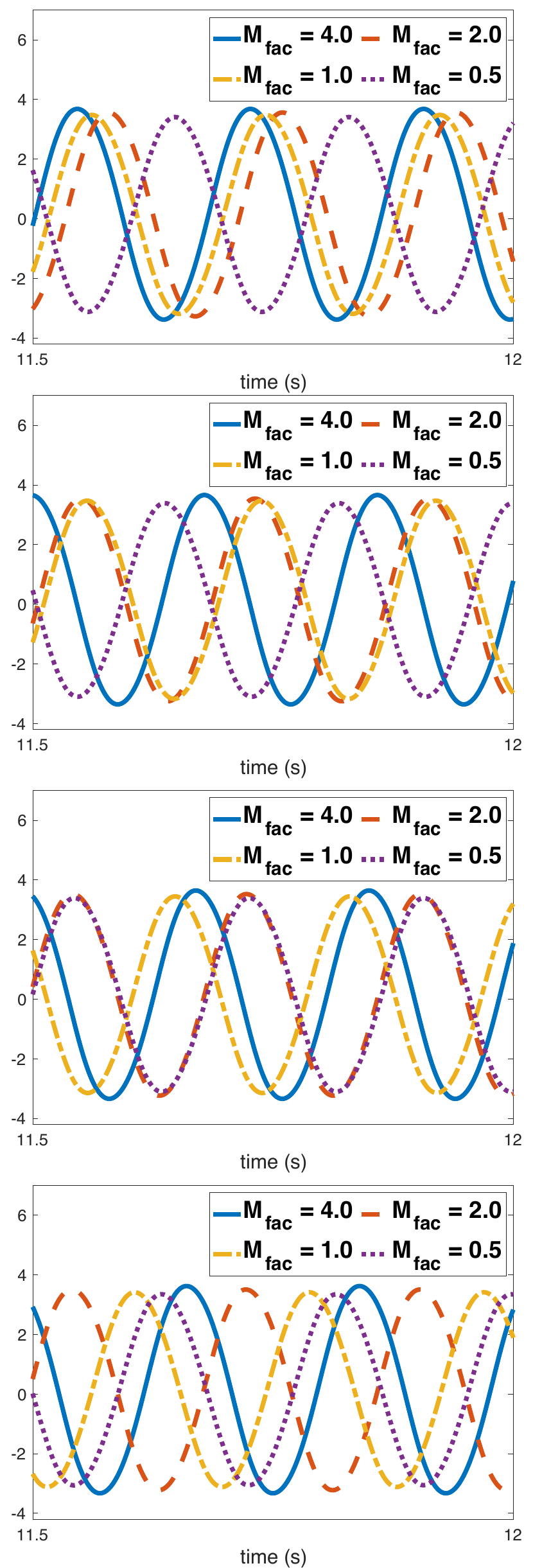}
\hspace{-0.09in}
\caption{$y$-displacement ($A_y$) of the point $A$ for different values of $\mfac$ for the modified Turek-Hron benchmark using different IB and B-spline kernels at a Cartesian resolution of $N = 64$. Figures in the rightmost panels show the periodic oscillations between $t=11.5$ and $t=12$.}
\label{fig:A_y_comparison_N=64_appendix}
\end{figure}

This section details results for the modified Turek-Hron benchmark using various IB and B-spline kernels (see Figures~\ref{fig:A_x_comparison_N=64_appendix} and~\ref{fig:A_y_comparison_N=64_appendix}).
Table~\ref{table:nh_kinematics_N=64} reports the means for $A_x$ and $A_y$, which are $x$-, $y$-displacements of the point $A$, as well as the Strouhal numbers corresponding to the oscillations of $A_x$ and $A_y$ at periodic steady-state.
Here we look at relatively coarser resolution cases, in which the number of grid cells on coarsest grid level is $N = 64$. 
These results indicate that the three-point B-spline kernel is the only kernel that shows consistent Strouhal numbers for all values of $\mfac$ = 0.5, 1, 2, and 4 at $N=64$, and that it is less sensitive to changes in $\mfac$.
Other kernels show loss of accuracy as we refine the Lagrangian mesh for a fixed Eulerian grid that is relatively coarse.

\begin{table}[t!!]
    \setlength\tabcolsep{4.5pt}
    \scriptsize
	\centering	
	\caption{Results for the modified Turek-Hron benchmark with different IB and B-spline kernels and various values of $\mfac$ at a Cartesian resolution of $N = 64$.}
\begin{tabular}{c c c |c c | c c| c c|l}
\cline{2-9}
& \multicolumn{2}{ |c| }{$\mfac = 0.5$} & \multicolumn{2}{ |c| }{$\mfac = 1.0$} & \multicolumn{2}{ |c| }{$\mfac = 2.0$} & \multicolumn{2}{ |c| }{$\mfac = 4.0$} \\  
\cline{1-9}
\multicolumn{1}{ |c| }{Kernel} & $A_x$ $(\times 10^{-3})$ & $St_x$ & $A_x$ $(\times 10^{-3})$ & $St_x$ & $A_x$ $(\times 10^{-3})$ & $St_x$ & $A_x$ $(\times 10^{-3})$ & $St_x$ & \\ \cline{1-9}
\multicolumn{1}{ |c| }{IB (3-point)} & $-2.03 \pm 1.97$ & 10.4 & $-2.69 \pm 2.54$ & 10.4 & $-2.76 \pm 2.65$ & 10.8 & $-3.02 \pm 2.87$ & 10.8 &    \\ 
\cline{1-9}
\multicolumn{1}{ |c| }{IB (4-point)} & $-1.55 \pm 1.56$ & 10.4 & $-2.51 \pm 2.39$ & 10.4 & $-2.69 \pm 2.58$ & 10.8 & $-2.94 \pm 2.80$ & 10.8 &    \\ 
\cline{1-9}
\multicolumn{1}{ |c| }{IB (5-point)} & $-2.51 \pm 2.38$ & 10.4 & $-2.57 \pm 2.44$ & 10.4 & $-2.70 \pm 2.60$ & 10.8 & $-2.95 \pm 2.81$ & 10.8 &    \\ 
\cline{1-9}
\multicolumn{1}{ |c| }{IB (6-point)} & $-2.46 \pm 2.35$ & 10.4 & $-2.48 \pm 2.36$ & 10.4 & $-2.65 \pm 2.53$ & 10.4 & $-2.89 \pm 2.75$ & 10.8 &    \\ 
\cline{1-9}
\multicolumn{1}{ |c| }{B-spline (3-point)} & $-2.57 \pm 2.44$ & 10.8 & $-2.69 \pm 2.55$ & 10.8 & $-2.76 \pm 2.66$ & 10.8 & $-3.03 \pm 2.89$ & 10.8 &    \\ 
\cline{1-9}
\multicolumn{1}{ |c| }{B-spline (4-point)} & $-2.53 \pm 2.42$ & 10.4 & $-2.67 \pm 2.54$ & 10.4 & $-2.74 \pm 2.63$ & 10.8 & $-3.00 \pm 2.86$ & 10.8 &    \\ 
\cline{1-9}
\multicolumn{1}{ |c| }{B-spline (5-point)} & $-2.55 \pm 2.43$ & 10.4 & $-2.63 \pm 2.50$ & 10.8 & $-2.72 \pm 2.62$ & 10.8 & $-2.98 \pm 2.84$ & 10.8 &    \\ 
\cline{1-9}
\multicolumn{1}{ |c| }{B-spline (6-point)} & $-2.54 \pm 2.42$ & 10.4 & $-2.58 \pm 2.46$ & 10.4 & $-2.70 \pm 2.60$ & 10.8 & $-2.96 \pm 2.82$ & 10.8 &    \\ 
\hhline{=========}
\multicolumn{1}{ |c| }{Kernel} & $A_y$ $(\times 10^{-3})$ & $St_y$ & $A_y$ $(\times 10^{-3})$ & $St_y$ & $A_y$ $(\times 10^{-3})$ & $St_y$ & $A_y$ $(\times 10^{-3})$ & $St_y$ &\\
\cline{1-9}
\multicolumn{1}{ |c| }{IB (3-point)} & $1.37 \pm 29.2$ & 5.00 & $1.45 \pm 33.3$ & 5.00 & $1.42 \pm 34.1$ & 5.00 & $1.48 \pm 35.2$ & 5.00 &    \\ 
\cline{1-9}
\multicolumn{1}{ |c| }{IB (4-point)} & $1.04 \pm 25.9$ & 5.00 & $1.46 \pm 32.2$ & 5.00 & $1.41 \pm 33.5$ & 5.00 & $1.47 \pm 34.6$ & 5.00 &    \\ 
\cline{1-9}
\multicolumn{1}{ |c| }{IB (5-point)} & $1.39 \pm 32.2$ & 5.00 & $1.48 \pm 32.6$ & 5.00 & $1.41 \pm 33.6$ & 5.00 & $1.48 \pm 34.7$ & 5.00 &    \\ 
\cline{1-9}
\multicolumn{1}{ |c| }{IB (6-point)} & $1.39 \pm 31.9$ & 5.00 & $1.44 \pm 32.0$ & 5.00 & $1.41 \pm 33.2$ & 5.00 & $1.47 \pm 34.3$ & 5.00 &    \\ 
\cline{1-9}
\multicolumn{1}{ |c| }{B-spline (3-point)} & $1.41 \pm 32.7$ & 5.00 & $1.44 \pm 33.4$ & 5.00 & $1.41 \pm 34.2$ & 5.00 & $1.49 \pm 35.3$ & 5.00 &    \\ 
\cline{1-9}
\multicolumn{1}{ |c| }{B-spline (4-point)} & $1.41 \pm 32.5$ & 5.00 & $1.46 \pm 33.3$ & 5.00 & $1.42 \pm 34.0$ & 5.00 & $1.48 \pm 35.1$ & 5.00 &    \\ 
\cline{1-9}
\multicolumn{1}{ |c| }{B-spline (5-point)} & $1.39 \pm 32.5$ & 5.00 & $1.48 \pm 33.0$ & 5.00 & $1.42 \pm 33.8$ & 5.00 & $1.48 \pm 34.9$ & 5.00 &    \\ 
\cline{1-9}
\multicolumn{1}{ |c| }{B-spline (6-point)} & $1.41 \pm 32.4$ & 5.00 & $1.48 \pm 32.6$ & 5.00 & $1.41 \pm 33.7$ & 5.00 & $1.48 \pm 34.8$ & 5.00 &    \\ 
\cline{1-9}
\end{tabular}
	\label{table:nh_kinematics_N=64}	
\end{table} 

\end{document}